\pdfoutput=1
\documentclass[a4paper,10pt]{amsart}

\usepackage[utf8]{inputenc}

\usepackage[expansion=false]{microtype}
\usepackage{amssymb,amsmath,amsopn,amsxtra,amsthm,amsfonts}
\usepackage{xr}
\usepackage[mathcal]{euscript}
\usepackage{mathalfa}

\usepackage[english, status=final, nomargin, inline]{fixme}
\fxusetheme{color}

\usepackage[pdfusetitle,unicode]{hyperref}

\makeatletter
\let\c@equation\c@subsubsection

\makeatother

\newtheorem{cor}[subsubsection]{Corollary}
\newtheorem{lem}[subsubsection]{Lemma}
\newtheorem{prop}[subsubsection]{Proposition}

\newtheorem{thm}[subsubsection]{Theorem}

\newtheorem{defn}[subsubsection]{Definition}

\theoremstyle{definition}

\theoremstyle{remark}
\newtheorem{rem}[subsubsection]{Remark}

\newcommand{\secref}[1]{Section~\ref{#1}}
\newcommand{\ssecref}[1]{\sectsign\ref{#1}}
\newcommand{\sssecref}[1]{\ref{#1}}

\renewcommand{\eqref}[1]{(\ref{#1})}

\usepackage{tikz}
\usetikzlibrary{matrix}
\usepackage{tikz-cd}

\usepackage{fixltx2e}

\usepackage{xspace}

\usepackage{xifthen}

\usepackage{xparse}

\newcommand{\nc}{\newcommand}
\nc{\renc}{\renewcommand}
\nc{\ssec}{\subsection}
\nc{\sssec}{\subsubsection}
\nc{\on}{\operatorname}
\nc{\term}[1]{#1\xspace}

\DeclareMathSymbol{A}{\mathalpha}{operators}{`A}
\DeclareMathSymbol{B}{\mathalpha}{operators}{`B}
\DeclareMathSymbol{C}{\mathalpha}{operators}{`C}
\DeclareMathSymbol{D}{\mathalpha}{operators}{`D}
\DeclareMathSymbol{E}{\mathalpha}{operators}{`E}
\DeclareMathSymbol{F}{\mathalpha}{operators}{`F}
\DeclareMathSymbol{G}{\mathalpha}{operators}{`G}
\DeclareMathSymbol{H}{\mathalpha}{operators}{`H}
\DeclareMathSymbol{I}{\mathalpha}{operators}{`I}
\DeclareMathSymbol{J}{\mathalpha}{operators}{`J}
\DeclareMathSymbol{K}{\mathalpha}{operators}{`K}
\DeclareMathSymbol{L}{\mathalpha}{operators}{`L}
\DeclareMathSymbol{M}{\mathalpha}{operators}{`M}
\DeclareMathSymbol{N}{\mathalpha}{operators}{`N}
\DeclareMathSymbol{O}{\mathalpha}{operators}{`O}
\DeclareMathSymbol{P}{\mathalpha}{operators}{`P}
\DeclareMathSymbol{Q}{\mathalpha}{operators}{`Q}
\DeclareMathSymbol{R}{\mathalpha}{operators}{`R}
\DeclareMathSymbol{S}{\mathalpha}{operators}{`S}
\DeclareMathSymbol{T}{\mathalpha}{operators}{`T}
\DeclareMathSymbol{U}{\mathalpha}{operators}{`U}
\DeclareMathSymbol{V}{\mathalpha}{operators}{`V}
\DeclareMathSymbol{W}{\mathalpha}{operators}{`W}
\DeclareMathSymbol{X}{\mathalpha}{operators}{`X}
\DeclareMathSymbol{Y}{\mathalpha}{operators}{`Y}
\DeclareMathSymbol{Z}{\mathalpha}{operators}{`Z}

\nc{\sA}{\ensuremath{\mathcal{A}}\xspace}
\nc{\sB}{\ensuremath{\mathcal{B}}\xspace}
\nc{\sC}{\ensuremath{\mathcal{C}}\xspace}
\nc{\sD}{\ensuremath{\mathcal{D}}\xspace}
\nc{\sE}{\ensuremath{\mathcal{E}}\xspace}
\nc{\sF}{\ensuremath{\mathcal{F}}\xspace}
\nc{\sG}{\ensuremath{\mathcal{G}}\xspace}
\nc{\sH}{\ensuremath{\mathcal{H}}\xspace}
\nc{\sI}{\ensuremath{\mathcal{I}}\xspace}
\nc{\sJ}{\ensuremath{\mathcal{J}}\xspace}
\nc{\sK}{\ensuremath{\mathcal{K}}\xspace}
\nc{\sL}{\ensuremath{\mathcal{L}}\xspace}
\nc{\sM}{\ensuremath{\mathcal{M}}\xspace}
\nc{\sN}{\ensuremath{\mathcal{N}}\xspace}
\nc{\sO}{\ensuremath{\mathcal{O}}\xspace}
\nc{\sP}{\ensuremath{\mathcal{P}}\xspace}
\nc{\sQ}{\ensuremath{\mathcal{Q}}\xspace}
\nc{\sR}{\ensuremath{\mathcal{R}}\xspace}
\nc{\sS}{\ensuremath{\mathcal{S}}\xspace}
\nc{\sT}{\ensuremath{\mathcal{T}}\xspace}
\nc{\sU}{\ensuremath{\mathcal{U}}\xspace}
\nc{\sV}{\ensuremath{\mathcal{V}}\xspace}
\nc{\sW}{\ensuremath{\mathcal{W}}\xspace}
\nc{\sX}{\ensuremath{\mathcal{X}}\xspace}
\nc{\sY}{\ensuremath{\mathcal{Y}}\xspace}
\nc{\sZ}{\ensuremath{\mathcal{Z}}\xspace}

\nc{\bA}{\ensuremath{\mathbf{A}}\xspace}
\nc{\bB}{\ensuremath{\mathbf{B}}\xspace}
\nc{\bC}{\ensuremath{\mathbf{C}}\xspace}
\nc{\bD}{\ensuremath{\mathbf{D}}\xspace}
\nc{\bE}{\ensuremath{\mathbf{E}}\xspace}
\nc{\bF}{\ensuremath{\mathbf{F}}\xspace}
\nc{\bG}{\ensuremath{\mathbf{G}}\xspace}
\nc{\bH}{\ensuremath{\mathbf{H}}\xspace}
\nc{\bI}{\ensuremath{\mathbf{I}}\xspace}
\nc{\bJ}{\ensuremath{\mathbf{J}}\xspace}
\nc{\bK}{\ensuremath{\mathbf{K}}\xspace}
\nc{\bL}{\ensuremath{\mathbf{L}}\xspace}
\nc{\bM}{\ensuremath{\mathbf{M}}\xspace}
\nc{\bN}{\ensuremath{\mathbf{N}}\xspace}
\nc{\bO}{\ensuremath{\mathbf{O}}\xspace}
\nc{\bP}{\ensuremath{\mathbf{P}}\xspace}
\nc{\bQ}{\ensuremath{\mathbf{Q}}\xspace}
\nc{\bR}{\ensuremath{\mathbf{R}}\xspace}
\nc{\bS}{\ensuremath{\mathbf{S}}\xspace}
\nc{\bT}{\ensuremath{\mathbf{T}}\xspace}
\nc{\bU}{\ensuremath{\mathbf{U}}\xspace}
\nc{\bV}{\ensuremath{\mathbf{V}}\xspace}
\nc{\bW}{\ensuremath{\mathbf{W}}\xspace}
\nc{\bX}{\ensuremath{\mathbf{X}}\xspace}
\nc{\bY}{\ensuremath{\mathbf{Y}}\xspace}
\nc{\bZ}{\ensuremath{\mathbf{Z}}\xspace}

\nc{\dA}{\ensuremath{\mathds{A}}\xspace}
\nc{\dB}{\ensuremath{\mathds{B}}\xspace}
\nc{\dC}{\ensuremath{\mathds{C}}\xspace}
\nc{\dD}{\ensuremath{\mathds{D}}\xspace}
\nc{\dE}{\ensuremath{\mathds{E}}\xspace}
\nc{\dF}{\ensuremath{\mathds{F}}\xspace}
\nc{\dG}{\ensuremath{\mathds{G}}\xspace}
\nc{\dH}{\ensuremath{\mathds{H}}\xspace}
\nc{\dI}{\ensuremath{\mathds{I}}\xspace}
\nc{\dJ}{\ensuremath{\mathds{J}}\xspace}
\nc{\dK}{\ensuremath{\mathds{K}}\xspace}
\nc{\dL}{\ensuremath{\mathds{L}}\xspace}
\nc{\dM}{\ensuremath{\mathds{M}}\xspace}
\nc{\dN}{\ensuremath{\mathds{N}}\xspace}
\nc{\dO}{\ensuremath{\mathds{O}}\xspace}
\nc{\dP}{\ensuremath{\mathds{P}}\xspace}
\nc{\dQ}{\ensuremath{\mathds{Q}}\xspace}
\nc{\dR}{\ensuremath{\mathds{R}}\xspace}
\nc{\dS}{\ensuremath{\mathds{S}}\xspace}
\nc{\dT}{\ensuremath{\mathds{T}}\xspace}
\nc{\dU}{\ensuremath{\mathds{U}}\xspace}
\nc{\dV}{\ensuremath{\mathds{V}}\xspace}
\nc{\dW}{\ensuremath{\mathds{W}}\xspace}
\nc{\dX}{\ensuremath{\mathds{X}}\xspace}
\nc{\dY}{\ensuremath{\mathds{Y}}\xspace}
\nc{\dZ}{\ensuremath{\mathds{Z}}\xspace}

\nc{\bbA}{\ensuremath{\mathbb{A}}\xspace}
\nc{\bbB}{\ensuremath{\mathbb{B}}\xspace}
\nc{\bbC}{\ensuremath{\mathbb{C}}\xspace}
\nc{\bbD}{\ensuremath{\mathbb{D}}\xspace}
\nc{\bbE}{\ensuremath{\mathbb{E}}\xspace}
\nc{\bbF}{\ensuremath{\mathbb{F}}\xspace}
\nc{\bbG}{\ensuremath{\mathbb{G}}\xspace}
\nc{\bbH}{\ensuremath{\mathbb{H}}\xspace}
\nc{\bbI}{\ensuremath{\mathbb{I}}\xspace}
\nc{\bbJ}{\ensuremath{\mathbb{J}}\xspace}
\nc{\bbK}{\ensuremath{\mathbb{K}}\xspace}
\nc{\bbL}{\ensuremath{\mathbb{L}}\xspace}
\nc{\bbM}{\ensuremath{\mathbb{M}}\xspace}
\nc{\bbN}{\ensuremath{\mathbb{N}}\xspace}
\nc{\bbO}{\ensuremath{\mathbb{O}}\xspace}
\nc{\bbP}{\ensuremath{\mathbb{P}}\xspace}
\nc{\bbQ}{\ensuremath{\mathbb{Q}}\xspace}
\nc{\bbR}{\ensuremath{\mathbb{R}}\xspace}
\nc{\bbS}{\ensuremath{\mathbb{S}}\xspace}
\nc{\bbT}{\ensuremath{\mathbb{T}}\xspace}
\nc{\bbU}{\ensuremath{\mathbb{U}}\xspace}
\nc{\bbV}{\ensuremath{\mathbb{V}}\xspace}
\nc{\bbW}{\ensuremath{\mathbb{W}}\xspace}
\nc{\bbX}{\ensuremath{\mathbb{X}}\xspace}
\nc{\bbY}{\ensuremath{\mathbb{Y}}\xspace}
\nc{\bbZ}{\ensuremath{\mathbb{Z}}\xspace}

\nc{\mrm}[1]{\ensuremath{\mathrm{#1}}\xspace}
\nc{\mbf}[1]{\ensuremath{\mathbf{#1}}\xspace}
\nc{\mcal}[1]{\ensuremath{\mathcal{#1}}\xspace}
\nc{\msc}[1]{\ensuremath{\mathscr{#1}}\xspace}

\renc{\bar}[1]{\overline{#1}}

\let\sectsign\S
\let\S\relax

\nc{\sub}{\subset}
\nc{\too}{\longrightarrow}
\nc{\hook}{\hookrightarrow}
\nc*{\hooklongrightarrow}{\ensuremath{\lhook\joinrel\relbar\joinrel\rightarrow}}
\nc{\hooklong}{\hooklongrightarrow}
\nc{\twoheadlongrightarrow}{\relbar\joinrel\twoheadrightarrow}
\nc{\shiso}{\approx}
\nc{\isoto}{\xrightarrow{\sim}}
\nc{\isofrom}{\xleftarrow{\sim}}
\renc{\ge}{\geqslant}
\renc{\le}{\leqslant}
\renc{\geq}{\geqslant}
\renc{\leq}{\leqslant}

\nc{\id}{\mathrm{id}}

\DeclareMathOperator{\rk}{\mathrm{rk}}
\DeclareMathOperator{\Hom}{\mathrm{Hom}}
\nc{\uHom}{\underline{\smash{\Hom}}}
\DeclareMathOperator{\Maps}{\mathrm{Maps}}
\DeclareMathOperator{\Aut}{\mathrm{Aut}}
\DeclareMathOperator{\End}{\mathrm{End}}

\nc{\Pre}{\mathrm{PSh}{}}
\nc{\uEnd}{\underline{\smash{\End}}}

\renc{\lim}{\operatorname*{lim}}
\nc{\colim}{\operatorname*{colim}}
\nc{\Cofib}{\on{Cofib}}
\nc{\Fib}{\on{Fib}}
\nc{\initial}{\varnothing}
\nc{\op}{\mathrm{op}}

\let\bigcoprod=\coprod
\renc{\coprod}{\sqcup}

\nc{\bDelta}{\mbf{\Delta}}
\nc{\DM}{\mbf{DM}}
\nc{\eff}{\mathrm{eff}}
\nc{\veff}{\mathrm{veff}}
\nc{\cyc}{{\mrm{cyc}}}
\nc{\corr}{{\on{corr}}}
\nc{\fet}{{\mrm{f\acute et}}}
\nc{\fsyn}{{\mrm{fsyn}}}
\nc{\syn}{{\mrm{syn}}}
\nc{\Perf}{\mbf{Perf}}
\nc{\perf}{\on{perf}}
\nc{\oblv}{\on{oblv}}
\nc{\exact}{\on{exact}}
\def\C{\mrm{C}}
\nc{\F}{{\on{F}}}
\nc{\clopen}{{\mrm{clopen}}}
\nc{\B}{\mrm{B}}
\nc{\D}{\mrm{D}}
\nc{\Fin}{\on{Fin}}
\nc{\Cut}{\on{Cut}}
\nc{\Cart}{\on{Cart}}
\nc{\pairs}{\mathsf{pairs}}
\nc{\Pairs}{\mathrm{Pair}}
\nc{\Trip}{\mathrm{Trip}}
\nc{\Lab}{\mathrm{Lab}}
\nc{\coCart}{\mathrm{coCart}}
\nc{\RKE}{\mathrm{RKE}}
\nc{\strict}{\mathrm{strict}}
\nc{\Emb}{\mathrm{Emb}}
\nc{\Split}{\mathrm{Split}}
\nc{\Set}{\mathrm{Set}}
\nc{\sSets}{\mathrm{sSets}}
\nc{\pb}{\mathrm{pb}}
\nc{\fib}{\mathrm{fib}}
\nc{\diff}{\mrm{diff}}
\nc{\gp}{\mrm{gp}}
\nc{\chr}{\mrm{char}}
\nc{\mgp}{\mrm{mot-gp}}
\nc{\FSyn}{\mrm{FSyn}}
\nc{\FSYN}{\mathcal{FS}\mrm{yn}}
\nc{\FEt}{\mrm{FEt}}
\nc{\Spc}{\mrm{Spc}}
\nc{\Ob}{\mrm{Ob}}
\nc{\Spt}{\mrm{Spt}}
\nc{\T}{\bT}
\nc{\suspinf}{\Sigma^\infty}
\nc{\h}{\mrm{h}}
\nc{\uhom}{\underline{\mathrm{Hom}}}
\nc{\umap}{\underline{\mathrm{Maps}}}
\renc{\H}{\bH}
\nc{\Einfty}{{\sE_\infty}}
\nc{\Eone}{{\sE_1}}
\nc{\Stab}{\mrm{Stab}}
\nc{\lax}{{\mrm{lax}}}
\nc{\cocart}{{\mrm{cocart}}}
\nc{\Sch}{\on{Sch}}
\nc{\Fr}{\on{Fr}}
\nc{\A}{\mathbf{A}}
\nc{\N}{\mathbf{N}}
\nc{\Z}{\mathbf{Z}}
\nc{\Q}{\mathbf{Q}}
\nc{\Oo}{\mathcal{O}} 
\nc{\red}{{\on{red}}}
\nc{\Voev}{{\on{Voev}}}
\nc{\Corr}{\mrm{Corr}}
\nc{\Span}{\mathbf{Corr}}
\nc{\Gap}{\mrm{Gap}}
\nc{\Corrfr}{\Corr^{\fr}}
\nc{\Corrvfr}{\Corr^{\Vfr}}
\nc{\Spec}{\on{Spec}}
\nc{\Sm}{\on{Sm}}
\nc{\Gm}{\mathbf{G}_{\on{m}}}
\renc{\P}{\bP}
\nc{\nis}{\mathrm{nis}}
\nc{\zar}{\mathrm{zar}}
\nc{\et}{\mathrm{\acute et}}
\nc{\all}{\mathrm{all}}
\nc{\fold}{\mathrm{fold}}
\nc{\Fun}{\mathrm{Fun}}
\nc{\Ho}{\mathrm{Ho}}
\nc{\Segal}{\mathrm{Segal}}
\nc{\Mon}{\mrm{Mon}{}}
\nc{\Ab}{\mrm{Ab}}
\nc{\Sh}{\on{Sh}}
\nc{\M}{\mrm{M}}
\nc{\Lhtp}{L_{\A^1}}
\nc{\Lmot}{L_{\mrm{mot}}}
\nc{\mot}{\mrm{mot}}
\nc{\SH}{\mbf{SH}}
\nc{\RR}{\mbf{R}}
\nc{\CC}{\mbf{C}}
\nc{\Mod}{\mbf{Mod}}
\nc{\QCoh}{\mbf{QCoh}}
\nc{\MonUnit}{\mbf{1}}
\nc{\tr}{\on{tr}}
\nc{\vop}{\mrm{vop}}
\nc{\fr}{{\on{fr}}}
\nc{\Ar}{\mrm{Ar}}
\nc{\Vfr}{\on{Vfr}}
\nc{\frdiff}{{\on{frdiff}}}
\nc{\frGys}{\on{frGys}}
\nc{\SHfr}{\SH^{\fr}}
\nc{\SHfrdiff}{\SH^{\frdiff}}
\nc{\SHfrGys}{\SH^{\frGys}}
\nc{\InftyCat}{\infty\textnormal{-}\mrm{Cat}}
\nc{\TriCat}{\mathrm{TriCat}}
\nc{\Cat}{\mathrm{1\textnormal{-}Cat}}
\nc{\Th}{\on{Th}}
\def\G{\bG}
\nc{\CMon}{\mrm{CMon}{}}
\nc{\MGL}{\mrm{MGL}}
\nc{\Seg}{\mrm{Seg}{}}
\nc{\Tw}{\mrm{Tw}}
\nc{\sslash}{/\mkern-6mu/}
\nc{\PrL}{\mrm{Pr}^\mrm{L}}
\nc{\PrR}{\mrm{Pr}^\mrm{R}}
\nc{\pr}{\mrm{pr}}
\nc{\efr}{\mrm{efr}}
\nc{\nfr}{\mrm{nfr}}
\nc{\dfr}{\mrm{fr}}
\nc{\tfr}{\mrm{tfr}}
\nc{\Vect}{\mrm{Vect}}
\nc{\sVect}{\mrm{sVect}}
\nc{\fix}{\mrm{fix}}
\nc{\Hilb}{\mathrm{Hilb}}
\nc{\flci}{\mathrm{flci}}
\nc{\Isom}{\mathrm{Isom}}
\nc{\GL}{\mathrm{GL}}
\nc{\fin}{\mathrm{fin}}

\let\phi\varphi
\let\pt\ast
\let\emptyset\varnothing

\nc{\inftyCat}{\term{$\infty$-category}}
\nc{\inftyCats}{\term{$\infty$-categories}}

\nc{\inftyOneCat}{\term{$(\infty,1)$-category}}
\nc{\inftyOneCats}{\term{$(\infty,1)$-categories}}

\nc{\inftyGrpd}{\term{$\infty$-groupoid}}
\nc{\inftyGrpds}{\term{$\infty$-groupoids}}

\nc{\inftyTop}{\term{$\infty$-topos}}
\nc{\inftyTops}{\term{$\infty$-toposes}}

\nc{\inftyTwoCat}{\term{$(\infty,2)$-category}}
\nc{\inftyTwoCats}{\term{$(\infty,2)$-categories}}

\title{Motivic infinite loop spaces}

\author[E. Elmanto]{Elden Elmanto}
\address{Department of Mathematics\\
Harvard University\\
1 Oxford St.\\
Cambridge, MA 02138\\
USA}
\email{\href{mailto:elmanto@math.harvard.edu}{elmanto@math.harvard.edu}}
\urladdr{\url{https://www.eldenelmanto.com/}}
\thanks{E.E.\ and A.K.\ were supported by Institut Mittag-Leffler postdoctoral fellowships}

\author[M. Hoyois]{Marc Hoyois}
\address{Fakultät für Mathematik\\
Universität Regensburg\\
Universitätsstr. 31\\
93040 Regensburg\\
Germany}
\email{\href{mailto:marc.hoyois@ur.de}{marc.hoyois@ur.de}}
\urladdr{\url{http://www.mathematik.ur.de/hoyois/}}
\thanks{M.H.\ was partially supported by NSF grants DMS-1508096 and DMS-1761718}

\author[A. A. Khan]{Adeel A. Khan}
\address{\parbox{\linewidth}{%
IHES,
35 route de Chartres,
91440 Bures-sur-Yvette,
France\\%
Institute of Mathematics,
Academia Sinica,
Taipei 10617,
Taiwan}\vspace{0.2em}}
\email{\href{mailto:khan@ihes.fr}{khan@ihes.fr}}
\urladdr{\url{https://www.preschema.com}}

\author[V. Sosnilo]{Vladimir Sosnilo}
\address{Laboratory ``Modern Algebra and Applications''\\
St. Petersburg State University\\
14th line, 29B\\
199178 Saint Petersburg\\
Russia}
\email{\href{mailto:vsosnilo@gmail.com}{vsosnilo@gmail.com}}
\thanks{V.S.\ was supported by the grant of
the Government of the Russian Federation for the state support of scientific research carried out
under the supervision of leading scientists, agreement 14.W03.31.0030 dated 15.02.2018}

\author[M. Yakerson]{Maria Yakerson}
\address{Institute for Mathematical Research (FIM)\\
ETH Z\"urich \\
R\"amistr. 101\\  
8092 Z\"urich\\
Switzerland}
\email{\href{mailto:maria.yakerson@math.ethz.ch}{maria.yakerson@math.ethz.ch}}
\urladdr{\url{https://www.muramatik.com}}
\thanks{M.Y.\ was supported by SFB/TR 45 ``Periods, moduli spaces and arithmetic of algebraic varieties''}

\date{\today}

\begin{document}

\begin{abstract}
We prove a recognition principle for motivic infinite $\P^1$-loop spaces over a perfect field.
This is achieved by developing a theory of \emph{framed motivic spaces}, which is a motivic analogue of the theory of $\Einfty$-spaces.
A framed motivic space is a motivic space equipped with transfers along finite syntomic morphisms with trivialized cotangent complex in $K$-theory. Our main result is that grouplike framed motivic spaces are equivalent to the full subcategory of motivic spectra generated under colimits by suspension spectra. As a consequence, we deduce some representability results for suspension spectra of smooth varieties, and in particular for the motivic sphere spectrum, in terms of Hilbert schemes of points in affine spaces.
\end{abstract}

\dedicatory{In Memory of Vladimir Voevodsky}
\maketitle

\vspace{2em}
\parskip 0.2cm

\parskip 0pt
\tableofcontents

\parskip 0.2cm


\section{Introduction}
 
 This paper answers the following question in motivic homotopy theory: what sort of algebraic structure characterizes an infinite $\P^1$-loop space? Our answer is: transfers along framed finite syntomic morphisms. Before we elaborate on what we mean by this, let us briefly review the analogous question in classical homotopy theory.

\ssec{The recognition principle in ordinary homotopy theory}

\sssec{} 
Recall that a spectrum is a pointed space equipped with infinitely many successive deloopings. More precisely, it is sequence of pointed spaces $(X_0,X_1,X_2,\dots)$ with equivalences $X_i\simeq \Omega X_{i+1}$, where $\Omega X=\Hom_*(S^1,X)$ is the loop space of $X$.

If $E=(X_0,X_1,\dots)$ is a spectrum, its underlying space $X_0=\Omega^\infty E$ has a very rich structure. Being a loop space, it admits a multiplication $m\colon X_0 \times X_0 \rightarrow X_0$ which is associative and has inverses up to homotopy. Being a \emph{double} loop space, it has a second multiplication, which must agree with the first one and be commutative up to homotopy by the Eckmann–Hilton argument. Being an \emph{infinite} loop space means that $X_0$ has infinitely many associative multiplications that commute with one another: it is an \emph{$\Einfty$-space}.

The \emph{recognition principle} for spectra states that, if $E$ is connective, i.e., if each space $X_i$ is $(i-1)$-connected, then this additional structure on $X_0$ determines the spectrum $E$ up to homotopy. Several approaches were devised to make this precise. The first statement and proof of the recognition principle are due to Boardman and Vogt \cite{BV}. The name ``recognition principle'' was coined by May, who proved a version of it using his notion of operad \cite{may-1972}.

The point of view that we shall adopt is that of Segal \cite{segal1974categories}. Segal introduces the category $\Gamma$, which is the opposite of the category $\Fin_*$ of finite pointed sets. Alternatively, one can think of $\Gamma$ as the category whose objects are finite sets and whose morphisms are spans $X\leftarrow Y\rightarrow Z$ where $Y\to Z$ is injective. The \emph{Segal maps} are the injective maps $\{i\}\hook \{1,\dots,n\}$ in this category. The category $\Gamma$ is now the standard way to define an $\Einfty$-object in any $\infty$-category $C$ with finite products \cite[Definition 2.4.2.1]{HA}: it is a functor $X\colon \Gamma^\op\to C$ such that the Segal maps induce an equivalence $X(\{1,\dots,n\})\simeq \prod_{i=1}^n X(\{i\})$ for any $n$.

Segal's recognition principle is then the following equivalence of $\infty$-categories \cite[Proposition 3.4]{segal1974categories}:
\begin{equation} \label{eqn:spaces-spt}
\Spt_{\geq 0} \simeq \Pre_{\mathrm{Seg}}(\Gamma)^{\gp}.
\end{equation}
The left hand side is the $\infty$-category of connective spectra, and the right-hand side is the $\infty$-category of \emph{grouplike $\Einfty$-spaces}, i.e., presheaves of spaces on $\Gamma$ that satisfy Segal's condition and such that the induced monoid structure has inverses up to homotopy.

\sssec{}

\def\Mfd{\mathrm{Mfd}}

There are many variants of Segal's theorem. Let us mention one with a more geometric flavor, involving the category $\Mfd$ of smooth manifolds. Every space $X$ gives rise to a presheaf $\h_X\colon \Mfd^\op\to\Spc$ sending $M$ to the mapping space $\Maps(M,X)$, which can be viewed as a (very) generalized cohomology theory for smooth manifolds. The presheaf $\h_X$ has two special properties: it is \emph{local} (i.e., it is a sheaf) and it is \emph{$\RR$-homotopy invariant}. 
Using the fact that smooth manifolds are locally contractible, it is easy to see that the construction $X\mapsto\h_X$ defines an equivalence of $\infty$-categories
\begin{equation}\label{eqn:spc-manifolds}
\Spc \simeq \Pre_{\mathrm{loc,\RR}}(\Mfd).
\end{equation}

If $X=\Omega^\infty E$, then the cohomology theory $\h_X$ acquires additional features. For instance, if $f\colon M\to N$ is a finite covering map between smooth manifolds, there is transfer map \[f_*\colon \h_X(M)\to \h_X(N)\] induced by Atiyah duality. 
This enhanced functoriality makes $\h_X$ into a presheaf on a $2$-category $\Span^\fin(\Mfd)$, whose objects are smooth manifolds and whose morphisms are \emph{correspondences} $M\leftarrow N\rightarrow P$, where $M\leftarrow N$ is a finite covering map. It turns out that for a sheaf $\sF$ on $\Mfd$, an extension of $\sF$ to $\Span^\fin(\Mfd)$ is just another way of encoding an $\Einfty$-structure on $\sF$ (see \cite[\sectsign C.1]{norms}).
From this point of view, Segal's equivalence \eqref{eqn:spaces-spt} becomes an equivalence
\begin{equation}\label{eqn:Segal-manifolds}
\Spt_{\geq 0} \simeq \Pre_{\mathrm{loc,\RR}}(\Span^\fin(\Mfd))^\gp.
\end{equation}
In other words, connective spectra are equivalent to homotopy invariant sheaves on smooth manifolds with transfers along finite covering maps, with the additional condition that the monoid structure induced by these transfers is grouplike.
One can also show that the smash product of connective spectra is induced by the Cartesian product of smooth manifolds via Day convolution.

\begin{rem}
	The fact that transfers along finite covering maps characterize infinite loop spaces was conjectured by Quillen before being \emph{disproved} by Kraines and Lada \cite{kraines-lada}. This is because Quillen did not have the language to express the higher coherences necessary for the validity of the result, and therefore did not demand them.
\end{rem}

\ssec{The motivic recognition principle}

\sssec{}
In \cite{MorelAsterisque,MV}, Morel and Voevodsky construct the $\infty$-category $\H(S)$ of \emph{motivic spaces} over a base scheme $S$, as an algebro-geometric analog of the classical homotopy theory of spaces. 
A motivic space over $S$ is by definition a presheaf of spaces on the category $\Sm_S$ of smooth $S$-schemes that is local with respect to the Nisnevich topology and $\A^1$-homotopy invariant:
\[
\H(S) = \Pre_{\nis,\A^1}(\Sm_S).
\]
This definition can be compared with~\eqref{eqn:spc-manifolds}.
The $\infty$-category $\SH(S)$ of \emph{motivic spectra} over $S$ is a certain stabilization of $\H(S)$ with respect to the pointed projective line $\P^1$, introduced by Voevodsky \cite{Voevodsky:1998}. His construction imitates that of the $\infty$-category of spectra in topology: a motivic spectrum is a sequence of motivic spaces $(X_0,X_1,X_2,\dots)$ with equivalences $X_i\simeq \Hom_*(\P^1,X_{i+1})$ expressing $X_{i+1}$ as a $\P^1$-delooping of $X_{i}$.

\sssec{} Our main theorem is a direct analog of~\eqref{eqn:Segal-manifolds} in motivic homotopy theory. 
We construct an $\infty$-category $\Span^\fr(\Sm_S)$ of \emph{framed correspondences} between smooth $S$-schemes, and we prove the following result:

\begin{thm}[Motivic Recognition Principle, Theorem~\ref{thm:main}]
\label{thm:main-in-intro}
Let $k$ be a perfect field. There is an equivalence of symmetric monoidal $\infty$-categories 
\[\SH^{\veff}(k)\simeq \Pre_{\nis,\A^1}(\Span^{\fr}(\Sm_k))^{\gp} .\]
\end{thm}

Here, $\SH^{\veff}(k)$ is the $\infty$-category of \emph{very effective} motivic spectra \cite[Section 5]{SpitzweckOstvaer}, which  is the appropriate analog of connective spectra in motivic homotopy theory: it is the subcategory of $\SH(k)$ generated under colimits and extensions by $\P^1$-suspension spectra. 
To obtain the \emph{stable} $\infty$-category generated under colimits by $\P^1$-suspension spectra, i.e., the $\infty$-category of \emph{effective} motivic spectra \cite[Section 2]{Voevodsky:2002}, it suffices to replace presheaves of spaces with presheaves of spectra: 
\[\SH^{\eff}(k)\simeq \Pre_{\nis,\A^1}(\Span^{\fr}(\Sm_k), \Spt) .\]
We will discuss the $\infty$-category $\Span^{\fr}(\Sm_k)$ in more details below. 

\sssec{}
The problem of formulating a recognition principle for infinite $\P^1$-loop spaces was first posed by Voevodsky: in his collection of open problems in stable motivic homotopy theory \cite[Section 8]{Voevodsky:2002}, he mentions the ``hypothetical theory of operadic description of $\P^1$-loop spaces''. While many of the conjectures in \emph{loc. cit.} have been proved by other means, such a theory has been elusive. In fact, an \emph{operadic} approach to this problem remains so.

Voevodsky himself took important steps towards a recognition principle in several unpublished notes from the early 2000's (for example \cite{voevodsky2001notes}). A major next step was taken by Garkusha and Panin in \cite{garkusha2014framed}, where they use Voevodsky's ideas to explicitly compute $\P^1$-suspension spectra of smooth schemes. As we shall explain, our result relies heavily on their work.

\sssec{}
We now explain the idea behind the $\infty$-category $\Span^\fr(\Sm_k)$ that appears in Theorem~\ref{thm:main-in-intro}.
By design, the infinite $\P^1$-suspension functor
\[
\Sigma^\infty_{\P^1}\colon \H(S)_* \to \SH(S)
\]
sends $\P^1$ to a $\otimes$-invertible object. More generally, for every vector bundle $V$ over $S$, its fiberwise one-point compactification or Thom space $\Th_S(V)$ admits a $\otimes$-inverse $\Th_S(-V)$ in $\SH(S)$. Even more generally, for any element $\xi\in K(S)$ in the $K$-theory space of $S$, there is a well-defined Thom spectrum $\Th_S(\xi)\in\SH(S)$.

\sssec{}
Just as in Grothendieck's original theory of motives, this $\P^1$-inversion procedure leads to powerful duality theorems, which imply in particular that the cohomology theory represented by a motivic spectrum admits certain transfers.
To describe them, we need to introduce some notation.
 For a motivic spectrum $E\in\SH(S)$, a morphism $f\colon X\to S$, and an element $\xi\in K(X)$, let us write
\[
E(X,\xi) = \Maps_{\SH(X)}(\mathbf 1_X, \Th_X(\xi)\wedge f^*E) \in \Spc
\]
for the $\xi$-twisted $E$-cohomology space of $X$. By definition, the underlying space $\Omega^\infty_{\P^1}E$ of $E$ is the presheaf $X\mapsto E(X,0)$ on $\Sm_S$.

The motivic version of Atiyah duality, proved by J.~Ayoub \cite{Ayoub}, implies that for any smooth proper morphism $f\colon Y \to X$ between $S$-schemes and any $\xi\in K(X)$, there is a canonical Gysin transfer
\[
f_!\colon E(Y,f^*\xi+\Omega_f) \to E(X,\xi),
\]
where $\Omega_f$ is the relative cotangent sheaf. When $\Omega_f\simeq 0$, or equivalently when $f$ is finite étale, we obtain in particular a transfer $f_!\colon E(Y,0)\to E(X,0)$. Thus, $\Omega^\infty_{\P^1}E$ is a presheaf with finite étale transfers.

\sssec{}
The transfers discussed so far are completely analogous to those found in the cohomology of smooth manifolds.
 In this algebro-geometric setting, however, there exist even more transfers coming from the theory of \emph{motivic fundamental classes} \cite{LevineVFC,DJKFundamental}. Namely, if $E\in\SH(S)$ and $f\colon Y\to X$ is a proper local complete intersection morphism of $S$-schemes, there is a Gysin transfer
\[
f_!\colon E(Y,f^*\xi+\sL_f)\to E(X,\xi),
\]
where $\sL_f$ is the cotangent complex of $f$. Since $f$ is a local complete intersection, $\sL_f$ is a perfect complex and gives a well-defined element in $K(Y)$. 
While the cotangent complex $\sL_f$ is zero if and only if $f$ is finite étale, its image in $K$-theory may still be zero if $f$ is finite and ramified. A trivialization of $\sL_f$ in $K$-theory is what we will call a \emph{tangential framing} of $f$. A finite lci morphism that admits a tangential framing is necessarily \emph{syntomic}, i.e., it is flat in addition to being lci.

\sssec{}\label{sssec:Corrfr-intro}
In summary, if $E\in \SH(S)$, then for any finite syntomic morphism $f\colon Y\to X$ over $S$ with a trivialization $\alpha\colon \sL_f\simeq 0$ in $K(Y)$, there is a canonical transfer 
\[
(f,\alpha)_!\colon E(Y,0)\to E(X,0).
\]
The $\infty$-category $\Span^\fr(\Sm_S)$ encodes precisely these kinds of transfers: its objects are smooth $S$-schemes, and a morphism from $X$ to $Y$ is a correspondence
\begin{equation*}
  \begin{tikzcd}
     & Z \ar[swap]{ld}{f}\ar{rd}{g} & \\
    X &   & Y
  \end{tikzcd}
\end{equation*}
over $S$ where $f$ is finite syntomic, together with a tangential framing of $f$; it is an $\infty$-category because trivializations of $\sL_f$ in $K(Z)$ form an $\infty$-groupoid.
The construction of this $\infty$-category is the subject of \secref{sect:infinity-category} and is one of the main technical achievements of this paper, although it will be no surprise to the seasoned practitioners of higher category theory.

Theorem~\ref{thm:main-in-intro} can thus be informally summarized as follows: 
a motivic space is an infinite $\P^1$-loop space if and only if it admits a theory of fundamental classes for finite syntomic morphisms (and is grouplike).
The precise connection between Theorem~\ref{thm:main-in-intro} and fundamental classes will be discussed in a sequel to this paper \cite{EHKSY2}.

\sssec{}
Cohomology theories with framed finite syntomic transfers are quite common ``in nature''. By the motivic recognition principle, if such a cohomology theory is Nisnevich-local and $\A^1$-invariant on $\Sm_k$, then it is automatically representable by an effective motivic spectrum over $k$. For example:
\begin{itemize}
	\item If $\mathcal F$ is a presheaf with transfers in the sense of Voevodsky, then in particular $\mathcal F$ has framed finite syntomic transfers (see \ssecref{ssec:fr-to-cyc}). The corresponding motivic spectrum is the usual Eilenberg–Mac Lane spectrum associated with a presheaf with transfers.
	\item More generally, if $\mathcal F$ is a presheaf with Milnor–Witt transfers in the sense of Calmès and Fasel \cite{Calmes:2014ab}, then $\mathcal F$ has framed finite syntomic transfers (this is proved in the sequel \cite{EHKSY2}).
	\item Algebraic $K$-theory admits transfers along finite flat morphisms. The corresponding motivic spectrum is the effective cover of Voevodsky's $K$-theory spectrum $\mathrm{KGL}$.
	\item Witt theory and Hermitian $K$-theory admit transfers along finite syntomic morphisms with trivialized canonical bundle. Since the canonical bundle is the determinant of the cotangent complex, they admit in particular framed finite syntomic transfers.
	\item The universal example of a presheaf with finite syntomic transfers sends $X$ to the groupoid of finite syntomic $X$-schemes. In \cite{EHKSY3}, we show that the corresponding motivic spectrum is Voevodsky's algebraic cobordism spectrum $\mathrm{MGL}$.
\end{itemize}

\ssec{Framed correspondences}

\sssec{}
Even though the theory of motivic fundamental classes motivated our approach to the recognition principle, we do not actually use  it in our proof. 
Instead, the core of the proof is a comparison theorem between Voevodsky's notion of framed correspondence introduced in \cite{voevodsky2001notes} and the notion of tangentially framed correspondence described in~\sssecref{sssec:Corrfr-intro}, as we now explain.

\sssec{}\label{sssec:efr-intro}
If $X$ and $Y$ are $S$-schemes, a framed correspondence from $X$ to $Y$ in the sense of Voevodsky is a correspondence
\begin{equation*}
  \begin{tikzcd}
     & Z \ar[swap]{ld}{f}\ar{rd}{g} & \\
    X &   & Y
  \end{tikzcd}
\end{equation*}
where $f$ is finite, equipped with the following additional data:

\noindent{(i)}
an integer $n\geq 0$ and an embedding $Z\hook \A^n_X$,

\noindent{(ii)}
an \'etale neighborhood $U$ of $Z$ in $\A^n_X$,

\noindent{(iii)}
a morphism $\varphi\colon U \to \A^n$ such that $\varphi^{-1}(0)=Z$ as closed subschemes of $U$, and

\noindent{(iv)}
an $S$-morphism $U \to Y$ lifting $g$.

In this paper, we will call such a correspondence an \emph{equationally framed correspondence} from $X$ to $Y$, to emphasize that the framing is performed by the equations $\varphi_1=0,\dots,\varphi_n=0$. 
Two equationally framed correspondences from $X$ to $Y$ are considered equivalent if they only differ by refining the étale neighborhood $U$ or by increasing $n$. 

We denote by $\Corr^\efr_S(X,Y)$ the set of equationally framed correspondences from $X$ to $Y$. Voevodsky showed that there is a canonical isomorphism
\[ 
\Corr^\efr_S(X,Y) \simeq \colim_{n\to\infty} \Maps\left (X_+\wedge (\P^1)^{\wedge n},L_\nis\left(Y_+\wedge \frac{\A^{n}}{\A^n-0}\right)\right),
\] 
where $L_\nis$ is Nisnevich sheafification (see Appendix~\ref{app:voevodsky-lemma}). Given that $(\P^1)^{\wedge n}\simeq \A^n/(\A^n-0)$ in the pointed motivic homotopy category, $\Corr^\efr_S(X,Y)$ can be viewed as an ``$\A^1$-homotopy-free'' version of the mapping space 
\[
\Maps(\Sigma^\infty_{\P^1}X_+,\Sigma^\infty_{\P^1}Y_+)
\] 
in $\SH(S)$, which is given by the same formula with $L_\mot$ instead of $L_\nis$.

\sssec{}
Using Voevodsky's notion of framed correspondences, Ananyevskiy, Garkusha, Neshitov, and Panin developed the theory of \emph{framed motives} in a series of papers \cite{garkusha2014framed,hitr,agp,gnp}. In particular, they prove framed versions of some of the main theorems in Voevodsky's theory of motives \cite{voevodsky-2000}. 

One of the main outcomes of the theory of framed motives is that, when $S$ is the spectrum of a perfect field $k$, one can compute somewhat explictly the mapping space between $\P^1$-suspension spectra of smooth $k$-schemes using equationally framed correspondences:

\begin{thm}[Garkusha–Panin, see also Corollary~\ref{cor:main}]
	\label{thm:gp-intro}
	Let $k$ be a perfect field. For any $Y\in \Sm_k$, there is an equivalence
\[
\Maps(\Sigma^\infty_{\P^1}(-)_+, \Sigma^\infty_{\P^1}Y_+) \simeq L_\zar (\Lhtp\Corr^\efr_k(-,Y))^\gp
\]
of presheaves of spaces on $\Sm_k$.
\end{thm}

Here, $\Lhtp$ is the naive $\A^1$-localization functor defined by
\[
(\Lhtp\sF)(X) = \lvert \sF(X\times\Delta^\bullet) \rvert,
\]
where $\Delta^n\simeq\A^n$ is the standard algebraic $n$-simplex, and ``gp'' denotes group completion with respect to an explicit $\Einfty$-structure on $\Lhtp\Corr^\efr_k(-,Y)$.

\sssec{}
Let $\Corr^\fr_S(X,Y)$ denote the space of tangentially framed correspondences from $X$ to $Y$, which is the mapping space from $X$ to $Y$ in the $\infty$-category $\Span^\fr(\Sm_S)$. Given an equationally framed correspondence as in \sssecref{sssec:efr-intro}, the morphism $f$ is necessarily syntomic and the morphism $\phi$ determines a trivialization of its cotangent complex in $K(Z)$. We therefore have a forgetful map
\[
\Corr^\efr_S(X,Y) \to \Corr^\fr_S(X,Y),
\]
and we prove the following comparison theorem:

\begin{thm}[Corollary \ref{cor:efr-vs-dfr}]
	\label{thm:comparison-intro}
	Let $Y$ be a smooth $S$-scheme. Then the natural transformation
	\[
	\Corr^\efr_S(-,Y) \to \Corr^\fr_S(-,Y)
	\]
	is a motivic equivalence of presheaves on $\Sch_S$.
\end{thm}

In other words, Nisnevich-locally and up to $\A^1$-homotopy, an equational framing is equivalent to its induced tangential framing.

\sssec{}
Modulo the construction of the $\infty$-category $\Span^\fr(\Sm_S)$ and further technical details, Theorem~\ref{thm:main-in-intro} follows from Theorems \ref{thm:gp-intro} and~\ref{thm:comparison-intro}.

The reason that it is necessary to pass from equationally framed to tangentially framed correspondences in order to state the recognition principle is that the sets $\Corr^\efr_S(X,Y)$ are \emph{not} the mapping sets of a category. Indeed, the embedding $Z\hook \A^n_X$ that is part of the data of an equationally framed correspondence makes it impossible to compose such correspondences. The basic idea leading to the proof of Theorem~\ref{thm:comparison-intro} is that the space of embeddings of $Z$ into $\A^n_X$ becomes $\A^1$-contractible as $n\to\infty$, which suggests that these embeddings can be ignored. The main difficulty is then to show that the rest of the data in an equationally framed correspondence can also be simplified so as not to depend on the embedding of $Z$.

\sssec{} 
The proof of Theorem~\ref{thm:comparison-intro} goes through yet another type of framed correspondences, which is interesting in its own right. Given two $S$-schemes $X$ and $Y$, a \emph{normally framed correspondence} from $X$ to $Y$ is a correpondence
\begin{equation*}
  \begin{tikzcd}
     & Z \ar[swap]{ld}{f}\ar{rd}{g} & \\
    X &   & Y
  \end{tikzcd}
\end{equation*}
where $f$ is finite syntomic, equipped with the following additional data:

\noindent{(i)}
an integer $n\geq 0$ and an embedding $Z\hook \A^n_X$,

\noindent{(ii)}
a trivialization of the normal bundle of $Z$ in $\A^n_X$.

Two such correspondences are considered equivalent if they only differ by increasing $n$, and the set of normally framed correspondences from $X$ to $Y$ is denoted by $\Corr^\nfr_S(X,Y)$. Every equationally framed correspondence gives rise to a normally framed correspondence, which in turn gives rise to a tangentially framed correspondence. In the proof of Theorem~\ref{thm:comparison-intro}, we show that for $Y\in\Sm_S$ both forgetful maps
\[
\Corr^\efr_S(-,Y) \to \Corr^\nfr_S(-,Y) \to \Corr^\fr_S(-,Y)
\]
are motivic equivalences of presheaves on $\Sch_S$. 
We also show that Theorem~\ref{thm:gp-intro} remains true if $\Corr^\efr$ is replaced with either $\Corr^\nfr$ or $\Corr^\fr$ (see Corollary~\ref{cor:main}).

\sssec{}
The main feature of normally framed correspondences is that, for any $Y\in\Sm_S$, the presheaf $\Corr^\nfr_S(-,Y)$ is \emph{representable} by a smooth ind-$S$-scheme (at least if $Y$ is quasi-projective over $S$).
As a result, when $S$ is the spectrum of a perfect field, we obtain a model for the underlying space of the $\P^1$-suspension spectrum $\Sigma^\infty_{\P^1}Y_+$ as the group completion of the motivic homotopy type of this smooth ind-$S$-scheme. 

Let us describe this model more explicitly in the case $Y=S$, i.e., for the motivic sphere spectrum.
Let $\Hilb^{\fin}(\A^n_S)$ be the Hilbert scheme of points in $\A^n_S$ over $S$: its $T$-points are closed subschemes $Z\subset \A^n_T$ that are finite locally free over $T$. A \emph{framing} of such a closed subscheme $Z$ is an isomorphism $\phi\colon \sN_{Z}\simeq \sO_Z^n$, where $\sN_Z$ is the conormal sheaf of $Z$ in $\A^n_T$. Consider the presheaf
\[
\Hilb^\fr(\A^n_S)\colon \Sch_S^\op \to \Set
\] 
sending $T$ to the set of pairs $(Z,\phi)$, where $Z$ is a $T$-point of $\Hilb^\fin(\A^n_S)$ and $\phi$ is a framing of $Z$. Although the scheme $\Hilb^\fin(\A^n_S)$ is typically highly singular, it turns out that $\Hilb^\fr(\A^n_S)$ is representable by a smooth $S$-scheme, the \emph{Hilbert scheme of framed points} in $\A^n_S$.

Let $\Hilb^\fr(\A^\infty_S)$ denote the ind-$S$-scheme $\colim_n \Hilb^\fr(\A^n_S)$. We will show that its naive $\A^1$-localization $\Lhtp\Hilb^\fr(\A^\infty_S)$ admits a canonical $\Einfty$-structure. Intuitively, it is given by forming unions of finite subschemes of $\A^\infty$, provided they are disjoint, which can always be arranged up to $\A^1$-homotopy.

\begin{thm}[see Theorem \ref{thm:sphere-representability}]
	Let $k$ be a field. Then there is a canonical equivalence of $\Einfty$-spaces
	\[
	\Omega^\infty_{\P^1}\mathbf 1_k \simeq L_\zar(\Lhtp \Hilb^\fr(\A^\infty_k))^\gp.
	\] 
\end{thm}

We can regard this theorem as an algebro-geometric analogue of the description of the topological sphere spectrum in terms of framed $0$-dimensional manifolds and cobordisms.

\subsection{Outline of the paper} 
Each section starts with a summary of its own contents, so we only give here a very brief outline. In~\secref{sect:framesmain}, we study the notions of equationally framed, normally framed, and tangentially framed correspondences, and we show that all three notions are motivically equivalent. In~\secref{sect:recogmain}, we prove the recognition principle, taking the construction of the $\infty$-category $\Span^{\fr}(\Sm_S)$ for granted. More precisely, we give a minimal axiomatic description of the $\infty$-category $\Span^{\fr}(\Sm_S)$ that suffices to prove the recognition principle. In~\secref{sect:infinity-category}, we construct the $\infty$-category $\Span^{\fr}(\Sm_S)$. This section is somewhat technical --- we first describe a general procedure for constructing $\infty$-categories of ``labeled correspondences'', which we then apply to construct $\Span^{\fr}(\Sm_S)$. \secref{sect:app} gives various applications of the recognition principle.

\ssec{Conventions and notation}

We will use the language of $\infty$-categories following \cite{HTT} and \cite{HA}.
We denote by $\InftyCat$ the $\infty$-category of $\infty$-categories, by $\Spc$ the $\infty$-category of spaces (i.e., $\infty$-groupoids), and by $\Spt$ the $\infty$-category of spectra.

If $C$ is an $\infty$-category, we denote by $\Pre(C)$ the $\infty$-category of presheaves on $C$, by $\Pre_\Sigma(C)$ the full subcategory of presheaves that transform finite coproducts into finite products \cite[\sectsign 5.5.8]{HTT}, and by $\Pre_\tau(C)$ the full subcategory of $\tau$-local presheaves for a Grothendieck topology $\tau$ on $C$, in the sense of \cite[Definition 6.2.2.6]{HTT}. We denote by $L_\tau\colon \Pre(C)\to \Pre_\tau(C)$ the $\tau$-localization functor.
We write $\Maps_\C(x,y)$ or $\Maps(x,y)$ for spaces of morphisms in $C$, and we reserve the notation $\Hom(x,y)$ for internal hom objects in a symmetric monoidal $\infty$-category.
We write $\C^{\simeq}$ for the $\infty$-groupoid obtained from $\C$ by discarding non-invertible morphisms.

When $C$ is a suitable category of schemes, we denote by $\Lhtp\colon \Pre(C) \to \Pre(C)$ the functor
\[
(\Lhtp\sF)(X) = \colim_{n\in\Delta^\op} \sF(X\times\A^n),
\]
which is the reflection onto the full subctegory of $\A^1$-invariant sheaves. We denote by $L_\mot$ the motivic localization functor, which is the reflection onto the full subcategory of Nisnevich-local $\A^1$-invariant presheaves.

$\Sch$ will denote the category of schemes. If $S$ is a scheme, $\Sch_S$ is the category of $S$-schemes and $\Sm_S$ that of smooth $S$-schemes. We write $\QCoh(S)$ for the stable $\infty$-category of quasi-coherent sheaves on $S$, $\Perf(S)\subset\QCoh(S)$ for the stable $\infty$-category of perfect complexes, $\Vect(S)\subset\Perf(S)$ for the groupoid of locally free sheaves, and $\Vect_n(S)\subset\Vect(S)$ for the groupoid of locally free sheaves of rank $n$.

\ssec{Acknowledgments}

This work is born out of the many days and nights of fruitful discussions when all five of us were, at various times, participants of the program  ``Algebro-Geometric and Homotopical Methods" at the Institut Mittag-Leffler in Djursholm, Sweden. It is with immense pleasure that we thank the hospitable staff, all the participants, and the organizers Eric Friedlander, Lars Hesselholt, and Paul Arne {\O}stv{\ae}r. 
In an official capacity or otherwise, we were all students of Marc Levine. We would like to thank him for his invaluable encouragement and support. Elmanto would additionally like to thank Marc Levine and his adviser, John Francis, for supporting his visit to Essen in 2016 when shadows of this work were first conceived during discussions with Khan and Sosnilo.
Khan would like to thank Denis-Charles Cisinski and Frédéric Déglise for many helpful discussions about fundamental classes, Gysin maps, and much more.
Yakerson would like to thank Tom Bachmann for answering infinitely many of her questions.
We would also like to thank Alexey Ananyevskiy, Grigory Garkusha, Alexander Neshitov, and Ivan Panin for their work on framed motives which is crucial to our work and for patiently answering our questions, and Tom Bachmann for pointing out a mistake in an earlier version of Appendix~\ref{app:finite-fields}.


\section{Notions of framed correspondences} \label{sect:framesmain}

In this section, we study and compare three notions of \emph{framing} of a correspondence
\begin{equation*}
  \begin{tikzcd}
     & Z \ar[swap]{ld}{f}\ar{rd}{h} & \\
    X &   & Y
  \end{tikzcd}
\end{equation*}
over a base scheme $S$, where $f$ is finite and lci:
\begin{itemize}
	\item \emph{Equational framing} (\ssecref{ssec:voev}): the scheme $Z$ is explicitly cut out by equations in a henselian thickening of $Z$ over $X$ and $Y$. This is the notion introduced by Voevodsky in \cite[\sectsign 2]{voevodsky2001notes} and taken up by Garkusha and Panin in \cite{garkusha2014framed}. The main feature of equationally framed correspondences is that they are directly related to mapping spaces in the stable motivic homotopy $\infty$-category; see Remark~\ref{rem:suspension}. 
	\item \emph{Normal framing} (\ssecref{ssec:nfr}): this is simply an embedding $Z\hook\A^n_X$ with trivialized normal bundle. This version is the most useful to obtain geometric descriptions of the underlying spaces of motivic spectra. Indeed, if $Y$ is quasi-projective over $S$, then the functor sending $X$ to the set of normally framed correspondences from $X$ to $Y$ is ind-representable by quasi-projective $S$-schemes (smooth if $Y$ is smooth); see \ssecref{ssec:rep-results}.
	\item \emph{Tangential framing} (\ssecref{ssec:dfr}): this is a trivialization of the cotangent complex of $f$ in the $K$-theory of $Z$. Unlike the previous two, tangentially framed correspondences are the mapping spaces of an $\infty$-category whose objects are $S$-schemes, which we will construct in \secref{sect:infinity-category}.
\end{itemize} 
The main result of this section is that the presheaves parametrizing all three notions of framed correspondences are motivically equivalent (Corollary~\ref{cor:vfr-vs-nfr2} and Corollary~\ref{cor:nfr-vs-dfr}). 

Throughout this section, we fix an arbitrary base scheme $S$.


\ssec{Equationally framed correspondences}
\label{ssec:voev}

\sssec{}

Recall the following definition from \cite[\sectsign 2]{voevodsky2001notes} (see also \cite[Definition 2.1]{garkusha2014framed}):

\begin{defn} \label{defn:Voevodsky framed correspondence}
Let $X,Y\in\Sch_{S}$ and $n\geq 0$.
An \emph{equationally framed correspondence of level $n$} from $X$ to $Y$ over $S$ consists of the following data:

\noindent{\em(i)}
A closed subscheme $Z \subset \A^n_X$ finite over $X$.

\noindent{\em(ii)}
An \'etale neighborhood $U$ of $Z$ in $\A^n_X$, i.e., an \'etale morphism of algebraic spaces $U \to \A^n_X$ inducing an isomorphism $U\times_{\A^n_X}Z\simeq Z$.

\noindent{\em(iii)}
A morphism $\varphi\colon U \to \A^n$ such that $\varphi^{-1}(0)=Z$ as closed subschemes of $U$.

\noindent{\em(iv)}
An $S$-morphism $g\colon U \to Y$.
\end{defn}

Two equationally framed correspondences $(Z,U,\varphi,g)$ and $(Z',U',\varphi',g')$ are considered equivalent if $Z = Z'$ and the morphisms $g$ and $g'$ (resp. $\varphi$ and $\varphi'$) agree in an \'etale neighborhood of $Z$ refining $U$ and $U'$. We denote by $\Corr^{\efr,n}_S(X,Y)$ the set of equationally framed correspondences of level $n$, modulo this equivalence relation. Note that we have a functor
\[
\Corr^{\efr,n}_S(-,-)\colon \Sch_S^\op \times \Sch_S \to \Set.
\]

An equationally framed correspondence $(Z,U,\phi,g)$ from $X$ to $Y$ contains in particular a span
\begin{equation*}
  \begin{tikzcd}
     & Z \ar[swap]{ld}{f}\ar{rd}{h} & \\
    X &   & Y,
  \end{tikzcd}
\end{equation*}
where $f$ is finite and $h$ is the restriction of $g\colon U\to Y$ to $Z$. 

\begin{rem}\label{rem:voevodsky-lemma}
	Definition~\ref{defn:Voevodsky framed correspondence} is motivated by the following computation of Voevodsky: there is a natural isomorphism
\[
\Corr^{\efr,n}_S(X,Y) \simeq \Maps(X_+\wedge (\P^1)^{\wedge n},L_\nis(Y_+\wedge \T^{\wedge n})),
\]
where $\P^1$ is pointed at $\infty$ and $\T=\A^1/\A^1-0$, and similarly with the étale topology instead of the Nisnevich topology. We give a proof of this result in Appendix~\ref{app:voevodsky-lemma}.  In Voevodsky's formulation, the étale neighborhood $U$ in the definition of an equationally framed correspondence is assumed to be a scheme. This makes no difference if $X$ is affine, or if $X$ is qcqs geometrically unibranch with finitely many irreducible components and $Y$ is separated over $S$ (see Lemma~\ref{lem:henselization}(ii,iv)), but in general there may be too few schematic étale neighborhoods.
\end{rem}

\sssec{} \label{defn:hefr}

\def\pr{\mathrm{pr}}

Let $X,Y\in\Sch_{S}$ and let $n\geq 0$. The \emph{suspension morphism}
	\[\sigma_{X,Y}\colon \Corr^{\efr,n}_S(X,Y)\hook \Corr^{\efr,n+1}_S(X,Y)\]
sends $(Z,U,\varphi,g)$ to $(Z\times\{0\},U\times\A^1,\phi\times\id_{\A^1},g\circ \pr_U)$, where $\pr_U\colon U\times\A^1\to U$ is the projection onto the first factor.
	
We set
\[
\Corr^\efr_S(X,Y) = \colim(\Corr^{\efr,0}_S(X,Y) \xrightarrow{\sigma_{X,Y}} \Corr^{\efr,1}_S(X,Y)  \xrightarrow{\sigma_{X,Y}} \ldots).
\]
Note that $\sigma_{X,Y}$ is natural in both $X$ and $Y$, hence we have a functor
\[
\Corr^\efr_S(-,-)\colon \Sch_S^\op\times \Sch_S\to\Set.
\]
We will denote by $\h^{\efr,n}_S(Y)$ and $\h^\efr_S(Y)$ the presheaves $X\mapsto \Corr^{\efr,n}_S(X,Y)$ and $X\mapsto \Corr^\efr_S(X,Y)$ on $\Sch_S$. 

\begin{prop}\label{prop:fr-descent}
Let $Y\in\Sch_S$ and $n\geq 0$.

\noindent{\em(i)}
$\h^{\efr,n}_S(Y)$ is an étale sheaf and $\h^\efr_S(Y)$ is a sheaf for the quasi-compact étale topology.

\noindent{\em(ii)}
$\h^{\efr,n}_S(Y)$ and $\h^{\efr}_S(Y)$ satisfy closed gluing (see Definition~\ref{defn:closedgluing}).

\noindent{\em(iii)}
Let $R\hook Y$ be a Nisnevich (resp.\ étale) covering sieve generated by a single map. Then $\h^{\efr,n}_S(R)\to \h^{\efr,n}_S(Y)$ and $\h^\efr_S(R)\to \h^\efr_S(Y)$ are Nisnevich (resp.\ étale) equivalences.
\end{prop}

\begin{proof}
	In the notation of \sssecref{sssec:Q}, $\Corr^{\efr,n}_S(X,Y)$ is the fiber of the pointed map
	\[
	Q(\A^n\times Y,(\A^n-0)\times Y)((\P^1)^n\times X) \to \prod_{1\leq i\leq n} Q(\A^n\times Y,(\A^n-0)\times Y)((\P^1)^{n-1}\times X)
	\]
	induced by the $n$ embeddings at infinity $(\P^1)^{n-1}\hook (\P^1)^{n}$.
	Then (i) follows from Proposition~\ref{prop:VoevodskyLemma}, and (ii) from Proposition~\ref{prop:closedgluing}.
	Let us prove (iii). 
	 Note that $\h^{\efr,n}_S(R)(X)$ is the subset of $\h^{\efr,n}_S(Y)(X)$ consisting of those equationally framed correspondences $(Z,U,\phi,g)$ such that $g$ factors through $R$. It thus suffices to show that $\h^{\efr,n}_S(R)\to \h^{\efr,n}_S(Y)$ is a Nisnevich (resp.\ étale) epimorphism.
	  Refining the sieve $R$ if necessary, we may assume that it is generated be a single étale map.
	  Let $(Z,U,\phi,g)\in \h^{\efr,n}_S(Y)(X)$.
	 Let $R_U$ and $R_Z$ be the pullbacks of $R$ to $U$ and $Z$, and let $R_X$ be the sieve on $X$ generated by all maps $X'\to X$ whose pullback to $Z$ belongs to $R_Z$.
	 Since $R_Z$ is generated by a single étale map and $Z\to X$ is finite, $R_X$ is a Nisnevich (resp.\ étale) covering sieve of $X$. For $X'\to X$ in $R_X$, the map $Z\times_XX' \to Z$ belongs to $R_Z$, hence the étale neighborhood $U\times_XX'$ of $Z\times_XX'$ has a refinement in $R_U$. Thus, the pullback of $(Z,U,\phi,g)$ to $X'$ belongs to the subset $\h^{\efr,n}_S(R)(X')\subset \h^{\efr,n}_S(Y)(X')$, as desired.
\end{proof}

In particular, by Proposition~\ref{prop:fr-descent}(i), the presheaves $\h^{\efr,n}_S(Y)$ and $\h^\efr_S(Y)$ transform finite sums into finite products, i.e., they belong to $\Pre_\Sigma(\Sch_S)$.

\begin{rem}\label{rem:suspension}
Under the identification of Remark~\ref{rem:voevodsky-lemma}, the suspension morphism $\sigma_{X,Y}$ sends a pointed map $f\colon X_+\wedge (\P^1)^{\wedge n}\to Y_+\wedge \T^{\wedge n}$ to the composite
\[
X_+\wedge (\P^1)^{\wedge n+1}\xrightarrow{\id\wedge a} X_+\wedge (\P^1)^{\wedge n}\wedge \T \xrightarrow{f\wedge\id} Y_+\wedge \T^{\wedge n+1},
\]
where $a\colon \P^1\to \T$ is the collapse map $\P^1/\infty \to \P^1/\P^1-0 \simeq \A^1/\A^1-0$. In particular, if $X$ and $Y$ are smooth over $S$, there is a natural map
\begin{equation}\label{eqn:fr-to-SH}
	\Corr^\efr_S(X,Y) \to \Maps_{\SH(S)}(\Sigma^\infty_\T X_+,\Sigma^\infty_\T Y_+).
\end{equation}
\end{rem}

\sssec{}
\label{sssec:Corrfr-extended-functoriality}
	It follows from Remarks \ref{rem:voevodsky-lemma} and~\ref{rem:suspension} that $\Corr^{\efr,n}_S(-,-)$ and $\Corr^\efr_S(-,-)$ have canonical extensions to functors
	\[
	\Sch_{S+}^\op \times \Sch_{S+} \to \Set,
	\]
	where $\Sch_{S+}$ is the full subcategory of pointed $S$-schemes of the form $X_+$. Equivalently, $\Sch_{S+}$ is the category whose objects are $S$-schemes and whose morphisms are partially defined maps with clopen domains. It is also straightforward to describe this extended functoriality directly in terms of Definition~\ref{defn:Voevodsky framed correspondence}.

\begin{rem}\label{rem:GPadditivity}
	Proposition~\ref{prop:fr-descent}(iii) does not hold for sieves that are generated by more than one (or zero) maps. For instance, if $n\neq 1$, the canonical map
	\[
	\h^\efr_S(Y_1)\coprod\dotsb\coprod \h^\efr_S(Y_n) \to \h^\efr_S(Y_1\coprod \dotsb\coprod Y_n)
	\]
	is not an étale-local equivalence. Using the extended functoriality from \sssecref{sssec:Corrfr-extended-functoriality}, we instead have a more interesting map
	\[
	\h^\efr_S(Y_1\coprod \dotsb\coprod Y_n) \to \h^\efr_S(Y_1)\times\dotsb\times \h^\efr_S(Y_n).
	\]
	The additivity theorem of Garkusha and Panin \cite[Theorem 6.4]{garkusha2014framed} states that the latter map is an $\Lhtp$-equivalence. 
\end{rem}

\sssec{} \label{sssec:hefr-for-arbitrary} 
We define
\[
\h^{\efr,n}_S, \h^{\efr}_S\colon \Pre_\Sigma(\Sch_S)_* \to \Pre_\Sigma(\Sch_S)_*
\]
to be the sifted-colimit-preserving extensions\footnote{Some care must be taken because $\Sch_S$ is a large category: these extensions are only defined on presheaves that are small sifted colimits of representable presheaves, or else they take values in large presheaves. Alternatively, one can let $\Sch_S$ be a large enough but small category of $S$-schemes. Either way, this has no practical consequences.} of the functors
\[
\h^{\efr,n}_S, \h^{\efr}_S\colon \Sch_{S+} \to \Pre_\Sigma(\Sch_S)_*
\]
from \sssecref{sssec:Corrfr-extended-functoriality}.
	For example, if $Y \in \Sch_S$ and $W \subset Y$ is a subscheme, then the quotient $Y/W$ in $\Pre_{\Sigma}(\Sch_S)_*$ is a simplicial colimit of representable presheaves:
\begin{equation*}
\begin{tikzcd}
Y_+ 
& (Y \coprod W)_+ \arrow[l, shift left]
\arrow[l, shift right]
& (Y \coprod W \coprod W)_+ \arrow[l]
\arrow[l, shift left=
2
]
\arrow[l, shift right=
2
]
& \cdots.
\arrow[l, shift left]
\arrow[l, shift right]
\arrow[l, shift left=3]
\arrow[l, shift right=3]
\end{tikzcd}
\end{equation*}
Therefore, the presheaf $\h^{\efr,n}(Y/W)$ is the colimit of the simplicial diagram
\begin{equation*}
\begin{tikzcd}
\h^{\efr,n}_S(Y) 
& \h^{\efr,n}_S(Y\coprod W) \arrow[l, shift left]
\arrow[l, shift right]
& \h^{\efr,n}_S(Y\coprod W\coprod W) \arrow[l]
\arrow[l, shift left=
2
]
\arrow[l, shift right=
2
]
& \cdots.
\arrow[l, shift left]
\arrow[l, shift right]
\arrow[l, shift left=3]
\arrow[l, shift right=3]
\end{tikzcd}
\end{equation*}

\nc{\Fcal}{\mathcal{F}}

Furthermore, since sifted colimits commute with finite products, the additivity theorem of Garkusha–Panin recalled in Remark~\ref{rem:GPadditivity} implies that for any $\Fcal_1, \dotsc, \Fcal_k \in \Pre_{\Sigma}(\Sch_S)_*$, the map $\h^{\efr}_S(\Fcal_1\vee \cdots \vee \Fcal_k) \rightarrow \h^{\efr}_S(\Fcal_1) \times \cdots \times \h^{\efr}_S(\Fcal_k)$ is an $\Lhtp$-equivalence. 

\begin{rem} \label{rem:desc} 
	We already implicitly used this extension in the discussion of the descent properties of $\h^{\efr,n}_S$ in Proposition~\ref{prop:fr-descent}, in order to make sense of $\h^{\efr, n}_S(R)$ where $R \hookrightarrow Y$ is a sieve. Moreover, Proposition~\ref{prop:fr-descent}(iii) implies that the extension of $\h^{\efr,n}_S$ to $\Pre_{\Sigma}(\Sch_S)$ preserves local equivalences for the Nisnevich and étale topologies. Indeed, by \cite[Lemma 2.10]{norms}, it suffices to show that $\h^{\efr,n}_S$ sends covering sieves to local equivalences. For finitely generated sieves, this follows from the observation that the localization functor $\Pre(\Sch_S) \rightarrow \Pre_{\Sigma}(\Sch_S)$ sends a sieve generated by $U_1,\dotsc,U_n$ to the sieve generated by $U_1\coprod\cdots\coprod U_n$, and one can reduce to the case of finitely generated sieves as in the proof of Corollary~\ref{cor:vfr-vs-nfr2}.
	The same discussion applies to $\h^\efr_S$.
\end{rem}

\sssec{}\label{sssec:efr-module-structure}
If $X,Y\in\Sch_S$ and $n\geq 0$, there is pointed map
\[
X_+\wedge \h^{\efr,n}_S(Y) \to \h^{\efr,n}_S(X\times_S Y)
\]
sending a map $a\colon Z\to X$ and an equationally framed correspondence $(W,U,\phi,g)$ from $Z$ to $Y$ to the equationally framed correspondence $(W,U,\phi,(a\times g)\circ\gamma)$ where $\gamma\colon U\hookrightarrow Z\times_S U$ is the graph of $U\to Z$. Note that the squares
\begin{equation*}
  \begin{tikzcd}
    X_+\wedge \h^{\efr,n}_S(Y) \ar[hookrightarrow]{r}{\id\wedge \sigma_{Y}}\ar{d}
      & X_+\wedge\h^{\efr,n+1}_S(Y) \ar{d}
    \\
    \h^{\efr,n}_S(X\times_SY) \ar[hookrightarrow]{r}{\sigma_{X\times Y}}
      & \h^{\efr,n+1}_S(X\times_SY)
  \end{tikzcd}
\end{equation*}
commute, and we obtain a map 
\[
X_+\wedge \h^{\efr}_S(Y) \to \h^{\efr}_S(X\times_S Y)
\]
in the colimit.
These maps are clearly natural in $X_+$ and $Y_+$ and define right-lax $\Sch_{S+}$-linear structures on the functors $\h^{\efr,n}_S,\h^{\efr}_S\colon \Sch_{S+}\to \Pre_\Sigma(\Sch_{S})_\pt$. Extending to $\Pre_\Sigma$, we obtain right-lax $\Pre_\Sigma(\Sch_S)_\pt$-linear structure on the functors $\h^{\efr,n}_S,\h^{\efr}_S\colon\Pre_\Sigma(\Sch_S)_\pt\to \Pre_\Sigma(\Sch_S)_\pt$.

\sssec{}\label{sssec:efr-unit}
Note that $\h^{\efr,0}_S(Y)$ is the presheaf on $\Sch_S$ represented by $Y_+$. We therefore have a map
\begin{equation*}
Y_+\to\h^\efr_S(Y),
\end{equation*}
natural in $Y_+\in\Sch_{S+}$ and compatible with the $\Sch_{S+}$-linear structures, which extends to a $\Pre_\Sigma(\Sch_S)_*$-linear natural transformation $\id\to \h^{\efr}_S\colon \Pre_\Sigma(\Sch_S)_* \to \Pre_\Sigma(\Sch_S)_*$. The $\Pre_\Sigma(\Sch_S)_*$-linearity means in particular that the following triangle commutes for every $\sF,\sG\in\Pre_\Sigma(\Sch_S)_*$:
\[
\begin{tikzcd}
	\sF\wedge \sG \ar{dr} \ar{r} & \sF\wedge \h^{\efr}_S(\sG) \ar{d} \\
	& \h^{\efr}_S(\sF\wedge\sG).
\end{tikzcd}
\]

\sssec{}
We now investigate the scheme-theoretic properties of equationally framed correspondences. Recall that a morphism of affine schemes $f\colon \Spec S \rightarrow \Spec R$ is a \emph{relative global complete intersection} if there exists a presentation $S = R[x_1, \dotsc, x_n]/(f_1,\dotsc, f_c)$ and the nonempty fibers of $f$ have dimension $n-c$ \cite[Tag 00SP]{stacks}.

\begin{lem}\label{lem:gci}
	Let $X$ be an affine scheme, $U\to \A^n_X$ an affine étale morphism, and $Z\subset U$ a closed subscheme cut out by $c$ equations.
	Suppose that the the nonempty fibers of $Z\to \A^n_X \to X$ have dimension $n-c$. Then $Z\to X$ is a relative global complete intersection.
\end{lem}

\begin{proof}
	Choose a closed immersion $i\colon U\hook \A^{n+m}_X$ over $\A^n_X$. Since $U\to \A^n_X$ is étale, we obtain an isomorphism $\sN_i\simeq i^*(\Omega_{\A^{n+m}_X/\A^n_X})\simeq \sO^m_U$. 
	Choose functions $f_1$, \dots, $f_m$ on $\A^{n+m}_X$ lifting generators of $\sN_i$. By Nakayama's lemma, there exists a function $h$ on $\A^{n+m}_X$ such that $U$ is cut out by $f_1$, \dots, $f_m$ in $(\A^{n+m}_X)_h$. But then $U$ is cut out by $m+1$ equations in $\A^{n+m+1}_X$. Hence, $Z$ is cut out by $c+m+1$ equations in $\A^{n+m+1}_X$. By definition, this implies that $Z\to X$ is a relative global complete intersection.
\end{proof}

We use the notion of \emph{regular immersion} following \cite[VII, Definition 1.4]{SGA6} (as opposed to \cite[Definition 16.9.2]{EGA4-4}). Regular immersions in this sense are called Koszul-regular in \cite[Tag 0638]{stacks}.
Recall that a morphism $f\colon Y\to X$ is a \emph{local complete intersection} (lci) if, locally on $Y$, it is the composition of a regular immersion and a smooth morphism \cite[VIII, Definition 1.1]{SGA6}.

\begin{prop}\label{prop:lci-flat}
	Let $f\colon Y\to X$ be an lci morphism of virtual relative dimension $d\geq 0$.
	
	\noindent{\em(i)} 
	$f$ is flat if and only if its nonempty fibers have dimension $d$.
	
	\noindent{\em(ii)}
	Suppose $X$ and $Y$ affine. Then $f$ is a relative global complete intersection if and only if it is flat and $[\sL_f]=[\sO_Y^d]$ in $K_0(Y)$.
\end{prop}

\begin{proof}
	(i) Suppose $f$ flat. Since flat lci morphisms are stable under base change \cite[Tag 01UI]{stacks}, we can assume that $X$ is the spectrum of a field, in which case it is obvious that $Y$ has dimension $d$. To prove the converse, since the question is local on $X$ and $Y$, we can assume that $X$ and $Y$ are affine. Choose a closed immersion $i\colon Y\hook\A^n_X$ over $X$. Since $f$ is lci, $i$ is regular \cite[Tag 069G]{stacks}, and hence the conormal sheaf $\sN_i$ is locally free of rank $n-d$. Let $U\subset\A^n_X$ be an affine open subscheme where $\sN_i$ is free. By Nakayama's lemma, there exists a function $h$ on $U$ such that $Y_U$ is cut out by $n-d$ equations in $U_h$. Applying Lemma~\ref{lem:gci}, we deduce that $Y_U\to X$ is a relative global complete intersection, hence flat \cite[Tag 00SW]{stacks}. 
	
	(ii) In the above situation, $[\sL_f]=[\sO_Y^n]-[\sN_i]$ in $K_0(Y)$, hence $[\sL_f]=[\sO_Y^d]$ if and only if $\sN_i$ is stably free. In this case, we may assume by increasing $n$ that $\sN_i$ is free, so that $Y\to X$ itself is a relative global complete intersection. Conversely, if $f$ is a relative global complete intersection, then we can choose $i$ so that $\sN_i$ is free \cite[Tag 00SV]{stacks}.
\end{proof}

\sssec{} Recall that a morphism of scheme is \emph{syntomic} if it is flat and lci. A key property of this notion, which we already used in the proof of Proposition~\ref{prop:lci-flat}(i), is that syntomic morphisms are stable under arbitrary base change, unlike lci morphisms.

\begin{prop}\label{prop:framed-corr-properties}
Let $(Z, U, \varphi, g)$ be an equationally framed correspondence from $X$ to $Y$ over $S$, and let $i\colon Z\hook\A^n_X$ be the inclusion.

\noindent{\em(i)}
The conormal sheaf $\sN_i$ is free of rank $n$.

\noindent{\em(ii)}
The immersion $i$ is regular.

\noindent{\em(iii)}
The morphism $Z \to X$ is syntomic.

\noindent{\em(iv)}
If $X$ is affine, then $Z\to X$ is a relative global complete intersection.
\end{prop}

\begin{proof}
Let $i'\colon Z\hook U$ be the lift of $i$.
Since $U\to \A^n_X$ is étale, it induces an isomorphism of conormal sheaves $\sN_i\simeq \sN_{i'}$.
Moreover, the Cartesian square
  \begin{equation*}
    \begin{tikzcd}
      Z \ar[hookrightarrow]{r}{i'}\ar{d}
        & U \ar{d}{\varphi}
      \\
      \{0\} \ar[hookrightarrow]{r}
        & \A^n
    \end{tikzcd}
  \end{equation*}
  induces an epimorphism $\sO^n_Z\twoheadrightarrow \sN_{i'}$.
To prove (i), it remains to show that $\sN_i$ is locally free of rank $n$, which will follow from (ii).
Statements (ii) and (iii) are local on $X$, so we may assume that $X$ is affine. Then we may assume that $U$ is also affine. By Lemma~\ref{lem:gci}, we deduce that $Z\to X$ is a relative global complete intersection. In particular, it is syntomic \cite[Tag 00SW]{stacks} and $i$ is regular \cite[Tag 069G]{stacks}.
\end{proof}


\ssec{Normally framed correspondences}
\label{ssec:nfr}

\sssec{}
\label{sssec:nfr-intro}

We will now introduce a simpler notion of framed correspondences. This notion carries less data than equationally framed correspondences, but it turns out to be equivalent after applying motivic localization (Corollary~\ref{cor:vfr-vs-nfr2}). To motivate this definition, let $X,Y\in\Sch_S$, let $(Z,U,\phi,g)$ be an equationally framed correspondence of level $n$ from $X$ to $Y$ over $S$, and let $i\colon Z\hook\A^n_X$ be the inclusion.
The ideal $\sI$ of $Z$ in $U$ is by definition generated by the components of $\phi\colon U\to \A^n$. In particular, we have an epimorphism
\[
\tau\colon \sO^n_Z \to \sI/\sI^2\simeq \sN_i, \quad e_m\mapsto [\phi_m],
\]
which is necessarily an isomorphism since $\sN_i$ is locally free of rank $n$ (Proposition~\ref{prop:framed-corr-properties}(i)). In other words, $\phi$ induces a trivialization of the conormal sheaf of $i$. 

\begin{defn}\label{defn:normal-framed-corr}
	Let $X,Y\in\Sch_S$ and $n\geq 0$. A \emph{normally framed correspondence of level $n$} from $X$ to $Y$ over $S$ consists of the following data:

\noindent{\em(i)}
A span 
\begin{equation*}
  \begin{tikzcd}
     & Z \ar[swap]{ld}{f}\ar{rd}{h} & \\
    X &   & Y
  \end{tikzcd}
\end{equation*}
over $S$, where $f$ is finite syntomic.

\noindent{\em(ii)}
A closed immersion $i\colon Z\hook \A^n_X$ over $X$.

\noindent{\em(iii)}
A trivialization $\tau\colon \sO_Z^n\simeq\sN_i$ of the conormal sheaf of $i$.
\end{defn}

We denote by $\Corr^{\nfr,n}_S(X,Y)$ the set (or rather discrete groupoid) of normally framed correspondences of level $n$; it is a functor of $X$ and $Y$. The closed immersion $\A^n_X\hook \A^{n+1}_X$, $x\mapsto (x,0)$, induces an obvious suspension morphism
\[
\sigma_{X,Y}\colon\Corr^{\nfr,n}_S(X,Y) \hook \Corr^{\nfr,n+1}_S(X,Y),
\]
and we set
\[
\Corr^\nfr_S(X,Y) = \colim(\Corr^{\nfr,0}_S(X,Y) \xrightarrow{\sigma_{X,Y}} \Corr^{\nfr,1}_S(X,Y)  \xrightarrow{\sigma_{X,Y}} \ldots).
\]

Again, we have a functor
\[
\Corr^\nfr_S(-,-)\colon \Sch_S^\op\times \Sch_S\to\Set.
\]
We will denote by $\h^{\nfr,n}_S(Y)$ and $\h^\nfr_S(Y)$ the presheaves $X\mapsto \Corr^{\nfr,n}_S(X,Y)$ and $X\mapsto \Corr^\nfr_S(X,Y)$ on $\Sch_S$.

\begin{rem}
	In Definition~\ref{defn:normal-framed-corr}(i), we can replace ``syntomic'' by ``lci''. Indeed, (ii) and (iii) imply that $f$ is of virtual relative dimension $0$, hence flat by Proposition~\ref{prop:lci-flat}(i).
\end{rem}

\begin{prop}\label{prop:nfr-descent}
Let $Y\in\Sch_S$ and $n\geq 0$.

\noindent{\em(i)}
$\h^{\nfr,n}_S(Y)$ is an fpqc sheaf and $\h^{\nfr}_S(Y)$ is a sheaf for the quasi-compact fpqc topology.

\noindent{\em(ii)}
$\h^{\nfr,n}_S(Y)$ and $\h^{\nfr}_S(Y)$ satisfy closed gluing (see Definition~\ref{defn:closedgluing}).

\noindent{\em(iii)}
Let $R\hook Y$ be a Nisnevich (resp.\ étale) covering sieve generated by a single map. Then $\h^{\nfr,n}_S(R)\to \h^{\nfr,n}_S(Y)$ and $\h^\nfr_S(R)\to \h^\nfr_S(Y)$ are Nisnevich (resp.\ étale) equivalences.
\end{prop}

\begin{proof}
(i) is clear since finite syntomic maps, closed immersions, vector bundles, and representable presheaves satisfy fpqc descent. The proof of (iii) is similar to that of Proposition~\ref{prop:fr-descent}(iii). Let us prove (ii).
	Suppose $X=X_0\coprod_{X_{01}}X_1$ where $X_{01}\hook X_0$ and $X_{01}\hook X_1$ are closed immersions.
	If $Z\to X$ is flat, it follows from \cite[Lemme 4.4]{ferrand} that $Z=Z_0\coprod_{Z_{01}}Z_1$, where $Z_*=X_*\times_XZ$.
	Moreover, the canonical functor $\Vect(Z)\to \Vect(Z_0)\times_{\Vect(Z_{01})}\Vect(Z_1)$ is an equivalence of groupoids, by \cite[Theorem 16.2.0.1]{SAG}. The claim follows easily.
\end{proof}

In particular, by Proposition~\ref{prop:nfr-descent}(i), the presheaves $\h^{\nfr,n}_S(Y)$ and $\h^\nfr_S(Y)$ transform finite sums into finite products, i.e., they belong to $\Pre_\Sigma(\Sch_S)$.

\sssec{} 

As explained in \sssecref{sssec:nfr-intro}, every equationally framed correspondence $(Z,U,\phi,g)$ from $X$ to $Y$ gives rise to a normally framed correspondence $X\leftarrow Z\to Y$, where $Z\to Y$ is the restriction of $g\colon U\to Y$ to $Z$, $i\colon Z\hook \A^n_X$ is the given inclusion, and the trivialization $\tau\colon \sO_Z^n\simeq\sN_i$ is induced by $\phi$. This defines a map
\begin{equation}\label{eqn:vfr-to-nfr}
	\Corr^{\efr,n}_S(X,Y) \to \Corr^{\nfr,n}_S(X,Y),
\end{equation}
which is easily checked to be natural in both $X$ and $Y$.
Moreover, the squares
\begin{equation*}
  \begin{tikzcd}
    \Corr^{\efr,n}_S(X,Y) \ar[hookrightarrow]{r}{\sigma_{X,Y}}\ar{d}
      & \Corr^{\efr,n+1}_S(X,Y) \ar{d}
    \\
    \Corr^{\nfr,n}_S(X,Y) \ar[hookrightarrow]{r}{\sigma_{X,Y}}
      & \Corr^{\nfr,n+1}_S(X,Y)
  \end{tikzcd}
\end{equation*}
commute. Taking the colimit of~\eqref{eqn:vfr-to-nfr} over $n$, we obtain a natural map
\begin{equation}\label{eqn:vfr-to-nfr2}
	\Corr^{\efr}_S(X,Y) \to \Corr^{\nfr}_S(X,Y).
\end{equation}

\sssec{}
\label{sssec:Corrnfr-extended-functoriality}
The functors $\Corr^{\nfr,n}_S(-,-)$ and $\Corr^\nfr_S(-,-)$ can be extended to
\[
\Sch_{S+}^\op \times \Sch_{S+} \to \Set.
\]
In other words, normally framed correspondences are functorial for partially defined maps with clopen domains. This extended functoriality, in either variable, comes from the fact that clopen embeddings are simultaneously finite syntomic and closed immersions.
Moreover, the maps \eqref{eqn:vfr-to-nfr} and~\eqref{eqn:vfr-to-nfr2} are natural on $\Sch_{S+}^\op\times\Sch_{S+}$, with respect to the extended functoriality of equationally framed correspondences described in \sssecref{sssec:Corrfr-extended-functoriality}.
 We leave the details to the reader.
 
\sssec{}
\label{sssec:nfr-E_infinity}
Via the obvious functor 
 \[
\Sch_{S+}\times \Fin_\pt \to \Sch_{S+},\quad (Y_+,I_+)\mapsto (Y^{\coprod I})_+,
 \]
every functor on $\Sch_{S+}$ is automatically a $\Fin_\pt$-object.
In particular, $\h^{\efr}(Y)$ and $\h^{\nfr}(Y)$ are $\Fin_\pt$-objects and the natural transformation 
\begin{equation} \label{einftymap1}
\h^\efr_S(Y)\to\h^\nfr_S(Y)
\end{equation} 
is a map of $\Fin_\pt$-objects. The additivity theorem of Garkusha and Panin \cite[Theorem 6.4]{garkusha2014framed} implies that the $\Fin_\pt$-object $\Lhtp \h^\efr_S(Y)$ is an $\Einfty$-object. We similarly have:

\begin{prop}\label{prop:nfr-E_infinity}
	Let $Y_1,\dotsc,Y_k\in\Sch_S$. Then the canonical map
	\[
	\h^\nfr_S(Y_1\coprod\dotsb\coprod Y_k) \to \h^\nfr_S(Y_1)\times\dotsb\times\h^\nfr_S(Y_k)
	\]
	is an $\Lhtp$-equivalence.
	In particular, for every $Y\in\Sch_S$, the $\Fin_\pt$-object $\Lhtp\h^\nfr_S(Y)$ is an $\Einfty$-object.
\end{prop}

\begin{proof}
	This is obvious if $k=0$, so we may assume $k=2$. 
	If $i\colon Z\hook \A^n_X$ is a closed immersion and $a\colon \A^n_X\to\A^n_X$ is an automorphism, then there is an induced isomorphism $\sN_{a\circ i}\simeq \sN_i$. This defines an action of $\Aut_X(\A^n_X)$ on $\Corr^{\nfr,n}_S(X,Y)$, which is natural in $X$ and $Y$. In particular, we have an action of $GL_{n,S}$ on the functor $\h^{\nfr,n}_S(-)$. 
	
	Let $\alpha$ be the given map, which is the colimit of the maps
	\[
	\alpha_n\colon \h^{\nfr,n}_S(Y_1\coprod Y_2) \to \h^{\nfr,n}_S(Y_1)\times\h^{\nfr,n}_S(Y_2).
	\] 
	If $i_1\colon Z_1\hook \A^n_X$ and $i_2\colon Z_2\hook \A^n_X$ are closed immersions with trivialized conormal sheaves, then $(0,i_1)+(1,i_2)\colon Z_1\coprod Z_2\hook \A^1 \times \A^{n}_X$ is a closed immersion with trivialized conormal sheaf. This defines maps
	\[
	\beta_n\colon \h^{\nfr,n}_S(Y_1)\times\h^{\nfr,n}_S(Y_2) \to \h^{\nfr,n+1}_S(Y_1\coprod Y_2).
	\]
	If $n$ is even, we claim that $\Lhtp$ of the diagram
	\begin{equation*}
  \begin{tikzcd}
    \h^{\nfr,n}_S(Y_1\coprod Y_2) \ar{r}{\sigma}\ar{d}{\alpha_n}
      & \h^{\nfr,n+1}_S(Y_1\coprod Y_2) \ar{d}{\alpha_{n+1}}
    \\
    \h^{\nfr,n}_S(Y_1)\times\h^{\nfr,n}_S(Y_2) \ar{r}{\sigma\times\sigma} \ar{ur}{\beta_n}
      & \h^{\nfr,n+1}_S(Y_1)\times\h^{\nfr,n+1}_S(Y_2)
  \end{tikzcd}
	\end{equation*}
	can be made commutative in such a way that the composite homotopy is equivalent to the identity homotopy.
	This implies that the maps $\Lhtp\beta_{n}$, with $n$ even, fit in a sequence whose colimit is an inverse to $\Lhtp\alpha$.
	
	For $\epsilon\in\{0,1\}$, let $\rho_\epsilon\colon \h^{\nfr,n}_S(-) \to \h^{\nfr,n+1}_S(-)$ be the obvious map induced by the closed immersion $\A^n\hook \A^{n+1}$, $x\mapsto (\epsilon,x)$. Then $\sigma$ and $\rho_0$ differ by the action of a cyclic permutation matrix of order $n+1$ in $GL_{n+1}(\Z)$. Since $n$ is even, that matrix is $\A^1$-homotopic to the identity matrix, and hence we obtain an $\A^1$-homotopy $\lambda$ from $\sigma$ to $\rho_0$. 
	Define 
	\[
	\mu\colon\h^{\nfr,n}_S(Y_1\coprod Y_2)\to \h^{\nfr,n+1}_S(Y_1\coprod Y_2)^{\A^1}
	\]
	by sending the closed immersion $i_1+i_2\colon Z_1\coprod Z_2\hook \A^n_X$ with trivializations $\tau_1\colon\sO^n_{Z_1}\simeq\sN_{i_1}$ and $\tau_2\colon\sO^n_{Z_2}\simeq\sN_{i_2}$ to
	\[
	((\id, 0)\times i_1) + (\delta\times i_2)\colon (\A^1\times Z_1)\coprod (\A^1\times Z_2)\hook \A^1 \times \A^1 \times \A^n_X\simeq\A^1 \times\A^{n+1}_X
	\]
	with conormal sheaf trivialized by $\tau_1$, $\tau_2$, and
	\begin{align*}
		\sO_{\A^1}&\simeq \sN_{(\id,0)},\quad 1\mapsto 1\otimes t,\\
		\sO_{\A^1}&\simeq \sN_\delta,\quad 1\mapsto 1\otimes t -t\otimes 1.
	\end{align*}
	Then $\mu$ is an $\A^1$-homotopy from $\rho_0$ to $\beta_n \alpha_n$. 
	We thus obtain a sequence of $\A^1$-homotopies
	\[
	\sigma \stackrel{\lambda}\rightsquigarrow \rho_0 \stackrel\mu\rightsquigarrow \beta_n\alpha_n.
	\]
	
	Note that $\alpha_{n+1}\beta_n=\rho_0\times \rho_1$.
	Let $\nu_1\colon \h^{\nfr,n}_S(Y_1)\to \h^{\nfr,n+1}_S(Y_1)^{\A^1}$ be the constant $\A^1$-homotopy at $\rho_0$.
	Define
		\[
		\nu_2\colon \h^{\nfr,n}_S(Y_2)\to \h^{\nfr,n+1}_S(Y_2)^{\A^1}
		\]
		by sending the closed immersion $i_2\colon Z_2\hook \A^n_X$ with trivialization $\tau_2\colon\sO^n_{Z_2}\simeq\sN_{i_2}$ to
		\[
		\delta\times i_2\colon \A^1\times Z_2\hook \A^1 \times \A^1 \times \A^n_X\simeq\A^1 \times\A^{n+1}_X
		\]
		with conormal sheaf trivialized by $\tau_2$ and
		\[
			\sO_{\A^1}\simeq \sN_\delta,\quad 1\mapsto 1\otimes t -t\otimes 1.
		\]
		Then $\nu_2$ is an $\A^1$-homotopy from $\rho_0$ to $\rho_1$. We thus obtain a sequence of $\A^1$-homotopies
		\[
		\sigma\times\sigma \stackrel{\lambda\times\lambda}\rightsquigarrow \rho_0\times\rho_0 \stackrel{\nu_1\times\nu_2}\rightsquigarrow \rho_0\times\rho_1= \alpha_{n+1}\beta_n.
		\]
		We then observe that
		\[
		\alpha_{n+1}^{\A^1} \circ \lambda = (\lambda\times\lambda) \circ \alpha_n
		\quad\text{and}\quad
		\alpha_{n+1}^{\A^1} \circ \mu = (\nu_1\times\nu_2) \circ \alpha_n,
		\]
		which concludes the proof.
\end{proof}

\begin{rem}
	The proof of \cite[Theorem 6.4]{garkusha2014framed} only shows that the map
	\[
	\Lhtp\h^\efr_S(Y_1\coprod Y_2)(X) \to \Lhtp\h^\efr_S(Y_1)(X)\times\Lhtp\h^\efr_S(Y_2)(X)
	\]
	induces an isomorphism on $\pi_0\Maps(K,-)$ for every finite space $K$, which is not enough to conclude that it is an equivalence. 
	The proof of Proposition~\ref{prop:nfr-E_infinity} shows how to repair the argument.
\end{rem}

\sssec{} \label{sssec:einfty-nfr} Just as in the case of equationally framed correspondences in~\sssecref{sssec:hefr-for-arbitrary}, we can formally extend $\h^{\nfr,n}_S$ and $\h^{\nfr}_S$ to sifted-colimit-preserving functors
\[
\h^{\nfr,n}_S,\h^{\nfr}_S\colon \Pre_\Sigma(\Sch_S)_*\to \Pre_\Sigma(\Sch_S)_*.
\]
 It then follows from Proposition~\ref{prop:nfr-E_infinity} that, for any $\Fcal_1, \dotsc, \Fcal_k\in \Pre_\Sigma(\Sch_S)_*$, the canonical map 
\[\h^{\nfr}_S(\Fcal_1 \vee \cdots \vee \Fcal_k) \rightarrow \h^{\nfr}_S(\Fcal_1) \times \cdots \times \h^{\nfr}_S(\Fcal_k)
\]
is an $\Lhtp$-equivalence. Therefore, the functor $\Lhtp\h^\nfr_S$ lifts uniquely to
\[
\Lhtp\h^\nfr_S\colon\Pre_\Sigma(\Sch_S)_*\to \CMon(\Pre_\Sigma(\Sch_S)),
\]
and the map~\eqref{einftymap1} induces a natural transformation
\[
\Lhtp\h^\efr_S\to \Lhtp\h^\nfr_S\colon\Pre_\Sigma(\Sch_S)_*\to \CMon(\Pre_\Sigma(\Sch_S)).
\]

\begin{rem}
	It follows from Proposition~\ref{prop:nfr-descent}(iii) that the functors $\h^{\nfr,n}_S,\h^{\nfr}_S\colon \Pre_\Sigma(\Sch_S)_*\to \Pre_\Sigma(\Sch_S)_*$ preserve Nisnevich-local and étale-local equivalences (cf.\ Remark~\ref{rem:desc}).
\end{rem}

\sssec{}\label{sssec:nfr-module-structure}
If $X,Y\in\Sch_S$ and $n\geq 0$, there is a pointed map
\[
X_+\wedge \h^{\nfr,n}_S(Y) \to \h^{\nfr,n}_S(X\times_S Y)
\]
sending a map $Z\to X$ and a normally framed correspondence $Z\leftarrow W\rightarrow Y$ to the normally framed correspondence $Z\leftarrow W\to X\times_SY$, where the framing is unchanged and the second map is the composite $W\hookrightarrow Z\times_SW\to X\times_SY$. As in \sssecref{sssec:efr-module-structure}, this construction defines right-lax $\Pre_\Sigma(\Sch_S)_\pt$-linear structures on the functors $\h^{\nfr,n}_S,\h^\nfr_S\colon \Pre_\Sigma(\Sch_S)_\pt\to \Pre_\Sigma(\Sch_S)_\pt$.
Moreover, it is clear that the natural transformations $\h^{\efr,n}_S\to \h^{\nfr,n}_S$ and $\h^{\efr}_S\to \h^{\nfr}_S$ are $\Sch_{S+}$-linear, whence $\Pre_\Sigma(\Sch_S)_\pt$-linear. 
In particular, for every $\sF,\sG\in\Pre_\Sigma(\Sch_S)_\pt$, we have a commutative square
\[
\begin{tikzcd}
	\sF\wedge \h^{\efr}_S(\sG) \ar{r} \ar{d} & \h^{\efr}_S(\sF\wedge\sG) \ar{d} \\
	\sF\wedge \h^{\nfr}_S(\sG) \ar{r} & \h^{\nfr}_S(\sF\wedge\sG).
\end{tikzcd}
\]

\sssec{}\label{sssec:nfr-unit}
As in \sssecref{sssec:efr-unit}, we also have a $\Pre_\Sigma(\Sch_S)_*$-linear transformation $\id\to \h^{\nfr}_S\colon \Pre_\Sigma(\Sch_S)_*\to \Pre_\Sigma(\Sch_S)_*$, coming from the fact that $\h^{\nfr,0}_S\simeq\id$.

\sssec{}
We now prove that equationally framed correspondences and normally framed correspondences are equivalent from the point of view of motivic homotopy theory. The key step is the following result:

\begin{prop}\label{prop:extension}
Let $X,Y\in\Sch_{S}$, let $X_0\subset X$ be a closed subscheme, and let $n\geq 0$. 
Suppose that $X$ is affine and that $Y$ admits an étale map to an affine bundle over $S$. Then the map
\[
\Corr^{\efr,n}_S(X,Y) \to \Corr^{\efr,n}_S(X_0,Y) \times_{\Corr^{\nfr,n}_S(X_0,Y)} \Corr^{\nfr,n}_S(X,Y)
\]
induced by~\eqref{eqn:vfr-to-nfr} is surjective.
\end{prop}

\begin{proof}
	An element in the right-hand side consists of:
	\begin{itemize}
		\item a span $X\xleftarrow{f} Z\xrightarrow{h} Y$ with $f$ finite syntomic;
		\item an embedding $i\colon Z\hook \A^n_X$ over $X$ with a trivialization $\tau\colon \sO_Z^n \simeq \sN_i$;
		\item an equational framing of the induced span $X_0\xleftarrow{f_0} Z_0\xrightarrow{h_0} Y$:
		\[
		Z_0\stackrel{i_0}\hooklong \A^n_{X_0} \leftarrow U_0 \xrightarrow{(\varphi_0,g_0)} \A^n\times Y,
		\]
		where $i_0$ is restriction of $i$, $U_0\to \A^n_{X_0}$ is an affine étale neighborhood of $Z_0$, $g_0$ extends $h_0$, and $Z_0=\varphi_0^{-1}(0)$, such that the trivialization $\sO_{Z_0}^n\simeq \sN_{i_0}$ induced by $\tau$ coincides with that induced by $\varphi_0$.
	\end{itemize}
	The goal is to construct a suitable equational framing of $X\xleftarrow{f} Z\xrightarrow{h} Y$.
	Using Lemma~\ref{lem:lifting-neighborhoods}, we can lift the étale neighborhood $U_0$ of $Z_0$ in $\A^n_{X_0}$ to an étale neighborhood $U$ of $Z$ in $\A^n_X$. 
	Refining $U_0$ if necessary, we can assume that $U$ is affine (by Lemma~\ref{lem:henselization}(ii)).
	
	We first construct a simultaneous extension $g\colon U\to Y$ of $g_0\colon U_0\to Y$ and $h\colon Z\to Y$. Suppose first that $Y$ is an affine bundle over $S$. Since $U$ is affine, $U\times_SY\to U$ is a vector bundle over $U$. It follows from \cite[Lemma 2.5]{HoyoisCdh} that the restriction map
		\[
		\Maps_S(U,Y) \to \Maps_S(U_0,Y) \times_{\Maps_S(Z_0,Y)} \Maps_S(Z,Y) \simeq \Maps_S(U_0\coprod_{Z_0}Z,Y)
		\]
		is surjective, so the desired extension exists. In general, let $p\colon Y\to A$ be an étale map where $A$ is an affine bundle over $S$. By the previous case, there exists an $S$-morphism $U\to A$ extending $p\circ g_0$ and $p\circ h$. Then the étale map $U\times_{A}Y\to U$ has a section over $U_0\coprod_{Z_0}Z$, so there exists an open subset $U'\subset U\times_{A}Y$ that is an étale neighborhood of $U_0\coprod_{Z_0}Z$ in $U$. We can therefore replace $U$ by $U'$, and the projection $U'\to Y$ gives the desired extension.
	
	It remains to construct an extension $\phi\colon U\to\A^n$ of $\phi_0$ such that $\phi^{-1}(0)=Z$ and such that the induced trivialization $\sO_Z^n\simeq \sN_i$ equals $\tau$.
	Let $\sI\subset \sO(U)$ and $\sI_0\subset \sO(U_0)$ be the ideals defining $Z$ and $Z_0$. We observe that the map of $\sO(U)$-modules
	\begin{equation}\label{eq:surjective}
	\sI \to \sI_0 \times_{\sI_0/\sI_0^2} \sI/\sI^2 \simeq \sI_0\times_{\sN_{i_0}} \sN_i
	\end{equation}
	is surjective, since we have a morphism of short exact sequences
	\[
	\begin{tikzcd}
		0 \ar{r} & \sI^2 \ar{r} \ar{d} & \sI \ar{r} \ar{d} & \sI/\sI^2 \ar{r} \ar[equal]{d} & 0 \\
		0 \ar{r} & \sI_0^2 \ar{r} & \sI_0 \times_{\sI_0/\sI_0^2} \sI/\sI^2 \ar{r} & \sI/\sI^2 \ar{r} & 0
	\end{tikzcd}
	\]
	 and $\sI^2\to\sI_0^2$ is surjective. The morphism $\varphi_0$ and the trivialization $\tau$ define an $n$-tuple of elements in the right-hand side of~\eqref{eq:surjective}. By lifting these elements arbitrarily, we obtain a morphism $\varphi\colon U\to\A^n$ simultaneously lifting $\varphi_0$ and the trivialization $\tau$ of $\sN_i$.
	Moreover, by Nakayama's lemma, there exists a function $a$ on $U$ such that $\phi^{-1}(0)\cap U_a=Z$. Replacing $U$ by $U_a$ concludes the proof.
\end{proof}

The following corollary was also independently obtained by A.~Neshitov (see \cite[Corollary 4.9]{MWvsFramed} for a related statement at the level of connected components).

\begin{cor}\label{cor:vfr-vs-nfr}
	Suppose that $Y\in\Sm_{S}$ is a finite sum of schemes admitting étale maps to affine bundles over $S$. Then the $\Einfty$-map
	\[
	\Lhtp \h^\efr_S(Y) \to \Lhtp \h^\nfr_S(Y)
	\]
	induced by~\eqref{einftymap1}
	is an equivalence on affines. In particular, it induces an equivalence of presheaves of $\Einfty$-spaces
	\[
	L_\zar\Lhtp \h^\efr_S(Y) \simeq L_\zar\Lhtp \h^\nfr_S(Y).
	\]
\end{cor}

\begin{proof}
	By Remark~\ref{rem:GPadditivity} and Proposition~\ref{prop:nfr-E_infinity}, we can assume that $Y$ admits an étale map to an affine bundle over $S$.
	By Proposition~\ref{prop:extension}, for every $n\geq 0$, the map
	\[
	\h^\efr_S(Y)^{\A^n} \to \h^\efr_S(Y)^{\partial\A^n} \times_{\h^\nfr_S(Y)^{\partial\A^n}} \h^\nfr_S(Y)^{\A^n}
	\]
	is surjective on affines. 
	By Propositions~\ref{prop:fr-descent}(ii) and~\ref{prop:nfr-descent}(ii), both $\h^\efr_S(Y)$ and $\h^\nfr_S(Y)$ satisfy closed gluing. It follows from Lemma~\ref{lem:sset-as-schemes} that the map
	\[
	\h^\efr_S(Y)^{\A^\bullet} \to \h^\nfr_S(Y)^{\A^\bullet}
	\]
	is a trivial Kan fibration of simplicial sets when evaluated on any affine scheme. The equivalence of $\Einfty$-spaces follows because $L_\zar$ preserves finite products.
\end{proof}

\begin{cor}\label{cor:vfr-vs-nfr2}
	Let $\sF\in\Pre_\Sigma(\Sm_S)_*$. Then the map
	\[
	\Lmot \h^\efr_S(\sF) \to \Lmot \h^\nfr_S(\sF)
	\]
	induced by~\eqref{einftymap1} is an equivalence of presheaves of $\Einfty$-spaces.
\end{cor}

\begin{proof}
	We can assume that $\sF=Y_+$ for some $Y\in\Sm_S$.
	Write $Y$ as a filtered union of quasi-compact open subschemes $Y_i$. Then the map $\colim_i \h^\efr_S(Y_i)\to \h^\efr_S(Y)$ is an isomorphism on qcqs schemes (by Lemma~\ref{lem:henselization}(i)), and similarly for $\h^\nfr_S$. We can therefore assume that $Y$ is quasi-compact. Let $\{U_1,\dotsc,U_n\}$ be a finite open covering of $Y$ where each $U_i$ admits an étale map to an affine bundle over $S$, and let $U=U_1\coprod\dotsb\coprod U_n$. Then every iterated fiber product $U\times_Y\dotsb\times_YU$ is a finite sum of schemes admitting étale maps to affine bundles over $S$. We conclude using Propositions \ref{prop:fr-descent}(iii) and~\ref{prop:nfr-descent}(iii) and Corollary~\ref{cor:vfr-vs-nfr}. The equivalence of $\Einfty$-spaces follows because $\Lmot$ preserves finite products.
\end{proof}

\sssec{}
\label{sssec:Emb^fr}

We give a slightly different description of the set $\Corr^\nfr_S(X,Y)$ that will be useful later.
Let $\FSyn\subset \Fun(\Delta^1,\Sch)$ be the subcategory whose objects are finite syntomic morphisms and whose morphisms are Cartesian squares. For $Z\to X$ finite syntomic, let $\Emb_X(Z,\A^n_X)$ denote the set of closed immersions $Z\hook \A^n_X$ over $X$, and let
\[
\Emb_X(Z,\A^\infty_X) = \colim_{n\to\infty} \Emb_X(Z,\A^n_X).
\]
Note that these sets vary functorially with $Z\to X$:
\[
\Emb(-,\A^n),\; \Emb(-,\A^\infty)\colon \FSyn^\op\to\Set.
\]

Let $\sVect_0(Z)$ denote the groupoid of stable vector bundles of rank $0$ over $Z$. By definition, this is the colimit of the sequence
\[
\Vect_0(Z) \xrightarrow{\oplus\sO_Z} \Vect_1(Z) \xrightarrow{\oplus\sO_Z} \ldots.
\]
Then we have natural maps
\[
\Emb_X(Z,\A^n_X) \to \Vect_n(Z), \quad (i\colon Z\hook \A^n_X) \mapsto \sN_i,
\]
inducing
\[
\Emb_X(Z,\A^\infty_X) \to \sVect_0(Z).
\]
in the colimit.

Let $\Emb^\fr_X(Z,\A^n_X)$ denote the set of closed immersions $i\colon Z\hook \A^n_X$ equipped with a trivialization of their conormal sheaves, and let
\[
\Emb^\fr_X(Z,\A^\infty_X) = \colim_{n\to\infty} \Emb^\fr_X(Z,\A^n_X).
\]
In other words, we have Cartesian squares
\[
    \begin{tikzcd}
      \Emb^\fr_X(Z,\A^n_X) \ar{r} \ar{d}
        & \Emb_X(Z,\A^n_X) \ar{d}
      \\
      * \ar{r}{\sO_Z^n}
        & \Vect_n(Z),
    \end{tikzcd}
	 \qquad
    \begin{tikzcd}
      \Emb^\fr_X(Z,\A^\infty_X) \ar{r} \ar{d}
        & \Emb_X(Z,\A^\infty_X) \ar{d}
      \\
      * \ar{r}{0}
        & \sVect_0(Z),
    \end{tikzcd}
\]
which are natural in $(Z\to X)\in \FSyn^\op$.

The set of normally framed correspondences from $X$ to $Y$ can now be written as
\[
\Corr^\nfr_S(X,Y) = \colim_{X\stackrel f\leftarrow Z\to Y} \Emb^\fr_X(Z,\A^\infty_X),
\]
where the colimit is indexed by the groupoid of spans $X\xleftarrow{f} Z\to Y$ with $f$ finite syntomic.


\ssec{Tangentially framed correspondences}
\label{ssec:dfr}

\sssec{}
\label{sssec:cotangent-complex}

Recall that every morphism of schemes $f\colon Y\to X$ admits a \emph{cotangent complex} $\sL_f\in\QCoh(Y)$. The functoriality of the construction $f\mapsto\sL_f$ is described by the diagram
\begin{equation*}
  \begin{tikzcd}
     & \QCoh \ar{d}{p} \\
    \Fun(\Delta^1,\Sch)^\op \ar{r}{s} \ar{ur}{\sL} &  \Sch^\op,
  \end{tikzcd}
\end{equation*}
where $p$ is the coCartesian fibration classified by $X\mapsto\QCoh(X)$ and $s$ is the source functor.
Moreover, the section $\sL$ takes Tor-independent squares that are Cartesian (or more generally whose Cartesian gap is étale) to $p$-coCartesian edges.

\sssec{}
\label{sssec:cotangent-lci}

If $f\colon Y\to X$ is an lci morphism, its cotangent complex $\sL_f$ is perfect and concentrated in degrees $0$ and $1$.
In fact, given a factorization
\[
Y\stackrel i\hook V \stackrel p\to X
\]
where $i$ is a closed immersion and $p$ is smooth (which always exists locally on $Y$), we have a canonical equivalence
\[
\sL_f \simeq (\sN_i\to i^*(\Omega_p)).
\]

\sssec{}
\label{sssec:difffr}

In particular, if $f\colon Z\to X$ is finite syntomic and $i\colon Z\hook \A^n_X$ is a closed immersion over $X$, we have a canonical equivalence
\[
\sL_f \simeq (\sN_i\to \sO^n_Z).
\]
In this case, $\sL_f\simeq 0$ in $\Perf(Z)$ if and only if $f$ is finite étale. However, if we are given a trivialization $\tau\colon \sO_Z^n\simeq\sN_i$ as in the definition of a normally framed correspondence, we deduce that $\sL_f\simeq 0$ in the $K$-theory $\infty$-groupoid $K(Z)$. This motivates the following definition.

\begin{defn}\label{defn:diff-framed-corr}
	Let $X,Y\in\Sch_S$. A \emph{tangentially framed correspondence}, or simply a \emph{framed correspondence}, from $X$ to $Y$ over $S$ consists of the following data:

\noindent{\em(i)}
A span 
\begin{equation*}
  \begin{tikzcd}
     & Z \ar[swap]{ld}{f}\ar{rd}{h} & \\
    X &   & Y
  \end{tikzcd}
\end{equation*}
over $S$, where $f$ is finite syntomic.

\noindent{\em(ii)}
A trivialization $\tau\colon 0\simeq \sL_f$ of the cotangent complex of $f$ in the $K$-theory $\infty$-groupoid $K(Z)$.
\end{defn}

We denote by $\Corr^\dfr_S(X,Y)$ the $\infty$-groupoid of tangentially framed correspondences:
\[
\Corr^\dfr_S(X,Y)=\colim_{X\stackrel f\leftarrow Z\to Y} \Maps_{K(Z)}(0,\sL_f),
\]
the colimit being taken over the groupoid of spans $X\xleftarrow{f} Z\to Y$ with $f$ finite syntomic.

Explicitly, a $1$-morphism $(Z,f, h, \tau) \rightarrow (Z', f', h', \tau')$ between tangentially framed correspondences consists of:
\begin{itemize}
	\item an isomorphism $t\colon Z \to Z'$ making the following diagram commute:
\begin{equation*}
\begin{tikzcd}
 & Z \ar[swap]{dl}{f} \ar{dd}{t} \ar{dr}{h} & \\
 X  & & Y; \\
 & Z' \ar{ul}{f'} \ar[swap]{ur}{h'} & 
\end{tikzcd}
\end{equation*}
\item a path-homotopy $\beta\colon dt \circ\tau\to t^*(\tau')$ in $K(Z)$, where $dt\colon \sL_{f}\simeq t^*\sL_{f'}$ is the canonical equivalence.
\end{itemize}

By the functoriality of the cotangent complex recalled in \sssecref{sssec:cotangent-complex}, we obtain a functor
\[
\Corr^\dfr_S(-,-)\colon \Sch_S^\op\times\Sch_S\to \Spc.
\]
We will denote by $\h^\dfr_S(Y)$ the presheaf $X\mapsto \Corr^\dfr_S(X,Y)$ on $\Sch_S$. 

\begin{rem}
	In Definition~\ref{defn:diff-framed-corr}(i), we can replace ``syntomic'' by ``lci''. Indeed, (ii) implies that $f$ is of virtual relative dimension $0$, hence flat by Proposition~\ref{prop:lci-flat}(i).
\end{rem}

\begin{rem}
	If $X$ is affine, a span $X\xleftarrow{f} Z\to Y$ with $f$ finite lci admits a tangential framing if and only if $f$ is a relative global complete intersection. This follows from Proposition~\ref{prop:lci-flat}(ii).
\end{rem}

\begin{prop}\label{prop:corrfr-descent} 
Let $Y\in\Sch_S$.

\noindent{\em(i)}
$\h^{\dfr}_S(Y)$ is a Nisnevich sheaf on qcqs schemes.

\noindent{\em(ii)}
Let $R\hook Y$ be a Nisnevich (resp.\ étale) covering sieve generated by a single map. Then $\h^\dfr_S(R)\to \h^\dfr_S(Y)$ is a Nisnevich (resp.\ étale) equivalence.
\end{prop}

\begin{proof}
	(i) follows from the fact that $K$-theory is a Nisnevich sheaf on qcqs schemes. The proof of (ii) is similar to that of Proposition~\ref{prop:fr-descent}(iii).
\end{proof}

\sssec{}

We will now construct a natural map
\begin{equation}\label{eq:nfr-to-fr}
\Corr^\nfr_S(X,Y) \to \Corr^\dfr_S(X,Y),
\end{equation}
following the idea discussed in \sssecref{sssec:difffr}.

First, note that there is a natural map
\[
\sVect_0(-) \to K(-), \quad (\sE\in\Vect_n)\mapsto \sO^n-\sE.
\]
If $f\colon Z\to X$ is finite syntomic, it follows from \sssecref{sssec:difffr} that there is a canonical commutative square
\[
    \begin{tikzcd}
      \Emb_X(Z,\A^\infty_X) \ar{r} \ar{d}
        & * \ar{d}{\sL_f}
      \\
      \sVect_0(Z) \ar{r}
        & K(Z).
    \end{tikzcd}
\]
Combining this with the second Cartesian square of \sssecref{sssec:Emb^fr}, we obtain a commutative square
\[
    \begin{tikzcd}
      \Emb^\fr_X(Z,\A^\infty_X) \ar{r} \ar{d}
        & * \ar{d}{\sL_f}
      \\
      * \ar{r}{0}
        & K(Z),
    \end{tikzcd}
\]
whence a natural map
\begin{equation}\label{eqn:Emb-to-Maps}
\Emb^\fr_X(Z,\A^\infty_X) \to \Maps_{K(Z)}(0,\sL_f).
\end{equation}
Finally, taking the colimit of both sides over the groupoid of spans $X\xleftarrow{f}Z\to Y$ with $f$ finite syntomic, we obtain the desired map~\eqref{eq:nfr-to-fr}.

\sssec{}
\label{sssec:CorrF}

Both functors $\Corr^\nfr_S(-,-)$ and $\Corr^\dfr_S(-,-)$ are special cases of the following construction.
Recall that $\FSyn\subset \Fun(\Delta^1,\Sch)$ is the category whose objects are finite syntomic morphisms and whose morphisms are Cartesian squares. Define
\[
\Fun(\FSyn^\op,\Spc) \to \Fun(\Sch_S^\op\times\Sch_S,\Spc), \quad F\mapsto \Corr^F_S(-,-),
\]
by the formula
\[
\Corr^F_S(X,Y) = \colim_{X\stackrel f\leftarrow Z\to Y} F(f),
\]
where the colimit is indexed by the groupoid $\Corr^\fsyn_S(X,Y)$ of spans $X\xleftarrow{f} Z\to Y$ with $f$ finite syntomic.
We denote by $\h^F_S(Y)$ the presheaf $X\mapsto \Corr^F_S(X,Y)$ on $\Sch_S$.
When $F$ is the functor $\Emb^\fr(-,\A^\infty)$, we recover $\Corr^\nfr_S$, and when $F$ is the functor $\Maps_{K}(0,\sL)$, we recover $\Corr_S^\dfr$.
 
Given $F\colon \FSyn^\op\to \Spc$, we define a new functor $\Lhtp F\colon \FSyn^\op\to\Spc$ by
\[
(\Lhtp F)(Z\to X) = \colim_{n\in\Delta^\op} F(Z\times\A^n\to X\times \A^n).
\]

\begin{lem}\label{lem:Sing-on-frames}
	Let $F\colon \FSyn^\op\to \Spc$ be a functor and $Y\in\Sch_S$. Then the natural transformation $F\to \Lhtp F$ induces an equivalence
	\[
	\Lhtp \h^F_S(Y) \simeq \Lhtp \h^{\Lhtp F}_S(Y).
	\]
\end{lem}

\begin{proof}
	Note that $F\mapsto \Corr^F_S(-,-)$ preserves colimits.
	There is a canonical map $\h^{\Lhtp F}_S(Y)\to \Lhtp \h^F_S(Y)$ obtained as the colimit over $n\in\Delta^\op$ of the maps
	\[
	\psi_n^F\colon\Corr_S^{F(-\times\A^n)}(X,Y) \to \Corr_S^F(X\times\A^n,Y),\quad (Z\to X)\mapsto (Z\times\A^n\to X\times\A^n).
	\]
	It is clear that the composition
\begin{equation*} \label{localization}
	\h^F_S(Y) \to \h^{\Lhtp F}_S(Y) \to \Lhtp \h^F_S(Y)
\end{equation*}
	is the localization map. It remains to show that the other composition
\begin{equation} \label{othercomposite}
	\h^{\Lhtp F}_S(Y) \to \Lhtp \h^F_S(Y) \to \Lhtp\h^{\Lhtp F}_S(Y)
\end{equation}
	is also the localization map. Unpacking the definitions, the composite~\eqref{othercomposite} is the colimit over $n$ and $m$ of
	\[
	\Corr_S^{F(-\times\A^n)}(X,Y) \xrightarrow{\psi_n^F} \Corr_S^F(X\times\A^n,Y) \to \Corr_S^{F(-\times\A^m)}(X\times\A^n,Y).
	\]
	 Since $\Delta^\op$ is sifted, we can take $m=n$. Then we can rewrite this composition as
	\[
	\Corr_S^{F(-\times\A^n)}(X,Y) \xrightarrow{\pi_1^*} \Corr_S^{F(-\times\A^n\times\A^n)}(X,Y) \xrightarrow{\psi_n^{F(-\times\A^n)}} \Corr_S^{F(-\times\A^n)}(X\times\A^n,Y).
	\]
	On the other hand, the localization map $\h^{\Lhtp F}_S(Y) \to \Lhtp\h^{\Lhtp F}_S(Y)$ is the colimit over $n$ of
	\[
	\Corr_S^{F(-\times\A^n)}(X,Y) \xrightarrow{\pi_2^*} \Corr_S^{F(-\times\A^n\times\A^n)}(X,Y) \xrightarrow{\psi_n^{F(-\times\A^n)}} \Corr_S^{F(-\times\A^n)}(X\times\A^n,Y).
	\]
	Since $F\mapsto \Corr^F_S(-,-)$ preserves colimits, it will suffice to show that
	\[
	\pi_1^*,\pi_2^*\colon F(-\times\A^\bullet)\rightrightarrows F(-\times\A^\bullet\times\A^\bullet)
	\]
	induce the same map in the colimit over $\Delta^\op$. But both induce equivalences in the colimit, and they have a common retraction
	\[
	\delta^*\colon F(-\times\A^\bullet\times\A^\bullet)\to F(-\times\A^\bullet),
	\]
	where $\delta\colon \A^\bullet\to \A^\bullet\times \A^\bullet$ is degreewise the diagonal map.
	This concludes the proof.
\end{proof}

\sssec{}
\label{sssec:CorrF-E_infinity}

Let $\FSyn'\subset\Fun(\Delta^1,\Sch)$ be the subcategory whose objects are finite syntomic morphisms and whose morphisms are squares 
\[
\begin{tikzcd}
	Z' \ar{r} \ar{d} & Z \ar{d} \\
	X' \ar{r} & X
\end{tikzcd}
\]
whose Cartesian gap $Z'\to Z\times_XX'$ is a clopen embedding. For a functor $F\colon (\FSyn')^{\op}\to \Spc$, the construction of \sssecref{sssec:CorrF} can be extended to
\[
\Corr^F_S(-,-)\colon \Sch_{S}^\op \times \Sch_{S+} \to \Spc.
\]
In particular, $\h^F_S(Y)$ becomes a $\Fin_\pt$-object (see \sssecref{sssec:nfr-E_infinity}).

For $F=\Emb^\fr(-,\A^\infty)$, this partially recovers the extension of $\Corr^\nfr_S(-,-)$ discussed in \sssecref{sssec:Corrnfr-extended-functoriality}. 
The functor $\Maps_K(0,\sL)\colon \FSyn^\op\to \Spc$ also extends to $(\FSyn')^\op$ by the functoriality of the cotangent complex recalled in \sssecref{sssec:cotangent-complex}.
The natural transformation $\Emb^\fr(-,\A^\infty)\to\Maps_K(0,\sL)$ on $(\FSyn')^{\op}$ then shows that $\h^\nfr_S(Y)\to\h^\dfr_S(Y)$ is a map of $\Fin_*$-objects.
In general we have the following form of additivity:

\begin{lem} \label{lem:f-additivity}  
	Let $F\colon (\FSyn')^{\op} \rightarrow \Spc$ be a functor such that, for any finite collection of finite syntomic morphisms $\{Z_i\to X\}_{1\leq i\leq k}$, the map
	\[F(Z_1\coprod\cdots\coprod Z_k\to X) \rightarrow F(Z_1\to X) \times\cdots\times F(Z_k\to X)\]
	induced by the inclusions $Z_i\hookrightarrow Z_1\coprod\cdots\coprod Z_k$	is an equivalence. Then for every $Y_1, \dotsc, Y_k \in \Sch_S$, the natural map 
\[\h_S^F(Y_1 \coprod \cdots \coprod Y_k) \rightarrow \h_S^F(Y_1) \times \cdots \times \h_S^F(Y_k)\] 
is an equivalence. In particular, the $\Fin_\pt$-object $\h^F_S(Y)$ is an $\Einfty$-object.
\end{lem}

\begin{proof} 
There is an obvious equivalence of groupoids 
\[\prod_i \Corr^\fsyn_S(X,Y_i)\simeq \Corr^\fsyn_S(X,Y_1\coprod\cdots \coprod Y_k),\] 
and the claim follows by distributivity of finite products over colimits in $\Spc$.
\end{proof}

This lemma clearly applies to $F=\Maps_K(\sL,0)$, and we deduce that $\h^\dfr_S(Y)$ is a presheaf of $\Einfty$-spaces. Since $\Lhtp$ preserves finite products, we deduce that $\Lhtp \h^\nfr_S(Y)\to \Lhtp\h^\dfr_S(Y)$ is an $\Einfty$-map. 

\sssec{}\label{sssec:CorrF-pointing}
Let $F\colon\FSyn^\op\to\Spc$ be a functor. We make a couple of basic observations:
\begin{itemize}
	\item If $F$ preserves finite products, then the presheaf $\h^F_S(Y)$ belongs to $\Pre_\Sigma(\Sch_S)$, and we have
	\[
	\h^F_S\colon \Sch_S \to \Pre_\Sigma(\Sch_S).
	\]
	\item If $F(\emptyset\to X)$ is contractible for every $X$, then each $\infty$-groupoid $\Corr^F_S(X,Y)$ acquires a base point given by the empty correspondence, and the functor $\h^F_S$ can be promoted to
\[
\h^F_S\colon \Sch_S \to \Pre(\Sch_S)_\pt.
\]
\end{itemize}
If $F$ satisfies both of these conditions and moreover extends to $\FSyn'$ as in \sssecref{sssec:CorrF-E_infinity}, we obtain a functor
\[
\h^F_S\colon \Sch_{S+} \to \Pre_\Sigma(\Sch_S)_\pt,
\]
which we can formally extend to a sifted-colimit-preserving functor
\[
\h^F_S\colon \Pre_\Sigma(\Sch_S)_\pt \to \Pre_\Sigma(\Sch_S)_\pt.
\]
This construction recovers the extensions of $\h^{\nfr,n}_S$ and $\h^{\nfr}_S$ discussed in \sssecref{sssec:einfty-nfr}. 
We can also apply this to $F=\Maps_{K}(0,\sL)$ and to the natural transformation $\Emb^\fr(-,\A^\infty)\to\Maps_K(0,\sL)$, so that we obtain a natural tranformation
\begin{equation}\label{eqn:hnfr-to-hfr-extended}
\h^\nfr_S\to \h^\fr_S\colon \Pre_\Sigma(\Sch_S)_\pt \to \Pre_\Sigma(\Sch_S)_\pt.
\end{equation}

\begin{rem}
	It follows from Proposition~\ref{prop:corrfr-descent}(ii) that the functor $\h^{\fr}_S\colon \Pre_\Sigma(\Sch_S)_*\to \Pre_\Sigma(\Sch_S)_*$ preserve Nisnevich-local and étale-local equivalences (cf.\ Remark~\ref{rem:desc}).
\end{rem}

\sssec{}\label{sssec:dfr-module-structure}
Let $F\colon (\FSyn')^\op\to\Spc$ be a functor satisfying the conditions of \sssecref{sssec:CorrF-pointing}.
For $X,Y\in\Sch_S$, there is map
\[
X_+\wedge \h^{F}_S(Y) \to \h^{F}_S(X\times_S Y)
\]
defined exactly as in \sssecref{sssec:nfr-module-structure}, natural in $X_+$ and $Y_+$. In particular, we can formally extend it to a natural transformation of functors $\Pre_\Sigma(\Sch_S)_\pt\times \Pre_\Sigma(\Sch_S)_\pt\to \Pre_\Sigma(\Sch_S)_\pt$, which is natural in $F$. 

For example, the natural transformation $\Emb^\fr(-,\A^\infty)\to\Maps_K(0,\sL)$ yields for every $\sF,\sG\in\Pre_\Sigma(\Sch_S)_\pt$ a commutative square
\[
\begin{tikzcd}
	\sF\wedge \h^{\nfr}_S(\sG) \ar{r} \ar{d} & \h^{\nfr}_S(\sF\wedge\sG) \ar{d} \\
	\sF\wedge \h^{\fr}_S(\sG) \ar{r} & \h^{\fr}_S(\sF\wedge\sG),
\end{tikzcd}
\]
where the vertical maps are given by~\eqref{eqn:hnfr-to-hfr-extended}. The natural transformation $\id\to\h^{\nfr}_S$ from \sssecref{sssec:nfr-unit} is similarly induced by the transformation $\Emb^\fr(-,\A^0)\hook \Emb^\fr(-,\A^\infty)$. In particular, we have a commuting triangle
\[
\begin{tikzcd}
	\sF\wedge \sG \ar{dr} \ar{r} & \sF\wedge \h^{\fr}_S(\sG) \ar{d} \\
	& \h^{\fr}_S(\sF\wedge\sG).
\end{tikzcd}
\]

In fact, it is not difficult to construct a right-lax $\Pre_\Sigma(\Sch_S)_\pt$-linear structure on $\h^F_S$, but we will not need it.
For $F=\Maps_K(0,\sL)$, this structure is subsumed by the symmetric monoidal structure on the $\infty$-category of framed correspondences that we will construct in \secref{sect:infinity-category}.

\begin{lem} \label{lem:retract} 
	Let $F\colon (\FSyn')^\op\to\Spc$ be a functor satisfying the conditions of \sssecref{sssec:CorrF-pointing}. Then the functor $\h^F_S\colon \Pre_\Sigma(\Sch_S)_\pt\to \Pre_\Sigma(\Sch_S)_\pt$ preserves $\A^1$-homotopy equivalences.
\end{lem}

\begin{proof}
	It suffices to show that $\h^F_S$ preserves $\A^1$-homotopic maps. If $h\colon \A^1_+\wedge X\to Y$ is an $\A^1$-homotopy between $f$ and $g$, then the composite
	\[
	\A^1_+\wedge \h^F_S(X) \to \h^F_S(\A^1_+\wedge X) \xrightarrow{\h^F_S(h)} \h^F_S(Y)
	\]
	is an $\A^1$-homotopy between $\h^F_S(f)$ and $\h^F_S(g)$.
\end{proof}

\sssec{}
We now prove that normally framed correspondences and tangentially framed correspondences are equivalent from the point of view of motivic homotopy theory. The first step is the observation that the space of embeddings of a finite scheme in $\A^\infty$ is $\A^1$-contractible:

\begin{lem}\label{lem:contractible-Emb}
	The functor $L_{\A^1}\Emb(-,\A^\infty)\colon \FSyn^\op\to \Spc$ is contractible on affines.
\end{lem}

\begin{proof}
	Let $X$ be affine and let $f\colon Z\to X$ be finite syntomic. We will show that the simplicial set
	\[
	\Emb_{X\times\A^\bullet}(Z\times\A^\bullet,\A^\infty_{X\times\A^\bullet})
	\]
	is a contractible Kan complex by solving the lifting problem
	\[
   \begin{tikzcd}
     \partial\Delta^n \ar{r} \ar[hookrightarrow]{d}
       & \Emb_{X\times\A^\bullet}(Z\times\A^\bullet,\A^\infty_{X\times\A^\bullet}).
     \\
     \Delta^n \ar[dashed]{ur}
       &
   \end{tikzcd}
	\]
	The proof only uses that $f$ is finite.
	By Lemma~\ref{lem:sset-as-schemes}, it suffices to prove the following: for every finitely presented closed subscheme $X_0\subset X$, the pullback map
	\[
	\Emb_X(Z,\A^\infty_X) \to \Emb_{X_0}(Z_0,\A^\infty_{X_0})
	\]
	is surjective, where $Z_0=Z\times_XX_0$. Let $i_0\colon Z_0\hook \A^n_{X_0}$ be a closed immersion over $X_0$, given by $n$ functions $g_1$, \dots, $g_n$ on $Z_0$. Let $g_1'$, \dots, $g_n'$ be lifts of these functions to $Z$, and let $h_1$, \dots, $h_m$ be functions on $Z$ that generate the ideal of $Z_0$ in $Z$ as an $\sO(X)$-module (this is possible since $f$ is finite and $Z_0\subset Z$ is finitely presented). Then the $n+m$ functions $g_1'$, \dots, $g_n'$, $h_1$, \dots, $h_m$ define a closed immersion $i\colon Z\hook \A^{n+m}_X$ over $X$ lifting $i_0$.
\end{proof}

\begin{prop}\label{prop:nfr-to-dfr}
	The natural transformation
	\[
	L_{\A^1}\Emb^\fr(-,\A^\infty) \to L_{\A^1}\Maps_K(0,\sL)\colon \FSyn^\op\to\Spc
	\]
	induced by~\eqref{eqn:Emb-to-Maps} is an equivalence on affines.	
\end{prop}

\begin{proof}
	Recall that~\eqref{eqn:Emb-to-Maps} comes from a natural transformation of Cartesian squares
	\[
   \begin{tikzcd}
     \Emb^\fr_X(Z,\A^\infty_X) \ar{r} \ar{d}
       & \Emb_X(Z,\A^\infty_X) \ar{d}
     \\
     * \ar{r}{0}
       & \sVect_0(Z)
   \end{tikzcd}
	\quad\longrightarrow
   \begin{tikzcd}
     \Maps_{K(Z)}(0,\sL_f) \ar{r} \ar{d}
       & * \ar{d}{\sL_f}
     \\
     * \ar{r}{0}
       & K(Z).
   \end{tikzcd}
	\]
	We consider two cases. If $[\sL_f]\neq 0$ in $K_0(Z)$, then $\Maps_{K(Z)}(0,\sL_f)$ is empty and the result holds trivially.
	Suppose that $X$ is affine and that $[\sL_f]=0$ in $K_0(Z)$. 
	Then $\sL_f$ lives in the unit component $\tau_{\geq 1}K(Z)\subset K(Z)$.
	Since $Z$ is affine and $[\sL_f]=[\sO_Z^n]-[\sN_i]$ for any closed immersion $i\colon Z\hook \A^n_X$ over $X$, the conormal sheaf $\sN_i$ is stably free. It follows that the map $\Emb_X(Z,\A^\infty_X)\to \sVect_0(Z)$ lands in the component $\tau_{\geq 1}\sVect_0(Z)\subset \sVect_0(Z)$ of the trivial bundle. We may therefore rewrite the above Cartesian squares as follows:
	\[
   \begin{tikzcd}
     \Emb^\fr_X(Z,\A^\infty_X) \ar{r} \ar{d}
       & \Emb_X(Z,\A^\infty_X) \ar{d}
     \\
     * \ar{r}{0}
       & \tau_{\geq 1}\sVect_0(Z)
   \end{tikzcd}
	\quad\longrightarrow
   \begin{tikzcd}
     \Maps_{K(Z)}(0,\sL_f) \ar{r} \ar{d}
       & * \ar{d}{\sL_f}
     \\
     * \ar{r}{0}
       & \tau_{\geq 1}K(Z).
   \end{tikzcd}
	\]
	Note that $\tau_{\geq 1}\sVect_0\simeq BGL$. The map $\tau_{\geq 1}\sVect_0\to \tau_{\geq 1}K$ between the lower right corners is an $\Lhtp$-equivalence on affines, since it is a homology equivalence and $\Lhtp BGL$ is an $H$-space (because even permutation matrices are $\A^1$-homotopic to the identity). The map between the upper right corners is also an $\Lhtp$-equivalence by Lemma~\ref{lem:contractible-Emb}. Since the lower right corners are connected, it follows from \cite[Lemma 5.5.6.17]{HA} that $\Lhtp$ preserves these Cartesian squares, and we obtain the desired equivalence on the upper left corners.
\end{proof}

\begin{rem}
	The proof of Proposition~\ref{prop:nfr-to-dfr} shows that the functors $\Lhtp\Emb^\fr(-,\A^\infty)$ and $\Lhtp\Maps_K(0,\sL)$ are given on affines by the formula
\[
(f\colon Z\to X) \mapsto
\begin{cases}
\Maps_{KV(Z)}(0,\sL_f) & \text{if $[\sL_f]=0$ in $K_0(Z)$,} \\
\emptyset & \text{otherwise.}
\end{cases}
\]
where $KV=\Lhtp\tau_{\geq 1} K$ is Karoubi–Villamayor $K$-theory. Note that one cannot replace $KV$ by $\Lhtp K$, even if $X$ is regular, since $Z$ may not be $K_0$-regular.
\end{rem}

\begin{cor}\label{cor:nfr-vs-dfr}
	Let $\sF\in\Pre_\Sigma(\Sch_S)_*$. Then the $\Einfty$-map
	\[
	\Lhtp\h^\nfr_S(\sF) \to \Lhtp\h^\dfr_S(\sF)
	\]
	induced by~\eqref{eqn:hnfr-to-hfr-extended}
	is an equivalence on affines. In particular, it induces an equivalence of presheaves of $\Einfty$-spaces
	\[
	L_\zar\Lhtp \h^\nfr_S(\sF) \simeq L_\zar\Lhtp \h^\dfr_S(\sF).
	\]
\end{cor}

\begin{proof}
	We can assume that $\sF=Y_+$ for some $Y\in\Sch_S$. Then the result follows from Proposition~\ref{prop:nfr-to-dfr} and Lemma~\ref{lem:Sing-on-frames}.
\end{proof}

\begin{rem}
	One can show that Proposition~\ref{prop:extension} remains true if we replace $\Corr^{\efr,n}_S$ and $\Corr^{\nfr,n}_S$ by $\Corr^{\nfr}_S$ and $\Corr^{\dfr}_S$. However, we cannot deduce Corollary~\ref{cor:nfr-vs-dfr} from this because $\h^\dfr_S(Y)$ does not satisfy closed gluing for the boundary of the algebraic $n$-simplex.
\end{rem}

Combining Corollaries \ref{cor:vfr-vs-nfr2} and \ref{cor:nfr-vs-dfr}, we get:

\begin{cor}\label{cor:efr-vs-dfr}
	Let $\sF\in\Pre_\Sigma(\Sm_S)_*$. Then \eqref{einftymap1} and \eqref{eqn:hnfr-to-hfr-extended} induce equivalences of presheaves of $\Einfty$-spaces
	\[
	\Lmot \h^\efr_S(\sF) \simeq \Lmot \h^\nfr_S(\sF) \simeq \Lmot\h^\dfr_S(\sF).
	\]
\end{cor}


\section{The recognition principle} \label{sect:recogmain}

In this section, we prove the recognition principle for motivic spectra.
We start by reviewing the recognition principle for motivic $S^1$-spectra in~\ssecref{ssec:recog/s1recog}.
 We then construct in~\ssecref{ssec:recog/framed-h} the $\infty$-category of (tangentially) framed motivic spaces, denoted by $\H^{\fr}(S)$. This $\infty$-category is built by applying standard constructions of motivic homotopy theory to the $\infty$-category $\Span^{\fr}(\Sm_S)$. While the latter is constructed in~\secref{sect:infinity-category}, we axiomatize the relevant features of this $\infty$-category in~\sssecref{sssec:corrfr-axioms}.
 One key property of $\H^{\fr}(S)$ is that every object has the structure of a presheaf of $\Einfty$-spaces and thus determines a presheaf of connective spectra.
 In~\ssecref{ssec:recog/framed-sh}, we construct the $\infty$-category of (tangentially) framed motivic spectra, $\SH^{\fr}(S)$, and establish its basic properties. The main result of this paper is Theorem~\ref{thm:main}; it states that, when $S$ is the spectrum of a perfect field, $\SH^{\fr}(S)$ is equivalent to the classical stable motivic homotopy $\infty$-category $\SH(S)$ and, furthermore, that grouplike objects in $\H^{\fr}(S)$ embed fully faithfully in $\SH^{\fr}(S)\simeq \SH(S)$. These results rely on the main theorem of~\cite{garkusha2014framed}, which we review in~\ssecref{ssec:recog/gp}.

Throughout this section, we fix a base scheme $S$. 
To simplify the exposition, all schemes will be assumed quasi-compact and quasi-separated (this is first used in Proposition~\ref{prop:H^fr-basics}(ii)).


\ssec{The \texorpdfstring{$S^1$}{S¹}-recognition principle}

\label{ssec:recog/s1recog}

In this section, we recall the recognition principle for motivic $S^1$-spectra over a perfect field, which is essentially due to F.~Morel.

\sssec{}
\label{sssec:Bnis}
Let $\C$ be a presentable $\infty$-category equipped with the Cartesian symmetric monoidal structure.
For every $n\geq 0$, there is an adjunction
\[
\bB^n_\C: \Mon_{\sE_n}(\C) \rightleftarrows \C_* : \Omega^n
\]
where $\bB^n_\C$ is the $n$-fold bar construction \cite[Remark 5.2.3.6 and Example 5.2.2.4]{HA}.
We also have an adjunction
\[
\bB^\infty_\C:\CMon(\C) \rightleftarrows \Stab(\C) : \Omega^\infty,
\]
where $\Stab(\C)=\C\otimes \Spt$ is the stabilization of $\C$.

When $\C$ is $\Pre(\Sm_S)$, $\Pre_\nis(\Sm_S)$, or $\H(S)$, we will write $\bB^n_\C$ as $\bB^n$, $\bB^n_\nis$, or $\bB^n_\mot$, so that
\[
\bB_\nis^n \simeq L_\nis\circ \bB^n\quad\text{and}\quad \bB_\mot^n \simeq \Lmot\circ \bB^n.
\]
We then have a commutative diagram of adjunctions
\begin{equation} \label{eq:tau-commute}
  \begin{tikzcd}
    \CMon(\Pre(\Sch_S)) \ar[shift left=0.5ex]{r}{\bB^\infty} \ar[shift right=0.5ex,swap]{dd}{L_\nis}
      & \Stab(\Pre_{\nis}(\Sch_S)))  \ar[shift right=0.5ex,swap]{dd}{L_\nis} \ar[shift left=0.5ex]{l}{\Omega^\infty}
    \\
     & 
    \\
     \CMon(\Pre_{\nis}(\Sch_S)) \ar[shift left=0.5ex]{r}{\bB^\infty_\nis} \ar[shift right=0.5ex,swap]{dd}{\Lmot} \ar[shift right=0.5ex]{uu}
      & \Stab(\Pre_{\nis}(\Sch_S))) \ar[shift right=0.5ex]{uu} \ar[shift right=0.5ex,swap]{dd}{\Lmot} \ar[shift left=0.5ex]{l}{\Omega^\infty}
    \\
     &
    \\
    \CMon(\H(S)) \ar[shift left=0.5ex]{r}{\bB^\infty_\mot}  \ar[shift right=0.5ex]{uu}
      & \SH^{S^1}(S) \ar[shift right=0.5ex]{uu} \ar[shift left=0.5ex]{l}{\Omega^\infty},
  \end{tikzcd}
 \end{equation}
where $\SH^{S^1}(S)=\Stab(\H(S))$.

\sssec{}\label{sssec:grouplike}
A monoid $X\in\Mon(\C)$ is called \emph{grouplike} if the shearing maps $X\times X\to X\times X$ are equivalences. We denote by $\Mon^\gp(\C)\subset\Mon(\C)$ and $\CMon^\gp(\C)\subset\CMon(\C)$ the full subcategories of grouplike monoids. If finite products distribute over colimits in $\C$, these inclusions have left adjoints called \emph{group completion} and denoted by $X\mapsto X^\gp$.

If $T$ is an $\infty$-topos, then a monoid $X\in\Mon(T)$ is grouplike if and only if its sheaf of connected components $\pi_0(X)$ is grouplike.
In particular, a monoid $X$ in $\Pre_\nis(\Sm_S)$ is grouplike if and only if $\pi_0^\nis(X)$ is a sheaf of groups.

\sssec{} \label{sssec:topos} 
Let $T$ be an $\infty$-topos.
For every $n\geq 0$, the adjunction
\[
\bB^n_T: \Mon_{\sE_n}(T) \rightleftarrows T_* : \Omega^n
\]
of \sssecref{sssec:Bnis}
restricts to an equivalence between grouplike $\sE_n$-algebras and pointed $n$-connective objects \cite[Theorem 5.2.6.15 and Remarks 5.2.6.12]{HA}.
It follows that the adjunction
\[
\bB^\infty_T:\CMon(T) \rightleftarrows \Stab(T) : \Omega^\infty
\]
restricts to an equivalence between grouplike commutative monoids and connective spectra:
\[
\CMon^\gp(T) \simeq \Stab(T)_{\geq 0}.
\]
This holds in particular for $T=\Pre_\nis(\Sm_S)$. However, this does not apply to $T=\H(S)$ (see Remark~\ref{rem:nongp}).

\sssec{}

To formulate the recognition principle for motivic $S^1$-spectra, we will need the following definitions that extend \cite[Definition 1.7]{Morel} to presheaves of spaces.

\begin{defn}\label{def:gp} \leavevmode
	
	\noindent{\em(i)}
	Let $X \in \Mon(\Pre_\nis(\Sm_S))$. We say that $X$ is \emph{strongly $\A^1$-invariant} if it is $\A^1$-invariant and if $\bB_{\nis} X$ is $\A^1$-invariant, i.e., if $\bB_{\nis} X \simeq \bB_{\mot} X$. We denote by $\Mon_\mot(\H(S))\subset\Mon(\H(S))$ the full subcategory of strongly $\A^1$-invariant monoids.
	
	\noindent{\em(ii)}
	Let $X \in \CMon(\Pre_\nis(\Sm_S))$. We say that $X$ is \emph{strictly $\A^1$-invariant} if $\bB_{\nis}^n X$ is $\A^1$-invariant for all $n \geq 0$, i.e., if $\bB^n_{\nis} X \simeq \bB^n_{\mot} X$. We denote by $\CMon_\mot(\H(S))\subset\CMon(\H(S))$ the full subcategory of strictly $\A^1$-invariant commutative monoids.
\end{defn}

\begin{rem} \label{rem:nongp}
	When $S$ is the spectrum of a field, the discrete presheaf $L_\nis\mathbf{Z}[\Gm]$ is an example of a commutative monoid in $\H(S)$ that is not strongly $\A^1$-invariant \cite[Lemma 4.6]{Utsav}. In particular, $L_\nis\mathbf{Z}[\Gm]$ cannot be a loop space in $\H(S)$. This remains true for an arbitrary scheme $S\neq\emptyset$, since there always exists an essentially smooth morphism $S'\to S$ where $S'_\mathrm{red}$ is the spectrum of a field.
	This shows that the results of \sssecref{sssec:topos} do not hold for $T=\H(S)$, and in particular that $\H(S)$ is not an $\infty$-topos.
\end{rem}

\sssec{} Let $\SH^{S^1}(S)_{\geq 0}$ be the full subcategory of $\SH^{S^1}(S)$ spanned by the $S^1$-spectra that are connective in $\Stab(\Pre_\nis(\Sm_S))$:
\[
\SH^{S^1}(S)_{\geq 0} = \SH^{S^1}(S)\cap \Stab(\Pre_\nis(\Sm_S))_{\geq 0} \subset \Stab(\Pre_\nis(\Sm_S)).
\]
 We then have the following tautological result:

\begin{prop} \label{prop:s1} 
	The adjunction
	\[
	\bB^\infty_\mot:\CMon(\H(S)) \rightleftarrows \SH^{S^1}(S): \Omega^\infty
	\]
	restricts to an equivalence $\CMon_\mot^\gp(\H(S)) \simeq \SH^{S^1}(S)_{\geq 0}$.
\end{prop}

\begin{proof} The $\infty$-category $\SH^{S^1}(S)$ is the full subcategory of $\Stab(\Pre_{\nis}(\Sm_S))$ spanned by the $\A^1$-invariant objects. Presenting $\Stab(\Pre_{\nis}(\Sm_S))$ as a limit in $\InftyCat$, we get the following commutative diagram where the vertical arrows are fully faithful inclusions:
\begin{equation} \label{eq:omega-stab}
\begin{tikzcd}
\SH^{S^1}(S) \ar{r} \ar{d} & \cdots   \ar{r}{\Omega}& \H(S)_{\pt} \ar{d} \ar{r}{\Omega} \ar{r} & \H(S)_{\pt}   \ar{d} \\
\Stab(\Pre_{\nis}(\Sm_S)) \ar{r} & \cdots \ar{r}{\Omega} &  \Pre_{\nis}(\Sm_S)_{\pt} \ar{r}{\Omega} \ar{r} &  \Pre_{\nis}(\Sm_S)_{\pt}.
\end{tikzcd}
\end{equation}
Thus, $\SH^{S^1}(S)_{\geq 0}$ is precisely the full subcategory of $\Stab(\Pre_{\nis}(\Sm_S))_{\geq 0}$ consisting of those connective spectra $E$ such that all the deloopings of the $\Einfty$-algebra $\Omega^\infty E$ are $\A^1$-invariant, i.e., such that $\Omega^\infty E\in \CMon_\mot(\H(S))$.
\end{proof}

\sssec{} In order to obtain an explicit description of either side of the equivalence of Proposition~\ref{prop:s1}, we will need to assume $S$ is the spectrum of a (perfect) field $k$. 

The following theorem gives a description of $\CMon_\mot(\H(k))$:

\begin{thm} [Morel]  \label{thm:morel-s1-full}
	Let $k$ be a perfect field.
	
	\noindent{\em(i)}
	Let $X \in \Mon(\H(k))$. Then $X$ is strongly $\A^1$-invariant if and only if $\pi_0^{\nis}(X)$ is strongly $\A^1$-invariant.
	
	\noindent{\em(ii)} 
	Let $X \in \CMon(\H(k))$. Then $X$ is strictly $\A^1$-invariant if and only if $\pi_0^{\nis}(X)$ is strongly $\A^1$-invariant.
\end{thm}

\begin{proof} 
	(i) Since $\pi_0^\nis$ commutes with group completion, we may assume that $X$ is grouplike.
	In this case we have $\Omega \bB_{\nis} X \simeq X$, whence $\pi_0^{\nis}(X) \simeq \pi^{\nis}_1(\bB_{\nis} X)$.
	If $X$ is strongly $\A^1$-invariant, then $\pi^{\nis}_1(\bB_{\nis} X) \simeq \pi_1^{\nis}(\Lmot \bB_{\nis} X)$, which is strongly $\A^1$-invariant by~\cite[Theorem 1.9]{Morel}.

Conversely, suppose that $\pi_0^{\nis}(X)$ is strongly $\A^1$-invariant. We want to show that $\bB_\nis X\simeq \Lmot\bB_\nis X$. We have equivalences of presheaves
 \begin{align*}
X & \simeq \Lmot X \\
 & \simeq  \Lmot \Omega \bB_{\nis} X \\
 & \simeq \Omega \Lmot \bB_{\nis} X.
 \end{align*}
Here, the first equivalence is because $X \in \H(k)$, the second equivalence is because $X$ is grouplike, and the third equivalence follows from~\cite[Theorem 6.46]{Morel} using the assumption that $\pi_0^{\nis}(X)$ is strongly $\A^1$-invariant. Hence $\bB_\nis X\simeq \bB_\nis \Omega \Lmot \bB_{\nis} X$.
 To conclude, we note that $\Lmot$ preserves connected objects by \cite[\sectsign2, Corollary 3.22]{MV}, and hence $\bB_{\nis} \Omega  \Lmot \bB_{\nis} X \simeq \Lmot \bB_{\nis} X$.

	(ii) Necessity follows from (i). Conversely, suppose that $\pi_0^{\nis}(X)$ is strongly $\A^1$-invariant. We want to prove that $\bB_{\nis}^n X \simeq \Lmot \bB_{\nis}^n X$ for all $n \geq 1$, and we proceed by induction on $n$. The case $n=1$ follows from (i), so we may assume $n\geq 2$.
	We have equivalences of presheaves
 \begin{align*}
\bB^{n-1}_{\nis}X & \simeq \Lmot \bB^{n-1}_{\nis}X  \\
 & \simeq  \Lmot \Omega \bB^{n}_{\nis} X \\
 & \simeq \Omega \Lmot  \bB^{n}_{\nis} X.
 \end{align*}
Here, the first equivalence follows from the inductive hypothesis, the second equivalence is because $\bB^{n-1}_{\nis}X$ is connected, and the third equivalence again is an application of \cite[Theorem 6.46]{Morel} using that $\pi_0^\nis(\bB^{n-1}_\nis X)\simeq *$ is strongly $\A^1$-invariant.  Hence, $\bB^n_\nis X\simeq \bB_\nis\Omega \Lmot  \bB^{n}_{\nis} X$.
To conclude, we again use \cite[\sectsign2, Corollary 3.22]{MV} to obtain an equivalence $\bB_{\nis} \Omega  \Lmot \bB^{n}_{\nis} X \simeq \Lmot \bB^{n}_{\nis} X$.
\end{proof}

The following proposition gives a description of $\SH^{S^1}(k)_{\geq 0}$:

\begin{prop}\label{prop:t-structure}
	Let $k$ be a field. Then $\SH^{S^1}(k)_{\geq 0}\subset \SH^{S^1}(k)$ is the full subcategory generated under colimits by $\Sigma^\infty_{S^1}X_+$ for $X\in\Sm_k$.
\end{prop}

\begin{proof}
	Let $C\subset \SH^{S^1}(k)$ be the full subcategory generated under colimits by $\Sigma^\infty_{S^1}X_+$ for $X\in\Sm_k$.
	By the stable $\A^1$-connectivity theorem \cite[Theorem 6.1.8]{morel-connectivity}, we have $C\subset \SH^{S^1}(k)_{\geq 0}$. For the converse, it suffices to show that $\Stab(\Pre_\nis(\Sm_k))_{\geq 0}$ is generated under colimits by suspension spectra of smooth $k$-schemes. Since the Nisnevich topology on $\Sm_k$ is hypercomplete, this follows from \cite[Corollary C.2.1.7]{SAG}.
\end{proof}

\begin{rem}
	Proposition~\ref{prop:t-structure} fails over more general bases $S$. Indeed, Ayoub proved in \cite{Ayoub:2006} that if $S$ is a scheme of finite type over an algebraically closed field with $\dim(S)\geq 2$, then $C\not\subset \SH^{S^1}(S)_{\geq 0}$.
\end{rem}

\begin{cor}\label{cor:s1}
	Let $k$ be a perfect field. 
	Then the adjunction
	\[
	\bB^\infty_\mot:\CMon(\H(k)) \rightleftarrows \SH^{S^1}(k): \Omega^\infty
	\]
	restricts to an equivalence between:
	\begin{itemize}
		\item the full subcategory of commutative monoids $X$ in $\H(k)$ such that $\pi_0^\nis(X)$ is a strongly $\A^1$-invariant sheaf of groups;
		\item the full subcategory of $\SH^{S^1}(k)$ generated under colimits by $\Sigma^\infty_{S^1}X_+$ for $X\in\Sm_k$.
	\end{itemize}
\end{cor}

\begin{proof}
	Combine Proposition~\ref{prop:s1}, Theorem~\ref{thm:morel-s1-full}, and Proposition~\ref{prop:t-structure}.
\end{proof}


\ssec{Framed motivic spaces}
\label{ssec:recog/framed-h}

The $\infty$-category $\Span^{\fr}(\Sm_S)$  can be viewed as the motivic incarnation of Segal's category $\Gamma=\Fin_*^\op$. While we only construct this $\infty$-category in~\sssecref{sssec:framedcorr}, we introduce it informally in \sssecref{sssec:corrfsyn}, and we give an axiomatic description in~\sssecref{sssec:corrfr-axioms}, isolating the properties necessary to prove the recognition principle. 

\sssec{} \label{sssec:corrfsyn} 
In \ssecref{ssec:dfr}, we defined for $X,Y\in\Sch_S$ the $\infty$-groupoid $\Corr^\fr_S(X,Y)$ of framed correspondences from $X$ to $Y$. Such a correspondence is a span
\begin{equation*}
  \begin{tikzcd}
     & Z \ar[swap]{ld}{f}\ar{rd}{g} & \\
    X &   & Y
  \end{tikzcd}
\end{equation*}
over $S$ with $f$ finite syntomic, together with a trivialization of the cotangent complex $\sL_f$ in the $K$-theory of $Z$.

It follows from basic properties of the cotangent complex and of $K$-theory that there is a natural way to compose framed correspondences. The main coherence problem posed by the recognition principle is to construct an $\infty$-category $\Span^\fr(\Sm_S)$ in which the mapping spaces are the $\infty$-groupoids $\Corr^\fr_S(X,Y)$. We will achieve this in \secref{sect:infinity-category}.

If we forget about the trivialization of the cotangent complex, we obtain the groupoid $\Corr^\fsyn_S(X,Y)$ of finite syntomic correspondences from $X$ to $Y$, which already made an appearance in~\sssecref{sssec:CorrF}. These groupoids are the mapping spaces of a $2$-category $\Span^\fsyn(\Sm_S)$. It admits a wide subcategory $\Span^\fet(\Sm_S)$ whose morphisms are finite étale correspondences. Since finite étale maps have \emph{canonical} trivializations of their cotangent complex in $K$-theory, we have a factorization
\[
\Span^\fet(\Sm_S)\to \Span^\fr(\Sm_S) \to \Span^\fsyn(\Sm_S)
\]
of the inclusion $\Span^\fet(\Sm_S)\subset \Span^\fsyn(\Sm_S)$.
The $2$-category $\Span^\fet(\Sm_S)$ has a further wide subcategory $\Span^\clopen(\Sm_S)$ whose morphisms are spans $X \hookleftarrow Z \rightarrow Y$ where $Z \hook X$ is a summand inclusion. This is now a $1$-category and there is an obvious equivalence of categories
\[
\Sm_{S+} \simeq \Span^\clopen(\Sm_S),
\]
whence a functor
\[
\Sm_{S+}\to \Span^\fr(\Sm_S).
\]
This functor is our starting point for an axiomatic description of the $\infty$-category $\Span^\fr(\Sm_S)$.

\sssec{} \label{sssec:corrfr-axioms}

For the remainder of this section, we will assume given a symmetric monoidal $\infty$-category $\Span^\fr(\Sm_S)$ and a symmetric monoidal functor
\[
\gamma\colon \Sm_{S+} \to \Span^\fr(\Sm_S),
\]
inducing an adjunction
\[
\gamma^*: \Pre(\Sm_{S+}) \rightleftarrows \Pre(\Span^\fr(\Sm_S)) : \gamma_*,
\]
and satisfying the following properties:
\begin{enumerate}
	\item The functor $\gamma$ is essentially surjective.
	\item There is a natural equivalence of functors $\h^\fr_S\simeq \gamma_*\gamma^*\colon \Sm_{S+} \to \Pre(\Sm_S)_*$. In particular, the mapping space $\Maps(\gamma(X_+),\gamma(Y_+))$ in $\Span^\fr(\Sm_S)$ is identified with $\Corr^\fr_S(X,Y)$.
	\item Under the equivalence of (2), the natural map $X_+\to \h^\fr_S(X)$ agrees with the unit map $X_+\to \gamma_*\gamma^*(X_+)$.
	\item Under the equivalence of (2), the natural map $X_+\wedge \h^\fr_S(Y) \to \h^\fr_S(X\times_S Y)$ agrees with the composite
	\[
	X_+\wedge \gamma_*\gamma^*(Y_+) \to \gamma_*\gamma^*(X_+)\wedge \gamma_*\gamma^*(Y_+) \to \gamma_*\gamma^*((X\times_S Y)_+),
	\]
	where the first map is the unit map and the second map is the right-lax monoidal structure of $\gamma_*\gamma^*$.
	\item The composition $\pi_0\Corr^\fr_S(X,Y)\times\pi_0\Corr^\fr_S(Y,Z)\to \pi_0\Corr^\fr_S(X,Z)$ in the homotopy category $\h\Span^\fr(\Sm_S)$ sends framed correspondences
	\begin{equation*}
	\begin{tikzcd}
	 & T \ar[swap]{dl}{f} \ar{dr}{g} & \\
	 X  & & Y
	\end{tikzcd}
	\qquad
	\begin{tikzcd}
	 & W \ar[swap]{dl}{h} \ar{dr}{k} & \\
	 Y  & & Z
	\end{tikzcd}
	\end{equation*}
	with trivializations $\alpha\in\pi_0\Maps_{K(T)}(0,\sL_f)$ and $\beta\in\pi_0\Maps_{K(W)}(0,\sL_{h})$ to the framed correspondence
	\[
	\begin{tikzcd}
	 & T\times_YW \ar[swap]{dl}{f\circ \pr_T} \ar{dr}{k\circ \pr_W} & \\
	 X  & & Z
	\end{tikzcd}
	\]
	with trivialization 
	\[
	\phi \circ (\pr_T^*\alpha\oplus(\psi \circ g^*\beta)) \in \pi_0\Maps_{K(T\times_YW)}(0,\sL_{f\circ\pr_T}),
	\]
	 where $\psi\colon g^*\sL_h\simeq\sL_{\pr_T}$ is the canonical isomorphism and $\phi\colon \pr_T^*\sL_f\oplus\sL_{\pr_T}\simeq \sL_{f\circ \pr_T}$ is the isomorphism in $\tau_{\leq 1}K(T\times_YW)$ induced by the canonical cofiber sequence $\pr_T^*\sL_f\to \sL_{f\circ \pr_T}\to \sL_{\pr_T}$.
\end{enumerate}

\begin{rem}
	We note that axiom (5) is not used in the proof of Theorem~\ref{thm:main}(ii). Moreover, to prove Theorem~\ref{thm:main} without the symmetric monoidal structure, only the $\Sm_{S+}$-linear structure of the functor $\gamma$ is needed.
\end{rem}

\sssec{}\label{sssec:semiadditive} We record some immediate consequences of the axioms~\sssecref{sssec:corrfr-axioms}.

\begin{lem}\label{lem:semiadditive}
	The functor $\gamma\colon \Sm_{S+}\to\Span^\fr(\Sm_S)$ preserves finite coproducts and the $\infty$-category $\Span^\fr(\Sm_S)$ is semiadditive.
\end{lem}

\begin{proof}
	By \sssecref{sssec:corrfr-axioms}(1)--(3), $\gamma$ preserves finite coproducts if and only if the functor $\h^\fr_S\colon \Sm_{S+}\to \Pre(\Sm_S)_*$ lands in $\Pre_\Sigma(\Sm_S)_*$. This is the case by~\sssecref{sssec:CorrF-pointing}.
	Semiadditivity similarly follows from Lemma~\ref{lem:f-additivity}.
\end{proof}

It follows from Lemma~\ref{lem:semiadditive} that the $\infty$-category $\Pre_\Sigma(\Span^\fr(\Sm_S))$ is semiadditive \cite[Corollary 2.4]{ggn}.

\begin{lem}\label{lem:bilinear}
	The tensor product in $\Span^\fr(\Sm_S)$ preserves finite coproducts in each variable.
\end{lem}

\begin{proof}
	Since $\gamma\colon \Sm_{S+}\to \Span^\fr(\Sm_S)$ is symmetric monoidal, essentially surjective, and preserves finite coproducts by Lemma~\ref{lem:semiadditive}, this follows from the distributivity of finite products over finite coproducts in $\Sm_S$.
\end{proof}

It follows from Lemma~\ref{lem:bilinear} that $\Pre_\Sigma(\Span^\fr(\Sm_S))$ acquires a unique symmetric monoidal structure whose tensor product preserves colimits in each variable and such that the Yoneda embedding \[\Span^\fr(\Sm_S)\hook \Pre_\Sigma(\Span^\fr(\Sm_S))\] is symmetric monoidal \cite[Proposition 4.8.1.10]{HA}.

\sssec{} \label{sssec:framed-motivic-spaces} We now define the $\infty$-category of framed motivic spaces.
If $\sF\in\Pre_\Sigma(\Span^{\fr}(\Sm_S))$, we say that $\sF$ is \emph{Nisnevich-local} (resp.\ \emph{$\A^1$-invariant}) if its restriction to $\Sm_{S}$ is. We denote by $\Pre_\nis(\Span^{\fr}(\Sm_S))$ and $\Pre_{\Sigma,\A^1}(\Span^{\fr}(\Sm_S))$ the corresponding full subcategories of $\Pre_\Sigma(\Span^{\fr}(\Sm_S))$.

\begin{defn}
A \emph{framed motivic space} over $S$ is an $\A^1$-invariant Nisnevich-local presheaf on $\Span^{\fr}(\Sm_S)$. 
\end{defn}

We denote by $\H^\fr(S) \subset \Pre_\Sigma(\Span^{\fr}(\Sm_S))$ the full subcategory of framed motivic spaces.

\sssec{} \label{sssec:framed-basics} Since $\A^1$-invariance and Nisnevich locality are described by small sets of conditions that are preserved by limits in $\Pre_{\Sigma}(\Span^{\fr}(\Sm_S))$, the inclusions
\[\Pre_\nis(\Span^{\fr}(\Sm_S)),\,\Pre_{\Sigma,\A^1}(\Span^{\fr}(\Sm_S)),\,\H^{\fr}(S) \hookrightarrow \Pre_{\Sigma}(\Span^{\fr}(\Sm_S))\] 
are accessible left localizations. We denote the respective localization functors by $L_\nis$, $\Lhtp$, and $\Lmot$.

\begin{prop} \label{prop:H^fr-basics}\leavevmode

\noindent{\em(i)} The $\infty$-category $\H^\fr(S)$ is generated under sifted colimits by objects of the form $\Lmot\gamma(X_+)$, where $X\in\Sm_S$ is affine.

\noindent{\em(ii)} For every $X\in\Sm_S$, $\Lmot\gamma(X_+)\in\H^\fr(S)$ is compact.

\noindent{\em(iii)} The $\infty$-category $\H^{\fr}(S)$ is semiadditive.

\noindent{\em(iv)} There exists a unique symmetric monoidal structure on $\H^{\fr}(S)$ such that the localization functor $\Lmot\colon \Pre_\Sigma(\Span^\fr(\Sm_S))\to \H^\fr(S)$ is symmetric monoidal.
\end{prop}

\begin{proof} Statement (i) follows by the same argument as in $\H(S)$ (see for example~\cite[Proposition 3.16]{hoyois-sixops} or \cite[Lemma 4.3.5]{khan-BN1}). Statement (ii) follows from the fact the Nisnevich-local $\A^1$-invariant presheaves are closed under filtered colimits in $\Pre_\Sigma(\Span^\fr(\Sm_S))$.
	 Statement (iii) is obvious since $\Pre_\Sigma(\Span^\fr(\Sm_S))$ is semiadditive and $\H^{\fr}(S)$ is closed under finite products. For (iv), it suffices to show that for every $\mathcal F\in\Pre_\Sigma(\Span^\fr(\Sm_S))$, the functor
\[
\mathcal F\otimes(-)\colon \Pre_\Sigma(\Span^\fr(\Sm_S))\to \Pre_\Sigma(\Span^\fr(\Sm_S))
\]
preserves $\Lmot$-equivalences \cite[Proposition 4.1.7.4]{HA}. We may assume that $\mathcal F=\gamma(X_+)$ for some $X\in\Sm_S$, and then the result follows from the fact that $X\times(-)$ preserves Nisnevich sieves and $\A^1$-projections.
\end{proof}

\sssec{}  \label{gamma-star} Since $\gamma\colon \Sm_{S+} \rightarrow \Span^\fr(\Sm_S)$ preserves finite coproducts (Lemma~\ref{lem:semiadditive}), it induces an adjunction
\begin{equation*}
\gamma^*: \Pre_\Sigma(\Sm_S)_{*} \simeq \Pre_\Sigma(\Sm_{S+}) \rightleftarrows   \Pre_{\Sigma}(\Span^\fr(\Sm_S)): \gamma_*,
\end{equation*}
where $\gamma^*$ is symmetric monoidal (for the symmetric monoidal structure on the right-hand side defined in \sssecref{sssec:semiadditive}).

Since the functor $\gamma_*$ preserves $\A^1$-invariant and Nisnevich-local objects by definition, we have a commutative diagram of adjunctions:
\begin{equation*}
    \begin{tikzcd}
    \Pre_{\Sigma}(\Sm_S)_{*} \ar[shift left=.5ex]{r}{\gamma^*}\ar[shift right=.5ex,swap]{dd}{\Lmot}
      &\Pre_{\Sigma}(\Span^{\fr}(\Sm_S))  \ar[shift right=.5ex,swap]{dd}{\Lmot} \ar[shift left=.5ex]{l}{\gamma_*}
    \\
     & \\
    \H(S)_{*} \ar[shift left=.5ex]{r}{\gamma^*}  \ar[shift right=.5ex]{uu}
      & \H^{\fr}(S) \ar[shift right=.5ex]{uu} \ar[shift left=.5ex]{l}{\gamma_*}.
    \end{tikzcd}
  \end{equation*}

\begin{lem} \label{lem:gamma-nis} For any $X\in\Sm_S$ and any finitely generated Nisnevich covering sieve $R \hook X$, $\gamma_*\gamma^*(R_+) \rightarrow \gamma_*\gamma^*(X_+)$ is an $L_{\nis}$-equivalence.
\end{lem}

\begin{proof}
	Note that the image of $R\hook X$ in $\Pre_\Sigma(\Sm_S)$ is a Nisnevich sieve generated by a single map.
	Given the equivalence $\h^\fr_S\simeq \gamma_*\gamma^*$ from \sssecref{sssec:corrfr-axioms}(2), the claim follows from Proposition~\ref{prop:corrfr-descent}(ii).
\end{proof}

\begin{lem} \label{lem:a1}  For any $X \in \Sm_S$, the map $\gamma_*\gamma^*((X \times \A^1)_+) \rightarrow \gamma_*\gamma^*(X_+)$ induced by the projection $X\times\A^1\to X$ is an $L_{\A^1}$-equivalence.
\end{lem}

\begin{proof} 
	Given the equivalence $\h^\fr_S\simeq \gamma_*\gamma^*$ from \sssecref{sssec:corrfr-axioms}(2), the claim follows from Lemma~\ref{lem:retract} for $F = \Maps_K(0,\sL)$.
\end{proof} 

\begin{prop}\label{prop:gamma-local-equiv}
	The functor $\gamma_*\colon \Pre_\Sigma(\Span^\fr(\Sm_S)) \to \Pre_\Sigma(\Sm_S)_*$ commutes with the localization functors $L_\nis$, $\Lhtp$, and $\Lmot$.
\end{prop}

\begin{proof}
	This follows from Lemmas \ref{lem:gamma-nis} and \ref{lem:a1} and \cite[Lemma 2.10]{norms}.
\end{proof}

\begin{prop} \label{prop:gamma-basic} The functor $\gamma_*\colon\H^{\fr}(S) \rightarrow \H(S)_{*}$ is conservative and preserves sifted colimits.
\end{prop}

\begin{proof}  
	Conservativity follows immediately from the fact that $\gamma$ is essentially surjective.
	Note that $\gamma_*\colon \Pre_{\Sigma}(\Span^\fr(\Sm_S)) \rightarrow \Pre_{\Sigma}(\Sm_{S+})$ preserves sifted colimits, since sifted colimits in $\Pre_\Sigma(\C)$ are computed pointwise. Since moreover $\gamma_*$ preserves $\Lmot$-equivalences by Proposition~\ref{prop:gamma-local-equiv}, it follows that $\gamma_*\colon\H^{\fr}(S) \rightarrow \H(S)_{*}$ preserves sifted colimits.
\end{proof}

\sssec{}\label{sssec:lift}
 Since $\H^\fr(S)$ is semiadditive by Proposition~\ref{prop:H^fr-basics}(ii), it follows from \cite[Remark 2.7]{GlasmanGoodwillie} that the functor $\gamma_*$ admits a unique lift as follows:
\begin{equation*}
\begin{tikzcd}
 & \CMon(\H(S)) \ar[d] \\
 \H^{\fr}(S)  \ar{r}{\gamma_*} \ar[dashed]{ur}& \H(S)_{*}.
\end{tikzcd}
\end{equation*}
In other words, every framed motivic space is in particular a presheaf of $\Einfty$-spaces.

\begin{rem}\label{rem:fold}
	 Let $\Span^{\fold}(\Sm_S)$ denote the wide subcatgory of $\Span^{\fet}(\Sm_S)$ where the backward maps are sums of fold maps $X^{\coprod n}\to X$. According to \cite[\sectsign C.1]{norms}, presheaves on this $2$-category that transform finite coproducts into finite products model presheaves of $\Einfty$-spaces on $\Sm_S$. Hence, the above factorization can be concretely described as restrictions along the functors
\[ 
\Sm_{S+}\simeq \Span^{\clopen}(\Sm_S) \rightarrow \Span^{\fold}(\Sm_S) \rightarrow \Span^\fr(\Sm_S).
\] 
\end{rem}


\ssec{Framed motivic spectra} 

\label{ssec:recog/framed-sh}

\sssec{}
Recall that $\T=\A^1/\A^1-0$ in $\Pre_\Sigma(\Sm_S)_*$. We denote by $\G$ the pointed presheaf $(\G_m, 1) \in \Pre_{\Sigma}(\Sm_S)_{*}$, and we let 
\[
\T^{\fr} = \gamma^*(\T)
\quad\text{and}\quad
\G^{\fr}= \gamma^*(\G)
\]
in $\Pre_{\Sigma}(\Span^\fr(\Sm_S))$. 

\sssec{}
We define the presentably symmetric monoidal $\infty$-category of \emph{framed motivic spectra} by formally inverting in $\H^\fr(S)$ the object $\T^\fr$ with respect to the monoidal product:
\[\SH^{\fr}(S) = \H^{\fr}(S)[ (\T^{\fr})^{-1}];\]
see \cite[Definition 2.6]{Robalo}.
We denote the canonical adjunction by
\[ 
\Sigma^{\infty}_{\T,\fr}: \H^{\fr}(S) \rightleftarrows \SH^{\fr}(S): \Omega^{\infty}_{\T,\fr}.
\]

\begin{lem} \label{lem:3-sym} 
	The cyclic permutation of $\T^{\fr}\otimes \T^{\fr}\otimes \T^{\fr}$ in $\H^{\fr}(S)$ is homotopic to the identity.
\end{lem}

\begin{proof} Since the functor $\gamma^*\colon \H(S)_* \rightarrow \H^{\fr}(S)$ is symmetric monoidal, this follows from the standard fact that the cyclic permutation of $\T\wedge\T\wedge\T$ in $\H(S)_*$ is homotopic to the identity.
\end{proof}

Lemma~\ref{lem:3-sym} implies that $\SH^{\fr}(S)$ is computed in $\Mod_{\H^\fr(S)}(\PrL)$ as the colimit of the sequence
\[
\H^{\fr}(S) \xrightarrow{\Sigma_\T} \H^{\fr}(S) \xrightarrow{\Sigma_\T} \H^{\fr}(S) \xrightarrow{\Sigma_\T} \cdots,
\]
where $\Sigma_\T=\T^\fr\otimes(-)$ \cite[Corollary 2.22]{Robalo}.
In particular, the underlying $\infty$-category of $\SH^{\fr}(S)$ is the limit of the tower
\begin{equation}\label{eqn:SHfr}
\H^{\fr}(S) \xleftarrow{\Omega_\T} \H^{\fr}(S) \xleftarrow{\Omega_\T} \H^{\fr}(S) \xleftarrow{\Omega_\T} \cdots ,
\end{equation}
where $\Omega_\T=\Hom(\T^\fr,-)$ is right adjoint to $\Sigma_\T$.
Thus, a framed motivic spectrum is a collection of framed motivic spaces $(\sF_n)_{n \geq 0}$ together with equivalences $\sF_n \simeq \Omega_\T(\sF_{n+1})$ for each $n \geq 0$.

\begin{prop} \label{prop:SH^fr-cpct} 
	\leavevmode
	
	\noindent{\em(i)}
	The $\infty$-category $\SH^\fr(S)$ is stable.
	
	\noindent{\em(ii)}
	The $\infty$-category $\SH^{\fr}(S)$ is generated under sifted colimits by objects of the form $(\T^{\fr})^{\otimes n} \otimes \Sigma^\infty_{\T,\fr}\gamma^*(X_+)$, where $X \in \Sm_S$ is affine and $n \leq 0$.
	
	\noindent{\em(iii)}
	For every $X\in\Sm_S$, $\Sigma^\infty_{\T,\fr}\gamma^*(X_+)$ is compact in $\SH^\fr(S)$.
\end{prop}

\begin{proof}
	Stability of $\SH^\fr(S)$ follows from the equivalence $\T^\fr\simeq S^1\wedge \G^\fr$ in $\H^\fr(S)$.
	Assertion (ii) follows from Proposition~\ref{prop:H^fr-basics}(i) and \cite[Lemma 6.3.3.7]{HTT}.
	By Proposition~\ref{prop:H^fr-basics}(i,ii), $\H^\fr(S)$ is compactly generated and $\Sigma_\T$ preserves compact objects. It follows that the limit of the sequence~\eqref{eqn:SHfr} can be computed in $\mrm{Pr}^\mrm{R}_\omega$ \cite[Proposition 5.5.7.6]{HTT}, so that $\Sigma^\infty_{\T,\fr}$ preserves compact objects.
	Assertion (iii) now follows from Proposition~\ref{prop:H^fr-basics}(ii).
\end{proof}


\ssec{The Garkusha--Panin theorems}

\label{ssec:recog/gp}

To prove our recognition principle we will need, as input, a theorem of Garkusha and Panin \cite{garkusha2014framed}. This, in turn, relies on theorems of Ananyevskiy--Garkusha--Panin \cite{agp}, Garkusha--Neshitov--Panin \cite{gnp}, and Garkusha--Panin \cite{hitr}. We will also need some of these theorems later on. In this subsection, we will briefly review the statements of their theorems, as well as align them with our notation.

\sssec{} As reviewed in~\sssecref{sssec:Bnis}, we have an adjunction 
\[
\bB_{\mot}^{\infty}: \CMon(\H(S)) \rightleftarrows \SH^{S^1}(S):\Omega^\infty.
\] 
According to Remark~\ref{rem:GPadditivity}, for every $\sF\in\Pre_\Sigma(\Sm_S)_*$, $\Lhtp\h^\efr_S(\sF)$ and hence $L_\mot\h^\efr_S(\sF)$ are canonically $\Einfty$-algebras. Moreover, by~\eqref{eq:tau-commute}, we have
\[
L_\mot\bB^\infty \Lhtp\h^\efr(\sF) \simeq \bB^\infty_\mot L_\mot\h^\efr(\sF).
\]
 For $X\in\Sm_S$, the object
\[
\bB^\infty \Lhtp\h^\efr(X) \in\Stab(\Pre_\Sigma(\Sm_S))
\]
is called the ``framed motive'' of $X$ and denoted by $M_\fr(X)$ in \cite[Definition 5.2]{garkusha2014framed}.

\sssec{} \label{sssec:cancellation}
We now recall the main results of \cite{agp} and \cite{garkusha2014framed}.

Let $\sF\in\Pre_\Sigma(\Sm_S)_*$. The right-lax $\Sm_{S+}$-linear structure on the functor $\h^\efr_S$ defined in \sssecref{sssec:efr-module-structure} induces a natural map
\begin{equation}\label{eqn:Gprespectrum}
\G\wedge \bB_\mot^\infty L_\mot \h^\efr_S(\sF) \to \bB_\mot^\infty L_\mot \h^\efr_S(\G\wedge \sF).
\end{equation}
For $X\in\Sm_S$, the sequence of motivic $S^1$-spectra
\[
(\bB_\mot^\infty L_\mot \h^\efr_S(X_+), \bB_\mot^\infty L_\mot \h^\efr_S(\G\wedge X_+), \bB_\mot^\infty L_\mot \h^\efr_S(\G^{\wedge 2}\wedge X_+),\dots)
\]
is then a $\G$-prespectrum in $\SH^{S^1}(S)$, with structure maps given by~\eqref{eqn:Gprespectrum}. 
This is the motivic localization of the $\G$-prespectrum denoted by $M_\fr^\G(X)$ in \cite[Section 11]{garkusha2014framed}.
We will denote by
\[
M^{\efr}_S(X)\in\SH(S)
\]
its associated $\G$-spectrum. The maps
\[
\Sigma^\infty_{S^1} (\G^{\wedge n}\wedge X_+) \to \bB_\mot^\infty L_\mot \h^\efr_S(\G^{\wedge n}\wedge X_+)
\]
in $\SH^{S^1}(S)$ induced by $\sF\to \h^\efr_S(\sF)$ (see \sssecref{sssec:efr-unit}) define a natural map
\begin{equation}\label{eqn:GP-map}
	\Sigma^\infty_\T X_+ \to M^{\efr}_S(X)
\end{equation}
in $\SH(S)$. 

\begin{thm}[Ananyevskiy--Garkusha--Panin]\label{thm:cancel} 
	Suppose that $S$ is the spectrum of an infinite perfect field. 
	For every $X\in\Sm_S$ and $n\geq 0$, the map
	\[
	\bB_\mot^\infty L_\mot \h^\efr_S(\G^{\wedge n}\wedge X_+) \to \Omega_\G\bB_\mot^\infty L_\mot \h^\efr_S(\G^{\wedge n+1}\wedge X_+)
	\]
	adjoint to~\eqref{eqn:Gprespectrum} is an equivalence.
\end{thm}

\begin{proof}
	This is \cite[Theorem B]{agp}, complemented by \cite{DruzhininPanin} in characteristic $2$.
\end{proof}

\begin{thm}[Garkusha–Panin]
	\label{thm:GP}
	Suppose that $S$ is the spectrum of an infinite perfect field. For every $X\in\Sm_S$, the map~\eqref{eqn:GP-map} in $\SH(S)$ is an equivalence.
\end{thm}

\begin{proof}
	This is \cite[Theorem 11.1]{garkusha2014framed}, complemented by \cite{DruzhininPanin} in characteristic $2$.
\end{proof}

\sssec{} \label{sssec:hitr} 
One of the fundamental theorems of Voevodsky's theory of motives is that if $k$ is a perfect field and $\mathcal{F}$ is a $\A^1$-invariant presheaf of abelian groups with transfers on $\Sm_k$, then $L_{\nis}\mathcal{F}$ is strictly $\A^1$-invariant, in the sense that for every $i \geq 0$, the presheaf $H^i_{\nis}(-, L_{\nis}\mathcal{F})\colon \Sm_k^\op \rightarrow \mathrm{Ab}$ is $\A^1$-invariant \cite[Theorem 13.8]{mvw}. It follows that if $E$ is a presheaf of spaces or spectra with transfers on $\Sm_k$, then the canonical map $L_{\nis}\Lhtp E \rightarrow L_{\mot}E$ is an equivalence. 

The main result of \cite{hitr} is an analogous statement for presheaves of abelian groups with framed transfers. To state it, we consider the category $\Span^\efr_*(\Sm_S)$ defined as follows. The objects of $\Span^\efr_*(\Sm_S)$ are smooth $S$-schemes and the set of morphisms from $X$ to $Y$ is
\[
\bigvee_{n\geq 0}\Corr^{\efr,n}_S(X,Y).
\]
Composition is induced by the pairings
\begin{gather*}
\Corr^{\efr,n}_S(X,Y) \times \Corr^{\efr,m}_S(Y,Z) \to \Corr^{\efr,n+m}_S(X,Z), \\
 (f,g)\mapsto (\id_{\T^{\wedge n}}\wedge g) \circ (\id_{(\P^1)^{\wedge m}}\wedge f),
\end{gather*}
using the identification from Remark~\ref{rem:voevodsky-lemma}. 
Note that the wide subcategory of $\Span^\efr_*(\Sm_S)$ on the level $0$ correspondences is equivalent to $\Sm_{S+}$.

A presheaf $\sF\colon \Span^\efr_*(\Sm_S)^\op\to \Set$ is called \emph{radditive} if its restriction to $\Sm_{S+}$ belongs to $\Pre_\Sigma(\Sm_{S+})$.
The equationally framed correspondence
 \[\sigma_X = (0 \times X, \A^1\times X, \pr_{\A^1}, \pr_{X}) \in\Corr^{\efr,1}_S(X,X)\] 
 defines a morphism $\sigma_X\colon X\to X$ in $\Span^\efr_*(\Sm_S)$, and we say that $\sF$ is \emph{stable} if the map 
 \[\sigma_X^*\colon \mathcal{F}(X) \rightarrow \mathcal{F}(X)\] 
 is the identity for every $X \in \Sm_S$.

\begin{thm}[Garkusha--Panin] \label{thm:hitrmain} 
	Let $k$ be an infinite field and $\sF$ an $\A^1$-invariant stable radditive presheaf of abelian groups on $\Span^{\efr}_*(\Sm_k)$. Let $X$ be an essentially smooth semilocal $k$-scheme.
	
	\noindent{\em(i)}
	The Nisnevich sheaf $L_\nis\sF$ is $\A^1$-invariant.
	
	\noindent{\em(ii)}
	The canonical map $\sF(X)\to (L_\nis\sF)(X)$ is an isomorphism.
	
	\noindent{\em(iii)}
	If $k$ is perfect, then $H^i_\nis(-,L_\nis\sF)$ is $\A^1$-invariant for all $i\geq 0$.
	
	\noindent{\em(iv)}
	If $k$ is perfect, then $H^i_\nis(X,L_\nis\sF)=0$ for all $i\geq 1$.
\end{thm}

\begin{proof}
	By \cite[Theorem 2.1]{hitr}, the Nisnevich sheaf $L_\nis\sF$ has a canonical structure of presheaf on $\Span^{\efr}_*(\Sm_k)$, which is also stable and $\A^1$-invariant. In particular, (i) holds.
	To prove (ii), we can assume that $X$ is a semilocalization of an irreducible smooth affine $k$-scheme with generic point $\eta$.
	Then by \cite[Theorem 2.15(3)]{hitr}, the map $\sF(X)\to \sF(\eta)$ is injective.\footnote{In \emph{loc.\ cit.}, it is assumed that $X$ is a localization at a single point. This is only used in a reference to \cite[Lemma 10.1]{OjangurenPanin}, whose proof works in the semilocal case with trivial modifications.} Assertion (ii) now follows by considering the kernel and the cokernel of $\sF\to L_\nis\sF$. If $k$ is perfect, then by \cite[Theorems 16.10 and 16.11]{hitr}, complemented by \cite{DruzhininPanin} in characteristic $2$, $H^i_\nis(-,L_\nis\sF)$ has a structure of stable presheaf on $\Span^{\efr}_*(\Sm_k)$ and is $\A^1$-invariant. In particular, assertion (iii) holds, and assertion (iv) is an immediate consequence of (ii).
\end{proof}

There is a functor
\[
\lambda\colon \Span^{\efr}_*(\Sm_S) \to \h\Span^\fr(\Sm_S)
\]
given on mapping sets by
\[
\bigvee_{n\geq 0} \Corr^{\efr,n}_S(X,Y) \to \bigvee_{n\geq 0}\Corr^{\efr}_S(X,Y) \to \bigvee_{n\geq 0}\pi_0\Corr^{\fr}_S(X,Y)\xrightarrow{\nabla} \pi_0\Corr^{\fr}_S(X,Y).
\]
It is straightforward to check using \sssecref{sssec:corrfr-axioms}(5) that this is indeed a functor and that $\lambda(\sigma_X)=\id_X$ for all $X\in\Sm_S$.
Hence, for any $\sF\in\Pre_\Sigma(\Span^\fr(\Sm_S),\mathrm{Ab})$, the functor $\sF\circ\lambda$ is a stable radditive presheaf of abelian groups on $\Span^{\efr}_*(\Sm_S)$.

\begin{rem}
	One can similarly define a category $\Span^\nfr_*(\Sm_S)$ such that $\lambda$ factorizes as
	\[
	\Span^\efr_*(\Sm_S) \to \Span^\nfr_*(\Sm_S) \to \h\Span^\fr(\Sm_S).
	\]
\end{rem}

\begin{rem}
	Let $C\subset \Sm_S$ be the full subcategory of schemes that are finite sums of affine schemes admitting étale maps to affine bundles over $S$. By Corollaries~\ref{cor:vfr-vs-nfr} and \ref{cor:nfr-vs-dfr}, the assignment $\sF\mapsto \sF\circ\lambda$ defines an equivalence between $\Pre_{\Sigma,\A^1}(\Span^\fr(C),\Set)$ and the category of $\A^1$-invariant stable radditive presheaves of sets on $\Span^\efr_*(C)$. Since the inclusion $C\subset\Sm_S$ induces an equivalence of Nisnevich topoi and framed correspondences are compatible with the Nisnevich topology (in the sense of Propositions \ref{prop:fr-descent}(iii) and \ref{prop:corrfr-descent}(ii)), we deduce an equivalence between $\Pre_{\nis,\A^1}(\Span^\fr(\Sm_S),\Set)$ and the category of $\A^1$-invariant Nisnevich-local stable presheaves of sets on $\Span^\efr_*(\Sm_S)$.
\end{rem}

\begin{thm} \label{thm:hitrcor} 
	Let $k$ be a perfect field and let $\sF\in \Pre_\Sigma(\Span^\fr(\Sm_k))$ be grouplike. Then the canonical map
	\[
	\Lhtp \gamma_*\sF \to \Lmot\gamma_*\sF
	\]
	is an equivalence on essentially smooth semilocal $k$-schemes. In particular,
	\[
	L_\zar\Lhtp \gamma_*\sF \simeq \Lmot\gamma_*\sF.
	\]
\end{thm}

\begin{proof} 
	Suppose first that $k$ is infinite.
We will show that the map $\Lhtp \gamma_*\sF\to L_\nis\Lhtp \gamma_*\sF$ is an equivalence on semilocal schemes and that $L_\nis\Lhtp \gamma_*\sF$ is $\A^1$-invariant.
Using the Postnikov towers of $\Lhtp \gamma_*\sF$ and $L_\nis\Lhtp \gamma_*\sF$, it suffices to prove the following, where $X$ is semilocal:
\begin{itemize}
	\item the canonical map $\pi_n(\Lhtp \gamma_*\sF)(X)\to \pi_n^\nis(\Lhtp \gamma_*\sF)(X)$ is an isomorphism;
	\item for $1\leq i\leq n$, $H^i_\nis(X,\pi_n^\nis(\Lhtp \gamma_*\sF))=0$;
	\item for $0\leq i\leq n$, $H^i_{\nis}(-, \pi_n^\nis(\Lhtp \gamma_*\sF) )$ is $\A^1$-invariant.
\end{itemize}
 All three assertions follow from Theorem~\ref{thm:hitrmain} and Proposition~\ref{prop:gamma-local-equiv}, since $\pi_n(\Lhtp \sF)\circ\lambda$ is a stable radditive presheaf of abelian groups on $\Span^{\efr}_*(\Sm_k)$. 
 
 Suppose now that $k$ is finite, and let $X$ be an essentially smooth semilocal $k$-scheme.
  If $K/k$ is an infinite algebraic extension, $X_{K}$ is still semilocal.
  Since $\Lhtp$ and $\Lmot$ commute with essentially smooth base change, the map $\Lhtp \sF \to \Lmot\sF$ is an equivalence on $X_K$. Let $n\geq 0$ and let $\alpha$ be an element in the kernel or the cokernel of $\pi_n(\Lhtp \sF)(X) \to \pi_n(\Lmot\sF)(X)$. Then we can find finite extensions of $k$ of coprime degrees that kill $\alpha$. It then follows from Proposition~\ref{prop:coprime-injective} that $\alpha=0$.
\end{proof}

\begin{rem}\label{rem:hitrcor}
	By inspection of the proof, we see that the conclusions of Theorem~\ref{thm:hitrcor} hold for any presheaf $\sF\in\Pre_\Sigma(\Sm_k)$ with maps $\mu\colon \sF\times\sF \to \sF$ and $e\colon *\to \sF$ such that:
	\begin{itemize}
		\item $\mu(e,-)\simeq \id_\sF$;
		\item the shearing map $(\pi_1,\mu)\colon\sF\times\sF\to \sF\times\sF$ is an equivalence;
		\item the induced multiplication $\pi_0(\mu)$ on $\pi_0(\sF)$ is associative and commutative;
		\item for every $n\geq 0$, the presheaf of abelian groups $\pi_n(\Lhtp\sF,e)$ can be promoted to a stable presheaf on $\Span^\efr_*(\Sm_k)$.
	\end{itemize}
\end{rem}

\begin{cor}\label{cor:mot-gp}
Let $k$ be a perfect field and let $\sF\in \H^\fr(k)$. Then $\gamma_*\sF\in \CMon(\H(k))$ is strictly $\A^1$-invariant. In other words, there is a factorization
\[\begin{tikzcd}
   \H^\fr(k) \ar{r}{\gamma_*} \ar[dashed]{dr} & \CMon(\H(k)) \\
   & \CMon_\mot(\H(k)). \ar[hookrightarrow]{u}
\end{tikzcd}\]
\end{cor}

\begin{proof}
We must show that $L_\nis\bB^n\gamma_*\sF$ is $\A^1$-invariant for $n\geq 1$. Since $\bB^n\gamma_*\sF\simeq\gamma_*\bB^n\sF$ is connected, it is grouplike and thus Theorem~\ref{thm:hitrcor} applies.
\end{proof}

\begin{cor}
	If $k$ is a perfect field, then \[\H^\fr(k)^\gp=\Pre_{\zar,\A^1}(\Span^\fr(\Sm_k))^\gp\] as subcategories of $\Pre_\Sigma(\Span^\fr(\Sm_k))$.
\end{cor}

\begin{proof}
	This follows directly from Theorem~\ref{thm:hitrcor}.
\end{proof}


\ssec{The recognition principle} \label{ssec:recog/recog}

Having all the ingredients, we now come to the proof of the recognition principle. 
The equivalences
\[
L_\mot\T \simeq L_\mot\Sigma\G\quad\text{and}\quad L_\mot\T^{\fr} \simeq L_\mot \Sigma\G^{\fr}
\]
induce symmetric monoidal factorizations
\[
\Sigma^\infty_\T\simeq \Sigma^\infty_\G \circ \Sigma^\infty_{S^1}\quad\text{and}\quad \Sigma^\infty_{\T,\fr}\simeq \Sigma^\infty_{\G,\fr} \circ \Sigma^\infty_{S^1,\fr}.
\]
We then have a commutative diagram of adjunctions:
\begin{equation} \label{eq:rec-diagram}
  \begin{tikzcd}
    \Sm_{S+} \ar{dd}{\gamma} \ar{r}
      & \H(S)_{*} \ar[shift left=0.5ex]{r}{\Sigma^{\infty}_{S^1}}\ar[shift right=0.5ex,swap]{dd}{\gamma^*}
      & \SH^{S^1}(S)\ar[shift left=0.5ex]{r}{\Sigma^{\infty}_{\G}} \ar[shift right=0.5ex,swap]{dd}{\gamma^*}  \ar[shift left=0.5ex]{l}{\Omega^{\infty}_{S^1}}
      & \SH(S)\ar[shift right=0.5ex,swap]{dd}{\gamma^*}  \ar[shift left=0.5ex]{l}{\Omega^{\infty}_{\G}}
    \\
     & & &
    \\
    \Span^\fr(\Sm_S) \ar{r}
      & \H^{\fr}(S)  \ar[shift right=0.5ex,swap]{uu}{\gamma_*}  \arrow[shift left=0.5ex]{r}{\Sigma_{S^1,\fr}^{\infty}}
      & \SH^{S^1,\fr}(S) \ar[shift left=0.5ex]{r}{\Sigma_{\G,\fr}^{\infty}}  \ar[shift right=0.5ex,swap]{uu}{\gamma_*}   \ar[shift left=0.5ex]{l}{\Omega_{S^1,\fr}^{\infty}}
      & \SH^{\fr}(S).  \ar[shift right=0.5ex,swap]{uu}{\gamma_*}  \ar[shift left=0.5ex]{l}{\Omega_{\G,\fr}^{\infty}}
  \end{tikzcd}
  \end{equation}

\begin{prop} \label{prop:cons+colim}  
	The functors $\gamma_*\colon \SH^{S^1,\fr}(S) \rightarrow \SH^{S^1}(S)$ and $\gamma_*\colon \SH^{\fr}(S) \rightarrow \SH(S)$ are conservative and preserve colimits.
\end{prop}

\begin{proof} 
	It follows immediately from Proposition~\ref{prop:gamma-basic} and \cite[Proposition 5.5.7.6]{HTT} that these functors are conservative and preserve filtered colimits.
	They also preserve finite colimits since they are right adjoint functors between stable $\infty$-categories.
\end{proof}

\sssec{}
As a first step towards the recognition principle, we transport the results of \ssecref{ssec:recog/s1recog} to the framed setting.

\begin{lem}\label{lem:s1-commutes}
Let $\C$ and $\D$ be presentable $\infty$-categories in which finite products commute with sifted colimits, and let $f\colon \C\to \D$ be a functor that preserves finite limits and sifted colimits.
	 Then the following square commutes:
	\begin{equation*}
	\begin{tikzcd} 
	\CMon(\C) \ar[swap]{d}{\CMon(f)} \ar{r}{\bB^{\infty}_{\C}} & \Stab(\C) \ar{d}{\Stab(f)} \\
	\CMon(\D) \ar{r}{\bB^{\infty}_{\D}} & \Stab(\D).
	\end{tikzcd}
	\end{equation*}
\end{lem}

\begin{proof}
The assumption on $\C$ and $\D$ implies that $\CMon(\C)$ and $\CMon(\D)$ have sifted colimits and are accessible \cite[Corollaries 3.2.3.2 and 3.2.3.5]{HA}.
They also admit finite coproducts by \cite[Corollary 3.2.4.8]{HA}, hence are presentable.
Moreover, $\CMon(f)$ preserves all colimits since it preserves sifted colimits and finite coproducts.
The forgetful functors $\CMon(\C)\to\C$ and $\CMon(\D)\to\D$ induce equivalences on $\CMon$ \cite[Example 3.2.4.5]{HA} and on $\Stab$ (since stable $\infty$-categories are semiadditive). Replacing $f$ by $\CMon(f)$, we can therefore assume that $f$ preserves finite limits and all colimits, so that it has a right adjoint. Then the square of right adjoints obviously commutes.
\end{proof}

\begin{prop}\label{prop:s^1-rec-fr}
The functor $\Sigma^{\infty}_{S^1,\fr}\colon \H^{\fr}(S) \rightarrow \SH^{S^1,\fr}(S)$ is fully faithful when restricted to the full subcategory of grouplike strictly $\A^1$-invariant objects.
\end{prop}

\begin{proof} 
	We need to prove that for every $X \in \H^{\fr}(S)^{\gp}$ that is strictly $\A^1$-invariant, the unit map $X \rightarrow \Omega^{\infty}_{S^1, \fr}\Sigma^{\infty}_{S^1, \fr}X$ is an equivalence. Since $\gamma_*$ is conservative by Proposition~\ref{prop:gamma-basic}, it suffices to prove this after applying $\gamma_*$. By Lemma~\ref{lem:s1-commutes}, we then need to prove that the unit map $\gamma_*X \rightarrow \Omega^{\infty}_{S^1}\bB_{\mot}^{\infty}\gamma_*X$ is an equivalence, which follows from Proposition~\ref{prop:s1}.
\end{proof}

\begin{cor}\label{cor:s^1-rec-fr}
	If $k$ is a perfect field, the functor $\Sigma^{\infty}_{S^1,\fr}\colon \H^{\fr}(k)^\gp \rightarrow \SH^{S^1,\fr}(k)$ is fully faithful.
\end{cor}

\begin{proof}
	Combine Proposition~\ref{prop:s^1-rec-fr} and Corollary~\ref{cor:mot-gp}.
\end{proof}

\sssec{}
We now prove a cancellation theorem for framed motivic $S^1$-spectra. It is a direct analogue of Voevodsky's cancellation theorem in $\DM^{\eff}(k)$ \cite{voevodsky-cancel}: if $M\in \DM^{\eff}(k)$, then the unit map $M \rightarrow \Hom(\Z(1),M(1))$ is an equivalence.

\begin{thm}[Cancellation Theorem] \label{thm:cancellation}
	Let $k$ be a perfect field.
	 For every $M \in \SH^{S^1,\fr}(k)$, the unit map $\M \rightarrow \Omega_{\G}(\G^\fr\otimes M)$ is an equivalence.
\end{thm}

\begin{proof} 
	By Corollary~\ref{cor:infinite-conservative}(3), we can assume $k$ infinite.
	By Proposition \ref{prop:cons+colim}, it suffices to check that for any $X \in \Sm_k$, the canonical map 
\[\gamma_*\gamma^* \Sigma^{\infty}_{S^1}X_+ \rightarrow  \Omega_{\G}\gamma_*\gamma^*\Sigma^{\infty}_{S^1}(\G\wedge X_+) 
\]
is an equivalence.  
By~\sssecref{sssec:corrfr-axioms}(2,4), this map can be indentified with the map
\[
\bB^{\infty}_{\mot}\Lmot \h^{\fr}(X)  \to \Omega_{\G}\bB^{\infty}_{\mot}\Lmot \h^{\fr}(\G\wedge X_+)
\]
induced by the map $\G\wedge \h^\fr(X)\to \h^\fr(\G\wedge X_+)$ from \sssecref{sssec:dfr-module-structure}.
By \sssecref{sssec:nfr-module-structure} and \sssecref{sssec:dfr-module-structure}, we have a commutative square
\begin{equation*} \label{eqn:cancel-diagram1}
\begin{tikzcd}
\bB^{\infty}_{\mot}\Lmot \h^{\efr}(X) \ar{d} \ar{r} & \Omega_{\G}\bB_{\mot}^{\infty}\Lmot \h^{\efr}(\G\wedge X_+)   \ar{d}\\
\bB^{\infty}_{\mot}\Lmot \h^{\fr}(X)  \ar{r} & \Omega_{\G}\bB^{\infty}_{\mot}\Lmot \h^{\fr}(\G\wedge X_+).
\end{tikzcd}
\end{equation*}
The vertical maps are equivalences by Corollary \ref{cor:efr-vs-dfr}, and the top horizontal map is an equivalence by Theorem~\ref{thm:cancel}. Hence, the bottom horizontal map is an equivalence, as desired.
\end{proof}

\begin{cor} \label{cor:tan-framed-cancel}
	If $k$ is a perfect field, the functor $\Sigma^{\infty}_{\G,\fr}\colon \SH^{S^1,\fr}(k) \rightarrow \SH^{\fr}(k)$ is fully faithful.
\end{cor}

\begin{proof}
	Given $M\in\SH^{S^1,\fr}(k)$, the unit map $M\to\Omega^\infty_{\G,\fr}\Sigma^\infty_{\G,\fr}M$ is the colimit of the sequence
	\[
	M\to \Omega_\G(\G^\fr\otimes M) \to \Omega^2_\G((\G^{\fr})^{\otimes 2}\otimes M) \to \cdots.
	\]
	Every map in this sequence is an equivalence by Theorem~\ref{thm:cancellation}.
\end{proof}

\begin{cor}\label{cor:T-cancel}
	If $k$ is a perfect field, the functor $\Sigma^{\infty}_{\T,\fr}\colon \H^{\fr}(k)^\gp \rightarrow \SH^{\fr}(k)$ is fully faithful.
\end{cor}

\begin{proof}
	Combine Corollaries~\ref{cor:s^1-rec-fr} and \ref{cor:tan-framed-cancel}.
\end{proof}

\sssec{} Finally, we prove that the rightmost vertical adjunction in~\eqref{eq:rec-diagram} is an equivalence.

\begin{thm}[Reconstruction Theorem]
	\label{thm:identify} If $k$ is a perfect field, then the adjunction 
\begin{equation*}
\gamma^*: \SH(k) \rightleftarrows \SH^{\fr}(k): \gamma_*
\end{equation*}
is an equivalence of symmetric monoidal $\infty$-categories.
\end{thm}

\begin{proof} 
	Since the functor $\gamma_*$ is conservative by Proposition~\ref{prop:cons+colim}, it suffices to prove that $\gamma^*$ is fully faithful, i.e., that the unit transformation $\id \rightarrow \gamma_*\gamma^*$ is an equivalence. Since $\SH(k)$ is generated under colimits by the $\T$-desuspensions of smooth $k$-schemes and $\gamma_*$ preserves colimits by Proposition~\ref{prop:cons+colim}, we need only prove that the unit map
\[\Sigma^{-n}_\T\Sigma^{\infty}_\T X_+ \rightarrow \gamma_*\gamma^*\Sigma^{-n}_\T\Sigma^{\infty}_\T X_+\] 
is an equivalence for all $X \in \Sm_k$ and $n \geq 0$. Since the right-hand side is equivalent to $\Sigma^{-n}_\T\gamma_*\gamma^*\Sigma^{\infty}_\T X_+$, we can assume $n=0$. Moreover, by Corollary~\ref{cor:infinite-conservative}(5), we can assume $k$ infinite.

 Since $\gamma_*\gamma^*$ preserves motivic equivalences (Lemmas \ref{lem:gamma-nis} and \ref{lem:a1}), the motivic spectrum $\gamma_*\gamma^*\Sigma^{\infty}_\T X_+$ is the $\G$-spectrum associated with the $\G$-prespectrum
\[
(\bB^\infty_\mot L_\mot\gamma_*\gamma^*(X_+),\bB^\infty_\mot L_\mot\gamma_*\gamma^*(\G\wedge X_+),\dots)
\]
in $\SH^{S^1}(k)$.
By \sssecref{sssec:corrfr-axioms}(2)--(4), this $\G$-prespectrum can be identified with
\[
(\bB^\infty_\mot L_\mot\h^\fr(X_+),\bB^\infty_\mot L_\mot\h^\fr(\G\wedge X_+),\dots),
\]
in such a way that the unit map is induced by the maps $\G^{\wedge n}\wedge X_+\to\h^\fr(\G^{\wedge n}\wedge X_+)$.
By \sssecref{sssec:nfr-module-structure} and \sssecref{sssec:dfr-module-structure}, these maps factorize as
\[
\G^{\wedge n}\wedge X_+\to\h^\efr(\G^{\wedge n}\wedge X_+)\to\h^\fr(\G^{\wedge n}\wedge X_+),
\]
compatibly with the structure maps of the corresponding $\G$-prespectra. Since the second map is a motivic equivalence by Corollary \ref{cor:efr-vs-dfr}, we are reduced to the statement of Theorem~\ref{thm:GP}.
\end{proof}

\sssec{}\label{sssec:recognition}
Recall that there are full subcategories
\[
\SH^\veff(S)\subset \SH^\eff(S)\subset \SH(S)
\]
defined as follows:
\begin{itemize}
	\item $\SH^\eff(S)\subset\SH(S)$ is the full subcategory generated under colimits by $\Sigma^{-n}\Sigma^\infty_\T X_+$ for $X\in\Sm_S$ and $n\geq 0$ (see \cite[\sectsign 2]{Voevodsky:2002});
	\item $\SH^\veff(S)\subset\SH(S)$ is the full subcategory generated under colimits and extensions by $\Sigma^\infty_\T X_+$ for $X\in\Sm_S$ (see \cite[Definition 5.5]{SpitzweckOstvaer}).
\end{itemize}
The $\infty$-category $\SH^\eff(S)$ is stable and $\SH^\veff(S)$ is the nonnegative part of a $t$-structure on $\SH^\eff(S)$. 

\begin{thm}[Motivic Recognition Principle]
\label{thm:main}
Let $k$ be a perfect field. 
	
\noindent{\em(i)}
The functor 
\[ \gamma_*\Sigma^{\infty}_{\T,\fr}\colon \H^{\fr}(k)^{\gp} \rightarrow \SH(k)
\]
is fully faithful and induces an equivalence of symmetric monoidal $\infty$-categories
		\[
		\H^\fr(k)^\gp \simeq \SH^\veff(k).
		\]

\noindent{\em(ii)}
The functor 
\[ \gamma_*\Sigma^{\infty}_{\G,\fr}\colon \SH^{S^1,\fr}(k) \rightarrow \SH(k)
\]
is fully faithful and induces an equivalence of symmetric monoidal $\infty$-categories
	\[
	\SH^{S^1,\fr}(k) \simeq \SH^\eff(k).
	\]
\end{thm}

\begin{proof}
	That these functors are fully faithful follows from Corollary~\ref{cor:tan-framed-cancel}, Corollary~\ref{cor:T-cancel}, and Theorem~\ref{thm:identify}. It remains to identify the images.
	In (i), the image is the full subcategory of $\SH(k)$ generated under colimits by $\Sigma^\infty_\T X_+$ for $X\in\Sm_k$, which is all of $\SH^\veff(k)$ by \cite[Remark after Proposition 4]{BachmannSlices}.
	In (ii), the image is the full subcategory of $\SH(k)$ generated under colimits by $\Sigma^{-n}\Sigma^\infty_\T X_+$ for $X\in\Sm_k$ and $n\geq 0$, which by definition is $\SH^\eff(k)$.
\end{proof}

\begin{cor}\label{cor:sifted}
	Let $k$ be a perfect field. Then the functor
	\[
	\Omega^\infty_\T\colon \SH^\veff(k) \to \H(k)
	\]
	is conservative and preserves sifted colimits. In particular, it is monadic.
\end{cor}

\begin{proof}
	This follows immediately from Theorem~\ref{thm:main}(i) and Proposition~\ref{prop:gamma-basic}.
\end{proof}

Recall that we denote by $X^\gp$ the group completion of a commutative monoid $X$. By \cite[Lemma 5.5]{HoyoisCdh}, the localization functors $L_\zar$, $L_\nis$, $\Lhtp$, and $L_\mot$ (regarded as functors to the subcategory of local objects) commute with group completion.

\begin{cor}\label{cor:main}
	Suppose that $S$ is pro-smooth over a perfect field $k$, and let $\sF\in \Pre_\Sigma(\Sm_S)_*$ be in the image of the pullback functor $\Pre_\Sigma(\Sm_k)_*\to \Pre_\Sigma(\Sm_S)_*$.
	Then there are canonical equivalences
\begin{align*}
\Omega_\T^\infty\Sigma_\T^\infty \Lmot \sF &\simeq L_\zar\Lhtp \h^{\fr}_S(\sF)^\gp\\
&\simeq L_\zar(\Lhtp \h^{\nfr}_S(\sF))^\gp\\
&\simeq L_\zar(\Lhtp \h^{\efr}_S(\sF))^\gp.
\end{align*}
\end{cor}

\begin{proof}
Since all nine functors involved commute with pro-smooth base change, we can assume that $S=\Spec k$. 
In this case, we have
	\[
	\Omega^\infty_\T\Sigma^\infty_\T\Lmot \sF \simeq L_\mot\h^\fr_S(\sF)^\gp
	\]
	by Theorem~\ref{thm:main}(i). The first equivalence then follows from Theorem~\ref{thm:hitrcor} and the second from Corollary~\ref{cor:nfr-vs-dfr}. 
	For the last equivalence, since $\Lhtp$ preserves colimits, it is enough to show that the map $(\Lhtp \h^{\efr}_S(Y))^\gp\to (\Lhtp \h^{\nfr}_S(Y))^\gp$ is a Zariski equivalence for every $Y\in\Sm_S$. By Corollary~\ref{cor:vfr-vs-nfr2}, it is a motivic equivalence. Hence, we must show that
	\[
	L_\zar(\Lhtp \h^{\efr}_S(Y))^\gp \simeq (L_\mot\h^{\efr}_S(Y))^\gp.
	\]
	This follows from Remark~\ref{rem:hitrcor}, since $(\Lhtp \h^{\efr}_S(Y))^\gp$ is a presheaf on $\Span^\efr_*(\Sm_S)$ whose homotopy groups are stable by \cite[Lemma~3.1.4]{MuraThesis}. 
\end{proof}

\sssec{} It is worthwhile to specialize Corollary~\ref{cor:main} to $\sF=S$. Consider the presheaf
\[
\FSYN^\fr_S\colon \Sch^\op_S \to \Spc, \quad X\mapsto \FSyn^\fr_X,
\]
where $\FSyn^\fr_X=\Corr^\fr_S(X,S)$ is the $\infty$-groupoid of \emph{framed finite syntomic $X$-schemes}, i.e., finite syntomic $X$-schemes with trivialized cotangent complex in $K$-theory; this is an $\Einfty$-semiring under disjoint union and Cartesian product of $X$-schemes.

Noting that every field is pro-smooth of its prime subfield, which is perfect, the first equivalence of Corollary~\ref{cor:main} becomes:

\begin{thm}\label{thm:BPQ}
	Suppose that $S$ is pro-smooth over a field. Then there is a canonical equivalence of $\Einfty$-ring spaces
	\[
	\Omega^\infty_\T\mathbf 1_S \simeq L_\zar\Lhtp (\FSYN^\fr_S)^\gp.
	\]
\end{thm}

We can regard Theorem~\ref{thm:BPQ} as a motivic analog of the Barratt–Priddy–Quillen theorem, which describes the topological sphere spectrum as $\Omega^\infty \mathbf 1 \simeq \Fin^\gp$, where $\Fin$ is the groupoid of finite sets.

We shall revisit the second equivalence of Corollary~\ref{cor:main} in \ssecref{ssec:rep-results}.


\section{The \texorpdfstring{$\infty$}{∞}-category of framed correspondences}  \label{sect:infinity-category}

In this technical section, we construct the symmetric monoidal $\infty$-category $\Span^{\fr}(\Sch_S)$ of framed correspondences.
This is an $\infty$-category whose objects are $S$-schemes and whose mapping spaces are the $\infty$-groupoids $\Corr^{\fr}_S(X,Y)$ defined in \ssecref{ssec:dfr}. 
This construction will be an example of an $\infty$-category of \emph{labeled correspondences}, which is the subject of~\ssecref{ssec:labeled-corr} (see Definition~\ref{defn:f-labelled}). 
A labeled correspondence is a span whose left and/or right leg carries additional data, called a ``label''. The structure needed to coherently compose such labeled correspondence is a \emph{labeling functor}. In the case of interest, a finite syntomic morphism is to be labeled by a trivialization of its cotangent complex in $K$-theory, and we construct the desired labeling functor in~\ssecref{ssec:label-fr}. Finally, in~\ssecref{ssec:symmon}, we construct the symmetric monoidal structure on $\Span^\fr(\Sch_S)$.


\ssec{\texorpdfstring{$\infty$}{∞}-Categories of labeled correspondences}
\label{ssec:labeled-corr}

\sssec{} \label{sssec:intuition} 
We start with a heuristic discussion of the idea of labeled correspondences.
Let $\C$ be an $\infty$-category with pullbacks. In \cite[Section 3]{BarwickMackey}, Barwick introduces an $\infty$-category $\Span(\C)$ of \emph{spans} in $\C$. The objects of this $\infty$-category are the objects of $\C$ and $1$-morphisms between two objects $X$ and $Y$ are spans
\begin{equation*}
\begin{tikzcd}
 & Z \ar[swap]{dl}{f} \ar{dr}{g} & \\
 X  & & Y.
\end{tikzcd}
\end{equation*}
One composes spans using pullbacks and the higher morphisms partake in witnessing higher coherences for composition.
In addition, if $\M$ is a collection of morphisms in $\C$ that is stable under composition and pullbacks, we may demand that the left leg belongs to $\M$, which defines a wide subcategory $\Span(\C,\M)\subset\Span(\C)$. There is moreover a functor
\[
\C\to \Span(\C,\M)
\]
that sends a morphism $g\colon X\to Y$ to the span $X\xleftarrow{\id} X\xrightarrow{g}Y$.

Assume now that to every morphism $f$ in $\M$ is associated a space $\F(f)$ of ``labels''. A span as above together with a point in $\F(f)$ will be called an \emph{$\F$-labeled correspondence} from $X$ to $Y$.
Our goal is to construct an $\infty$-category $\Span^\F(\C,\M)$ whose morphisms are $\F$-labeled correspondences, together with a factorization of the above functor as
\[
\C\to\Span^\F(\C,\M) \to \Span(\C,\M).
\]
In our application, $\C$ will be the $1$-category of schemes, $\M$ the class of finite syntomic morphisms, and $\F(f)$ the $\infty$-groupoid of trivializations of the cotangent complex $\sL_f$ in $K$-theory.

Of course, such a construction will require some additional structure on the assignment $f\mapsto \F(f)$.
First of all, let $\Phi_1(\C,M)\subset\Fun(\Delta^1,\C)$ denote the subcategory where the objects are morphisms in $\M$ and the $1$-morphisms are pullback squares. Then $f\mapsto \F(f)$ should be a presheaf on $\Phi_1(\C,M)$. This guarantees that the space $\Corr^\F(X,Y)$ of $\F$-labeled correspondences from $X$ to $Y$ is a functor of $X$ and $Y$ (cf.\ \sssecref{sssec:CorrF}). However, this is clearly not enough structure to define \emph{composition} of $\F$-labeled correpondences.

\sssec{} \label{sssec:3-coherence} Suppose that we are trying to compose three $\F$-labeled correspondences. We can first perform the composition in $\Span(\C,\M)$:
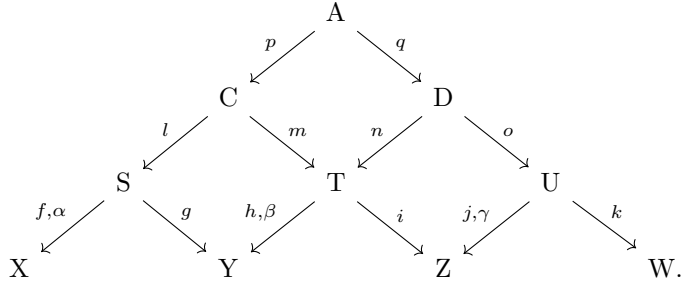
\begin{equation*}
\begin{tikzcd}
 & & &  A \ar[swap]{dl}{p} \ar{dr}{q}& & &\\
 & &C \ar[swap]{dl}{l} \ar{dr}{m}& &D\ar[swap]{dl}{n} \ar{dr}{o} & & \\
&S \ar[swap]{dl}{f,\alpha} \ar{dr}{g} & &\ar[swap]{dl}{h, \beta}T \ar{dr}{i}& &\ar[swap]{dl}{j, \gamma}U \ar{dr}{k} &\\
X& &Y&  &Z& &W.
\end{tikzcd}
\end{equation*}
Here, $\alpha\in \F(f)$, $\beta\in\F(h)$, and $\gamma\in\F(j)$ are given labels, and we want to construct a label in $\F(f\circ l \circ p)$. The functoriality of $\F$ for pullbacks determines labels $g^*(\beta)\in\F(l)$, $i^*(\gamma)\in\F(n)$, and $m^*i^*(\gamma)\in\F(p)$. To go further we need composition laws
\[
\circ\colon \F(f)\times\F(f') \to \F(f\circ f')
\]
for every pair of composable morphisms in $\M$. Given these, we have three \emph{a priori} different ways to obtain a label in $\F(f\circ l\circ p)$:
\[
\alpha\circ g^*(\beta\circ i^*(\gamma)),\quad \alpha\circ (g^*(\beta)\circ m^*i^*(\gamma)),\quad (\alpha\circ g^*(\beta))\circ m^*i^*(\gamma),
\]
and we want them to be equivalent. Thus, we should ask for given $2$-cells in the squares
\[
	\begin{tikzcd}
	\F(h) \times \F(n) \ar{d} \ar{r} & F(h\circ n) \ar{d} \\
	\F(l) \times \F(p) \ar{r} & F(l \circ p),
	\end{tikzcd}
	\qquad
\begin{tikzcd}
\F(f) \times \F(l) \times \F(p) \ar{d} \ar{r} & F(f\circ l) \times \F(p) \ar{d} \\
\F(f) \times \F(l \circ p) \ar{r} & F(f \circ l \circ p).
\end{tikzcd}
\]
In other words, composition of labels should be simultaneously associative and compatible with pullbacks. If $\C$ is a $1$-category and $\F$ is set-valued, then the data exhibited so far suffices to define the $(2,1)$-category $\Span^\F(C,M)$. But if $\C$ or $\F$ contains higher morphisms, we will need that associativity and compatibility with pullbacks hold up to coherent homotopy.

\sssec{} We begin our construction by recalling some basics from the theory of spans; references include \cite[Sections 2--5]{BarwickMackey} and \cite[Section 5]{haugseng}. For $p\geq 0$, let $\Sigma_p$ denote the 1-category whose objects are pairs $(i,j)$ with $0 \leq i \leq j \leq p$ and unique morphisms $(i, j) \rightarrow (i', j')$ whenever $i \leq i'$ and $j' \leq j$. This category is also called the \emph{twisted arrow category} of $\Delta^p$ \cite[\sectsign 5.2.1]{HA}. For example, $\Sigma_0=*$, $\Sigma_1$ is the ``walking span'' diagram
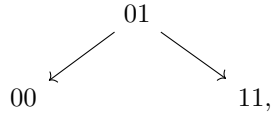
\begin{equation*}
\begin{tikzcd}
 & \ar{dl} 01  \ar{dr} & \\
00  & & 11,
\end{tikzcd}
\end{equation*}
and $\Sigma_2$ is the pyramid
\begin{equation*}
\begin{tikzcd}
 & &02 \ar{dl} \ar{dr}& &   \\
&01 \ar{dl} \ar{dr} & &\ar{dl} 12 \ar{dr}& \\
00& &11&  &22.
\end{tikzcd}
\end{equation*}

\sssec{} The categories $\Sigma_p$ assemble into a functor $\Sigma_{\bullet}\colon \Delta \rightarrow \Cat$; given a monotone map $\alpha\colon [p] \rightarrow [q]$, we get a map 
$$\alpha_*\colon \Sigma_p \rightarrow \Sigma_q,\quad \alpha_*(i,j) = (\alpha(i), \alpha(j)).$$
Inside $\Sigma_p$, there is a full subcategory $\Lambda_p \subset \Sigma_p$ spanned by objects $(i, j)$ where $j -i \leq 1$. This is displayed as the ``bottom of the pyramid":
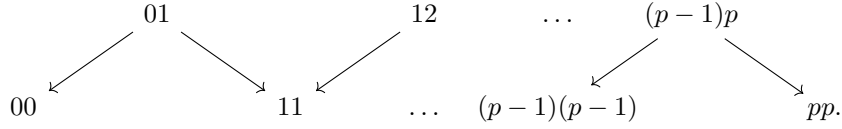
\begin{equation*}
\begin{tikzcd}[column sep={5em,between origins}]
 & 01 \ar{dl} \ar{dr} &  & 12 \ar{dl}  &  \dots &   (p-1)p \ar{dl} \ar{dr}  &  \\
00  &  & 11  &  \dots &  (p-1)(p-1) &  & pp .
\end{tikzcd}
\end{equation*}
 
\sssec{} Let $\C$ be an $\infty$-category. A functor $F\colon \Sigma_p \rightarrow \C$ is \emph{Cartesian} if it is the right Kan extension of its restriction to $\Lambda_p\subset\Sigma_p$. Equivalently, $F$ is Cartesian if for every $i\leq i'$ and $j\leq j'$, it sends the square
\[
\begin{tikzcd}
	ij' \ar{r}\ar{d} & i'j' \ar{d} \\ ij \ar{r} & i'j
\end{tikzcd}
\]
to a Cartesian square.
We denote by $\Span_p(\C) \subset \Fun(\Sigma_p, \C)^{\simeq}$ the full subcategory spanned by the Cartesian functors. An object of $\Span_p(\C)$ is called a \emph{$p$-span}.

\sssec{} Suppose given two collections of morphisms $M$ and $N$ in $\C$ that are stable under composition (in particular, contain identities) and under pullbacks along one another. Then we denote by $\Span_p(\C,M,N)\subset \Span_p(\C)$ the full subcategory consisting of the Cartesian functors $\Sigma_p\to \C$ that send left morphisms to $M$ and right morphisms to $N$.

The construction $(\C,M,N)\mapsto \Span_p(\C,M,N)$ is a functor on the $\infty$-category $\Trip$ of \emph{triples} defined as follows \cite[Section 5]{BarwickMackey}:
\begin{itemize}
	\item An object of $\Trip$ is a triple $(\C,M,N)$ where $\C$ is an $\infty$-category and $M$ and $N$ are collections of morphisms in $\C$ that are stable under composition and under pullbacks along one another (assumed to exist).
	\item A morphism $f\colon (\C,M,N)\to (\C',M',N')$ in $\Trip$ is a functor $f\colon \C\to \C'$ that preserves pullbacks of morphisms in $M$ along morphisms in $N$ and such that $f(M) \subset M'$ and $f(N) \subset N'$.
\end{itemize}

If $N$ is the collection of all morphisms in $\C$, the triple $(\C,M,N)$ will be abbreviated to $(\C,M)$. 

\sssec{}
The functor $\Sigma_\bullet\colon \Delta\to\Cat$ lifts to a functor
\[
\Delta\to\Trip,\quad [p]\mapsto (\Sigma_p,\Sigma_p^L,\Sigma_p^R),
\]
where $\Sigma_p^L$ is the collection of left morphisms and $\Sigma_p^R$ that of right morphisms. We can recast the definition of $\Span_p(\C,M,N)$ as
\[
\Span_p(\C,M,N) = \Maps_{\Trip}((\Sigma_p,\Sigma_p^L,\Sigma_p^R),(\C,M,N)),
\]
and in particular we obtain a simplicial $\infty$-groupoid $\Span_\bullet(\C,M,N)$.

The following is \cite[Proposition 5.9]{BarwickMackey}.

\begin{prop}[Barwick] \label{prop:span-complete} 
	For every triple $(\C,M,N)$, the simplicial space $\Span_{\bullet}(\C, M,N)$ is a complete Segal space.
\end{prop}

\sssec{}
Let $C$ be an $\infty$-category. We denote by $\coCart_C$ the subcategory of $\InftyCat_{/C}$ whose objects are the coCartesian fibrations and whose morphisms are the functors that preserve coCartesian edges. Morphisms in $\coCart_C$ will also be called \emph{strict morphisms of coCartesian fibrations} to emphasize the latter condition.
The $\infty$-categorical Grothendieck construction \cite[\sectsign3.2]{HTT} provides an equivalence of $\infty$-categories
\[
\Fun(C,\InftyCat) \simeq \coCart_C,\quad F\mapsto \int_C F.
\]

Moreover, by \cite[Corollary A.31]{GHN}, the Grothendieck construction is natural in $C$: there exists a functor
\[
\InftyCat^\op \to \Fun(\Delta^1,\InftyCat),\quad C\mapsto \int_C\colon \Fun(C,\InftyCat) \to \coCart_C.
\]

\begin{lem}\label{lem:enhanced-Grothendieck}
	Let $u\colon E\to \InftyCat$ be a Cartesian fibration, and let $E_{\sslash\InftyCat}\to E$ be the Cartesian fibration classified by $\Fun(u(-),\InftyCat)\colon E^\op\to\InftyCat$.
	Then there is a functor
	\[
	\int\colon E_{\sslash\InftyCat}\to E
	\]
	such that for every $X\in E$ and $F\colon u(X)\to\InftyCat$, $u(\int F)\simeq \int_{u(X)} F$.
\end{lem}

\begin{proof}
	The Grothendieck construction gives a functor
	\[
	\Fun(u(X'),\InftyCat) \stackrel\int\simeq \coCart_{u(X')} \subset \InftyCat{}_{/u(X')},
	\]
	natural in $X'\in E^\op$. 
	Applying the Cartesian version of the Grothendieck construction to this natural transformation, we obtain a morphism of Cartesian fibrations
	\[
	E_{\sslash\InftyCat} \to \Fun^\mathrm{cart}(\Delta^1,E)
	\]
	over $E$. The desired functor is the composition of this morphism with evaluation at $0\in\Delta^1$.
\end{proof}

\sssec{} \label{sssec:Phi}
 Given a triple $(\C, M,N)$ and $n\geq 0$, we denote by $\Phi_n(\C, M,N) \subset \Fun(\Delta^n,\C)$ the subcategory whose objects are the functors sending every edge of $\Delta^n$ to $M$ and whose morphisms are the Cartesian transformations with components in $N$. 

For every map $\alpha\colon [m]\to [n]$ in $\Delta$, the restriction functor $\alpha^*\colon \Fun(\Delta^n,\C)\to\Fun(\Delta^m,\C)$ sends $\Phi_n(\C,M,N)$ to $\Phi_m(\C,M,N)$.
Similarly, for every morphism of triples $f\colon (\C,M,N)\to (\C',M',N')$, the functor $f_*\colon \Fun(\Delta^n,\C)\to \Fun(\Delta^n,\C')$ sends $\Phi_n(\C,M,N)$ to $\Phi_n(\C',M',N')$. It follows that we have a functor
\[
\Phi_\bullet\colon \Trip \to\Fun(\Delta^\op,\InftyCat), \quad (\C,M,N)\mapsto ([n]\mapsto \Phi_n(\C,M,N)).
\]

\sssec{} \label{sssec:CorrF2}
Let $(\C,M,N)$ be a triple and let $F\colon \Phi_1(\C,M,N)^\op\to\Spc$ be a presheaf. For $X,Y\in\C$, denote by $\Corr(X,Y)$ the $\infty$-groupoid of spans $X\leftarrow Z\to Y$ with left leg in $M$ and right leg in $N$. 
The $\infty$-groupoid of \emph{$F$-labeled correspondences} from $X$ to $Y$ is then defined by the Cartesian square
\[
\begin{tikzcd}
\Corr^F(X,Y) \ar{r}\ar{d} & \Spc_* \ar{d}\\
\Corr(X,Y) \ar{r} & \Spc,	
\end{tikzcd}
\]
where the bottom horizontal map sends $X\xleftarrow{f} Z\to Y$ to $F(f)$. 
In other words, a point in $\Corr^{F}(X,Y)$ is a labeled span
\begin{equation*}
\begin{tikzcd}
 & Z \ar[swap]{dl}{f,\alpha} \ar{dr}{g} & \\
 X  & & Y
\end{tikzcd}
\end{equation*}
with $f\in M$, $g\in N$, and $\alpha\in F(f)$.
Equivalently,
\[
\Corr^F(X,Y) =\colim_{X\xleftarrow fZ\to Y}F(f).
\]
This defines a functor
\[
\Corr^F(-,-)\colon \C^\op\times\C \to \Spc.
\]

\sssec{} \label{sssec:unwind1} We now introduce the notion of labeling functor on a triple $(\C,M,N)$.

\begin{defn}\label{def:segal-presheaf}
	Let $X_\bullet\colon\Delta^\op\to\InftyCat$ be a simplicial $\infty$-category. A \emph{Segal presheaf} on $X_\bullet$ is a functor
	\[
	F\colon \int_{\Delta^\op} X_\bullet^\op \to\Spc
	\]
	satisfying the following condition: for every $n\geq 0$ and $\sigma \in X_n$, the map 
	\[
	F(\sigma) \rightarrow F(\rho_1^*(\sigma)) \times \cdots \times F(\rho_n^*(\sigma))
	\]
	 induced by the Segal maps $\rho_i\colon [1] \rightarrow [n]$ is an equivalence.
\end{defn}

We denote by $\Seg\Pre(X_\bullet)\subset \Fun(\int X_\bullet^\op,\Spc)$ the full subcategory of Segal presheaves.

\begin{defn} \label{def:labels} 
Let $(\C, M,N) \in \Trip$. A \emph{labeling functor} on $(\C,M,N)$ is a Segal presheaf on $\Phi_\bullet(\C,M,N)$. 
\end{defn}

A \emph{triple with labeling functor} $(\C,M,N;F)$ is a triple $(\C,M,N)$ with a labeling functor $F$. A morphism $(\C,M,N;F) \to (\C',M',N';F')$ is a morphism of triples $f\colon (\C,M,N)\to (\C',M',N')$ together with a morphism of Segal presheaves $F\to F'\circ f$.
We denote by $\Lab\Trip$ the $\infty$-category of triples with labeling functors. More precisely, $\Lab\Trip\to\Trip$ is the Cartesian fibration classified by
\[
\Trip^\op\to \InftyCat, \quad (\C,M,N)\mapsto \Seg\Pre(\Phi_\bullet(\C,M,N)).
\]
If $N$ is the collection of all morphisms of $\C$, we abbreviate $(\C,M,N;F)$ to $(\C,M;F)$.

\sssec{} Let us unpack some of the structure in a labeling functor $F$ on $(\C,M,N)$. The main part of $F$ is the presheaf $f\mapsto F(f)$ on $\Phi_1(\C,M,N)$, which defines the functor $(X,Y)\mapsto \Corr^F(X,Y)$ (see \sssecref{sssec:CorrF2}).

The Segal condition forces $F$ to be contractible on $\Phi_0(\C,M,N)$, which is the wide subcategory of $\C$ with morphisms in $N$.
The simplicial degeneracy map
\[
s_0\colon \Phi_0(\C,M,N) \to \Phi_1(\C,M,N)
\]
sends $X\in \C$ to $\id_X$, and hence we have a map $*\simeq F(X) \to F(\id_X)$, making $F(\id_X)$ a pointed space. In particular, we obtain an element $\id_X\in \Corr^F(X,X)$.

An object of $\Phi_2(\C,M,N)$ is a pair of composable morphisms $\bullet \stackrel{g}{\rightarrow} \bullet \stackrel{f}{\rightarrow} \bullet$ in $M$, and the simplicial face maps
\[
d_0,d_1,d_2\colon \Phi_2(\C,M,N)\to \Phi_1(\C,M,N)
\]
send $(f,g)$ to $f$, $f\circ g$, and $g$, respectively. We therefore have a span
\[
F(f)\times F(g) \xleftarrow{(d_0,d_2)} F(f,g) \xrightarrow{d_1} F(f\circ g),
\]
and the Segal condition forces the left map to be an equivalence. In particular, we obtain
\[
F(f)\times F(g) \to F(f\circ g).
\]
Using this, one can construct a composition law
\[
\Corr^F(X,Y) \times \Corr^F(Y,Z) \to \Corr^F(X,Z),
\]
and one can check that it is associative up to homotopies induced by the simplicial maps $d_1,d_2\colon \Phi_3(\C,M,N)\to \Phi_2(\C,M,N)$. The rest of the data provided by the Segal presheaf $F$ will guarantee that the $\infty$-groupoids $\Corr^F(X,Y)$ are in fact the mapping spaces of an $\infty$-category.

\sssec{} 
Let $(\C,M,N)\in\Trip$. By definition, an $n$-span $\sigma\in \Span_n(\C,M,N)$ is a morphism of triples
\[
\sigma\colon (\Sigma_n,\Sigma_n^L,\Sigma_n^R) \to (\C,M,N).
\]
If $F$ is a labeling functor on $(\C,M,N)$, we can consider the composition
\begin{align*}
\Span_n(\C,M,N) &= \Maps_{\Trip}((\Sigma_n,\Sigma_n^L,\Sigma_n^R), (\C,M,N)) \\
&\textstyle\xrightarrow{\int\Phi_\bullet^\op} \Fun_{\Delta^\op}(\int\Phi_\bullet(\Sigma_n,\Sigma_n^L,\Sigma_n^R)^\op, \int\Phi_\bullet(\C,M,N)^\op) \\
&\textstyle\xrightarrow{F} \Fun(\int\Phi_\bullet(\Sigma_n,\Sigma_n^L,\Sigma_n^R)^\op,\Spc)\xrightarrow{\lim} \Spc.
\end{align*}
By the Grothendieck construction, this classifies a map of $\infty$-groupoids $\Span^F_n(\C,M,N)\to \Span_n(\C,M,N)$. Explicitly, $\Span^F_n(\C,M,N)$ is defined by the Cartesian square
\[
\begin{tikzcd}
	\Span^F_n(\C,M,N) \ar{r} \ar{d} & \Fun(\int\Phi_\bullet(\Sigma_n,\Sigma_n^L,\Sigma_n^R)^\op,\Spc_\pt) \ar{d} \\
	\Span_n(\C,M,N) \ar{r} & \Fun(\int\Phi_\bullet(\Sigma_n,\Sigma_n^L,\Sigma_n^R)^\op,\Spc).
\end{tikzcd}
\]
Note that this square is functorial in $[n]\in \Delta^\op$, so that we have defined a simplicial space $\Span^F_\bullet(\C,M,N)$ over $\Span_\bullet(\C,M,N)$.

\sssec{}\label{sssec:Span-functoriality}
The construction $(\C,M,N;F)\mapsto \Span_\bullet^{F}(\C,M,N)$ can be promoted to a functor
\[
\Lab\Trip \to \Fun(\Delta^\op,\Spc).
\]
To construct this functor, we use \cite[Corollary 7.6]{GHN} to express the Grothendieck construction as a colimit over the twisted arrow $\infty$-category: if $E\to \C$ is the Cartesian fibration classified by $F\colon \C\to\InftyCat^\op$, then
\begin{equation}\label{eqn:lax-colimit}
E\simeq \colim_{(c\to c')\in\Tw(\C)} F(c')\times \C_{/c}. 
\end{equation}
In particular, the identity functor on $E$ corresponds to a natural family of functors
\begin{equation}\label{eqn:lax-colimit2}
F(c') \times \C_{/c} \to E.
\end{equation}

Recall that $\Lab\Trip\to\Trip$ is the Cartesian fibration classified by
\[
\Trip\to\InftyCat^\op,\quad (\C,M,N)\mapsto \Seg\Pre(\Phi_\bullet(\C,M,N)).
\]
Let us abbreviate a triple $(\C,M,N)$ to $\C$ and $\Seg\Pre(\Phi_\bullet(\C,M,N))$ to $\Seg(\C)$, as the construction will not depend on the definition of $\C\mapsto\Seg(\C)$. By~\eqref{eqn:lax-colimit}, we get
\[
\Lab\Trip \simeq \colim_{(\C\to\C')\in\Tw(\Trip)}\Seg(\C') \times \Trip_{/\C}.
\]
The fiberwise mapping space (see \sssecref{sssec:vop}) provides a functor
\[
\Maps(*,-)\colon \int_{\Delta^\op}\Seg(\Sigma_\bullet) \to \Spc, \quad ([n],F\in\Seg(\Sigma_n)) \mapsto \lim F.
\]
Since $\Span_\bullet(\C)=\Maps_{\Trip}(\Sigma_\bullet,\C)$, we obtain a functor
\[
\textstyle
\Seg(\C') \to  \Fun_{\Delta^\op}(\int\Span_\bullet(C'), \int\Seg(\Sigma_\bullet)) \xrightarrow{\Maps(*,-)} \Fun(\int\Span_\bullet(C'),\Spc),
\]
natural in $\C'\in\Trip^\op$. 
We also have
\[
\textstyle
\int\Span_\bullet(-)\colon \Trip_{/\C} \to (\coCart_{\Delta^\op})_{/\int\Span_\bullet(\C)},
\]
natural in $\C\in\Trip$.
Hence we get
\[
\textstyle
\Seg(\C')\times\Trip_{/\C} \to \Fun(\int\Span_\bullet(\C'),\Spc)\times (\coCart_{\Delta^\op})_{/\int\Span_\bullet(\C)},
\]
natural in $(\C,\C')\in\Trip\times\Trip^\op$. Using~\eqref{eqn:lax-colimit2}, we have a functor
\[
\textstyle
\Fun(\int\Span_\bullet(\C'),\InftyCat)\times (\coCart_{\Delta^\op})_{/\int\Span_\bullet(\C)} \to (\coCart_{\Delta^\op})_{\sslash\InftyCat},
\]
natural in $(\C\to \C')\in\Tw(\Trip)$.
Finally, we use the functor
\[
\int\colon (\coCart_{\Delta^\op})_{\sslash\InftyCat}\to \coCart_{\Delta^\op}
\]
from Lemma~\ref{lem:enhanced-Grothendieck} to obtain
\[
\Seg(\C')\times\Trip_{/\C} \to \coCart_{\Delta^\op}\simeq \Fun(\Delta^\op,\InftyCat), \quad (F,D\to C) \mapsto \Span_\bullet^{F}(D),
\]
natural in $(\C\to \C')\in\Tw(\Trip)$.

\sssec{}\label{sssec:Span-functoriality2}
As a special case of \sssecref{sssec:Span-functoriality}, if $(\C,M,N)\in\Trip$, we have a functor
\[
\Seg\Pre(\Phi_\bullet(\C,M,N)) \to \Fun(\Delta^\op,\Spc),\quad F\mapsto \Span^{F}_\bullet(\C,M,N).
\]

\sssec{}
The rest of this section will be dedicated to the proof of the following theorem.

\begin{thm}\label{thm:segal}
	Let $(\C,M,N;F)$ be a triple with labeling functor. Then the simplicial space $\Span_\bullet^{F}(\C,M,N)$ is a Segal space.
\end{thm}

\begin{rem}
	A labeling functor $F\in\Seg\Pre(\Phi_\bullet(\C,M,N))$ is designed to label the \emph{left} legs of the spans in $\Span_\bullet(\C,M,N)$. 
	By symmetry, if $G\in\Seg\Pre(\Phi_\bullet(\C,N,M))$, we also have a Segal space $\Span^G_\bullet(\C,M,N)$ whose morphisms are spans whose \emph{right} leg is labeled by $G$. We can moreover combine the two constructions to obtain a Segal space
	\[
	\Span^{F,G}_\bullet(\C,M,N) =\Span^F_\bullet(\C,M,N)\times_{\Span_\bullet(\C,M,N)}\Span^G_\bullet(\C,M,N).
	\]
\end{rem}

\sssec{} 
We begin with some combinatorial preliminaries.

\begin{lem} \label{lem:decompose}  
	Let $I$ be a poset and $K$ a category. Let $\mathrm{Sub}(K)$ be the poset of full subcategories of $K$, and let $I \rightarrow \mathrm{Sub}(K)$, $i\mapsto K_i$, be a morphism of posets. Suppose that:
	\begin{enumerate}
		\item each $K_i$ is a cosieve, i.e., if $x\in K_i$ and $x\to y$ is a morphism in $K$, then $y\in K_i$;
		\item $K = \bigcup_{i \in I} K_i$;
		\item for every $i,j\in I$, there exists $k\in I$ such that $k\leq i$, $k\leq j$, and $K_k=K_i\cap K_j$.
	\end{enumerate}	
	 Then, for any $\infty$-category $\C$ and any functor $F\colon K \rightarrow \C$, the canonical map
	 \[
	 \lim_K F\to \lim_{i\in I^\op}\lim_{K_i} F
	 \]
	 is an equivalence.
\end{lem} 

\begin{proof} 
	This is an application of \cite[Corollary 4.2.3.10]{HTT}. 
	According to \cite[Remark 4.2.3.9]{HTT}, it suffices to check that for every simplex $\sigma$ of $N(K)$, the poset
	\[
	I_\sigma=\{i\in I\:|\: \sigma\in N(K_i)\}
	\]
	is weakly contractible. 
	 By assumption (1), $I_\sigma=I_{\sigma(0)}$, so we can assume that $\sigma$ is an object of $K$.
	Assumption (2) then implies that $I_\sigma$ is nonempty, and assumption (3) implies that, for every $i,j\in I_\sigma$, there exists $k\in I_\sigma$ with $k\leq i$, $k\leq j$. Hence, $I_\sigma$ is cofiltered and in particular weakly contractible.
\end{proof}

\begin{lem}\label{lem:initial-object}
	Let $m,n\geq 0$ and let $M$ (resp.\ $N$) be the collection of maps of the form $(f,\id)$ (resp.\ of the form $(\id,f)$) in $\Delta^m\times\Delta^n$. Then the category $\int\Phi_\bullet(\Delta^m\times\Delta^n,M,N)^\op$ has an initial object given by the $m$-simplex $\Delta^m\simeq \Delta^m\times\{n\}\hookrightarrow \Delta^m\times\Delta^n$.
\end{lem}

\begin{proof}
	Let $\sigma\colon \Delta^k\to \Delta^m\times\Delta^n$ be an object in $\int\Phi_\bullet(\Delta^m\times\Delta^n,M,N)^\op$. Since every edge of $\Delta^k$ is sent to $M$, there is a unique $0\leq i\leq n$ and a unique $k$-simplex $\tau\colon \Delta^k \to\Delta^m$ such that $\sigma$ has the form
	\[
	\Delta^k \xrightarrow{\tau} \Delta^m \simeq \Delta^m\times\{i\} \hookrightarrow \Delta^m\times\Delta^n.
	\]
	By the coCartesian factorization system \cite[Example 5.2.8.15]{HTT}, a map from $\id\times n$ to $\sigma$ consists of $\tau'\colon \Delta^k\to\Delta^m$ and a map $\tau\circ(\id\times i)\to \tau'\circ(\id\times n)$ in $\Phi_k(\Delta^m\times\Delta^n,M,N)$, i.e., a Cartesian transformation with components in $N$. Since $\Delta^n$ is a poset, there is at most one such transformation, and there is none unless $\tau'=\tau$.
\end{proof}

\sssec{}
Let $F$ be a labeling functor on $(\C,M,N)$ and let $\sigma\in \Span_n(\C,M,N)$ be an $n$-span. 
Then there is an induced labeling functor 
\[F_\sigma\colon \int\Phi_\bullet(\Sigma_n,\Sigma_n^L,\Sigma_n^R)^\op\to\Spc\]
 on $(\Sigma_n,\Sigma_n^L,\Sigma_n^R)$, and by definition of $\Span_n^F(\C,M,N)$ we have a Cartesian square
\[
\begin{tikzcd}
	\lim F_\sigma \ar{r} \ar{d} & \Span_n^F(\C,M,N) \ar{d} \\
	* \ar{r}{\sigma} & \Span_n(\C,M,N).
\end{tikzcd}
\]
In other words, $\lim F_\sigma$ is the space of $F$-labelings of $\sigma$. The following key proposition expresses that an $F$-labeling of an $n$-span is uniquely determined by the labels of the $n$ left maps at the bottom of the pyramid.

\begin{prop}\label{prop:segal}
	Let $n\geq 0$ and let $F$ be a labeling functor on $(\Sigma_n,\Sigma_n^L,\Sigma_n^R)$. 
	For $0\leq i\leq n-1$, let $f_i$ be the morphism $(i,i+1) \to (i,i)$ in $\Sigma_n$.
	Then the canonical map
	\[
	\lim F \to F(f_0)\times \cdots \times F(f_{n-1})
	\]
	is an equivalence.
\end{prop}

\begin{proof}
	For $0\leq i\leq j\leq n$, let $\Sigma_n^{ij}$ be the full subcategory $\Sigma_n$ spanned by $(p,q) \in \Sigma_n$ such that $p \leq i$ and $q \geq j$.
	In other words, $\Sigma_n^{ij}\simeq \Delta^{n-j}\times\Delta^i$ is the sieve in $\Sigma_n$ generated by the vertex $(i,j)$. There is a functor 
	\[\Sigma_n \rightarrow \Cat_{/\Sigma_n},\quad (i,j) \mapsto \Sigma_n^{ij},\]
	sending each morphism to an inclusion of posets. 
	We let 
	\[
	\Sigma_{n}^{ij,L}=\Sigma_n^L\cap \Sigma_n^{ij}\quad\text{and}\quad\Sigma_{n}^{ij,R}=\Sigma_n^R\cap \Sigma_n^{ij}.
	\]
	We now apply Lemma~\ref{lem:decompose} with $K=\int\Phi_\bullet(\Sigma_n,\Sigma_n^L,\Sigma_n^R)^\op$, $I=\Sigma_n$, and
	\[
	I\to\mathrm{Sub}(K),\quad (i,j)\mapsto K_{ij}=\int\Phi_\bullet(\Sigma_n^{ij},\Sigma_n^{ij,L},\Sigma_n^{ij,R})^\op.
	\]
	Assumption (1) follows from the fact that $\Sigma_n^{ij}\subset \Sigma_n$ is a sieve, assumption (2) follows from the fact that $\Sigma_n = \bigcup_{ij} \Sigma_n^{ij}$, and assumption (3) follows from the fact that $(i,j) \mapsto \Sigma_n^{ij}$ preserves binary products. We deduce that
	\[
	\lim F \simeq \lim_{(i,j)\in\Sigma_n^\op} \lim_{K_{ij}} F.
	\]
	By Lemma~\ref{lem:initial-object}, $K_{ij}$ has an initial object $\sigma_{ij}\colon \Delta^{n-j}\to \Sigma_n^{ij}$ given by the sequence of left morphisms
	\[
	(i,n) \to (i,n-1) \to \cdots \to (i,j),
	\]
	hence $\lim_{K_{ij}}F\simeq F(\sigma_{ij})$.
Moreover, since the inclusion $\Lambda_n\subset \Sigma_n$ is cofinal, we have
	\begin{equation}\label{eqn:limLambda}
	\lim F \simeq \lim_{(i,j)\in\Lambda_n^\op} F(\sigma_{ij}).
	\end{equation}
	Finally, since $F$ is a Segal presheaf, the canonical map
	\[
	F(\sigma_{ij}) \to F((i,n) \to (i,n-1)) \times\cdots \times F((i,j+1)\to (i,j))
	\]
	is an equivalence.
	
	We now prove the proposition by induction on $n$. 
	If $n=0$, then $F$ is a functor $\Delta^\op \to\Spc$ such that $F(0)$ is contractible (by the Segal condition), hence $\lim F\simeq F(0)$ is contractible. If $n\geq 1$, we have 
	\begin{align*}
	\lim F & \simeq F(\sigma_{00}) \times_{F(\sigma_{01})} \lim_{(i,j)\in\Lambda_n^\op,i\geq 1} F(\sigma_{ij}) & \text{(by \eqref{eqn:limLambda})} \\
	&\simeq F(\sigma_{00}) \times_{F(\sigma_{01})} F(f_1)\times \cdots\times F(f_{n-1}) & \text{(by the induction hypothesis)} \\
	&\simeq (F(f_0)\times F(\sigma_{01})) \times_{F(\sigma_{01})} F(f_1)\times \cdots \times F(f_{n-1}) & \text{(by the Segal condition)} \\
	&\simeq F(f_0)\times \cdots \times F(f_{n-1}), &
	\end{align*}
	as desired.
\end{proof}

\begin{proof}[Proof of Theorem~\ref{thm:segal}]
	 For $n\geq 0$, consider the square
	\begin{equation}\label{eq:Segal}
	\begin{tikzcd}
	\Span^F_n(\C,M,N) \ar{r} \ar{d} & \Span^F_1(\C,M,N) \times_{\Span^F_0(\C,M,N)}\cdots\times_{\Span^F_0(\C,M,N)}\Span^F_1(\C,M,N) \ar{d} \\
	\Span_n(\C,M,N) \ar{r} & \Span_1(\C,M,N) \times_{\Span_0(\C,M,N)}\cdots\times_{\Span_0(\C,M,N)}\Span_1(\C,M,N).
	\end{tikzcd}
	\end{equation}
	We must show that the upper horizontal map is an equivalence. Since $\Span_\bullet(\C,M,N)$ is a Segal space, the lower horizontal map is an equivalence, so it suffices to show that this square is Cartesian. Fix $\sigma\in \Span_n(\C,M,N)$, and let $F_\sigma$ be the induced labeling functor on $(\Sigma_n,\Sigma_n^L,\Sigma_n^R)$. 
	Since $\Span^F_0(\C,M,N)\simeq \Span_0(\C,M,N)$, the vertical fiber of~\eqref{eq:Segal} over $\sigma$ is then the restriction map
	\[
	\lim F_\sigma \to \lim F_{\rho_1^*(\sigma)} \times \cdots\times \lim F_{\rho_n^*(\sigma)}.
	\]
	It is an equivalence by Proposition~\ref{prop:segal}.
\end{proof}

\sssec{} \label{sssec:complete-segal-space}
For a Segal space $E_{\bullet}$, we call $E_0$ its \emph{space of objects} and $E_1$ its \emph{space of morphisms}; for each $x, y \in E_0$ the space of maps between them is 
\[\Maps_{E_{\bullet}}(x, y)=\{x\} \times_{E_0} E_1 \times_{E_0} \{y\}.\]
We denote by $E_{\bullet}^+$ the complete Segal space associated with $E_\bullet$ \cite[Section 14]{rezk-css}. The natural map of simplicial spaces $\eta\colon E_{\bullet} \to E_{\bullet}^+$ is a Dwyer–Kan equivalence in the sense of \cite[Section 7.4]{rezk-css}, so that for each pair of objects $x, y \in E_0$, the map between mapping spaces
 $$\Maps_{E_{\bullet} }(x, y) \rightarrow \Maps_{E_{\bullet}^+}(\eta x, \eta y)$$ 
is an equivalence \cite[Section 14, Statement 3]{rezk-css}.
 
 \begin{defn} \label{defn:f-labelled}
 	Let $(\C,M,N;F)$ be a triple with labeling functor. The \emph{$\infty$-category of $F$-labeled correspondences} is the complete Segal space $\Span_\bullet^{F}(\C,M,N)^+$. 
 \end{defn}
 
 Unpacking the definition, we find an equivalence
 \[
 \Maps_{\Span_\bullet^{F}(\C,M,N)}(X,Y) \simeq \Corr^{F}(X,Y)
 \]
 natural in $(X,Y)\in\C^\op\times \C$ (see \sssecref{sssec:CorrF2}).
 
 \begin{rem}\label{rmk:not-complete}
 We note that $\Span^{F}_{\bullet}(\C,M)$ is usually \emph{not} a complete Segal space. Indeed, for $X \in \Span^{F}_0(\C,M) \simeq \C^\simeq$, any invertible element $\alpha$ of the $\sA_\infty$-space $F(\id_X)$ defines an invertible $F$-labeled correspondence from $X$ to $X$, but if $\alpha\not\simeq 1$ then this correspondence is not in the image of the degeneracy map $s_0 \colon \C^\simeq \rightarrow \Span^{F}_1(\C,M)$.
 \end{rem}


\ssec{The labeling functor for tangential framings}
\label{ssec:label-fr}

The formalism introduced in \ssecref{ssec:labeled-corr} allows us to construct $\infty$-categories of correspondences $\Span^F(\C,M)$ where we label the left leg of a span by some additional data. 
For the purposes of the recognition principle in motivic homotopy theory, we are interested in the example where $(\C,M)=(\Sch,\fsyn)$ and the additional data is a tangential framing.

\sssec{} \label{sssec:hCorr^fr}
Let ``$\syn$'' denote the collection of syntomic morphisms of schemes. 
Our labeling functor will be a Segal presheaf
 \[
 \fr\colon \int_{\Delta^\op}\Phi_\bullet(\Sch, \syn)^{\op} \rightarrow \Spc
 \]
that extends the functor
\[
\Phi_1(\Sch,\syn)^\op \to \Spc, \quad (f\colon Y\to X) \mapsto \Maps_{K(Y)}(0,\sL_f).
\]

\begin{rem}
	Note that $\Maps_{K(Y)}(0,\sL_f)=\emptyset$ if $f$ has positive relative dimension, so the labeling functor $\fr$ is only relevant for quasi-finite syntomic morphisms. 
	A more interesting labeling functor on $(\Sch,\syn)$, which we will not attempt to construct, would send $f\colon Y\to X$ to $\Maps_{K(Y)}(0,\sL_f-\sO_Y^{\rk\sL_f})$.
\end{rem}

To give an idea of what is involved in the construction, consider a pair of composable syntomic morphisms $Z\xrightarrow gY\xrightarrow{f}X$, viewed as an object of $\Phi_2(\Sch,\syn)$.
By the Segal condition, we must have an equivalence 
\[
\fr(f,g)\simeq \fr(f)\times\fr(g)=\Maps_{K(Y)}(0,\sL_f)\times\Maps_{K(Z)}(0,\sL_g).
\]
The simplicial map $d_1\colon \Phi_2(\Sch,\syn)\to\Phi_1(\Sch,\syn)$ should then induce
\begin{equation}\label{eqn:composition-of-trivializations}
\Maps_{K(Y)}(0,\sL_f)\times\Maps_{K(Z)}(0,\sL_g) \to \Maps_{K(Z)}(0,\sL_{f\circ g}).
\end{equation}
In other words, given trivializations of $\sL_f$ and $\sL_g$ in $K$-theory, we must specify a trivialization of $\sL_{f\circ g}$ in $K$-theory. 
In this situation, we have the fundamental cofiber sequence
\[
g^*(\sL_f) \to \sL_{f\circ g} \to \sL_g
\]
in $\Perf(Z)$,
inducing a canonical equivalence $\phi\colon g^*(\sL_f)\oplus \sL_g\simeq \sL_{f\circ g}$ in $K(Z)$. The map~\eqref{eqn:composition-of-trivializations} is then the composition
\begin{multline*}
\Maps_{K(Y)}(0,\sL_f)\times\Maps_{K(Z)}(0,\sL_g) \xrightarrow{g^*\times\id} \Maps_{K(Z)}(0,g^*(\sL_f))\times\Maps_{K(Z)}(0,\sL_g)\\
\xrightarrow{\oplus}\Maps_{K(Z)}(0,g^*(\sL_f)\oplus\sL_g) \xrightarrow{\phi} \Maps_{K(Z)}(0,\sL_{f\circ g}).
\end{multline*}

\sssec{}
Let us recall the notion of $n$-gapped objects for $n \geq 0$. Let
\[
\Ar_n=\Fun(\Delta^1, \Delta^n)
\]
be the category of arrows in $\Delta^n$. Explicitly, this has objects $(i, j)$ for $0 \leq i \leq j \leq n$ and a unique morphism $(i, j) \rightarrow (i', j')$ whenever $i \leq i'$ and $j \leq j'$. This category is displayed as:

\begin{equation*}
\begin{tikzcd}
(0,0) \ar{r} & (0,1) \ar{r} \ar{d} & (0,2) \ar{r} \ar{d} & \cdots \ar{r} & (0,n) \ar{d} \\
 & (1,1) \ar{r}  & (1,2) \ar{r} \ar{d} & \cdots \ar{r} & (1,n) \ar{d} \\
 &  & (2,2) \ar{r}  & \cdots \ar{r} & (2,n) \ar{d} \\
 &  &   & \ddots & \vdots\ar{d} \\
 &  &    & & (n,n).
\end{tikzcd}
\end{equation*}

\sssec{}  \label{sssec:ngapusual} Suppose that $\C$ is a pointed $\infty$-category with finite colimits.

\begin{defn} The $\infty$-category of \emph{$n$-gapped objects} in $\C$, $\Gap(n, \C)$, is the full subcategory of $\Fun(\Ar_n, \C)$ spanned by those functors $F\colon \Ar_n \rightarrow \C$ such that:
\begin{enumerate}
\item $F(i,i)$ is a zero object for all $i$;
\item For $i \leq j \leq k$, the diagram
\begin{equation*}
\begin{tikzcd}
F(i,j) \ar{r} \ar{d} & F(i, k) \ar{d}\\
F(j, j)\ar{r} & F(j, k)
\end{tikzcd}
\end{equation*}
is coCartesian.
\end{enumerate}
\end{defn}
 
These conditions imply that the functor $\Gap(n, \C) \rightarrow \Fun(\Delta^{n-1}, \C)$, given by restriction along $\Delta^{n-1}\hookrightarrow \Ar_n$, $i\mapsto (0,i+1)$, is an equivalence (here, $\Delta^{-1}=\emptyset$). Its inverse is given by left Kan extension.

If $\InftyCat^{\mathrm{pt,rex}}$ denotes the $\infty$-category of pointed $\infty$-categories with finite colimits and functors that preserve finite colimits, we have a functor
\[
\Gap\colon\Delta^\op\times\InftyCat^{\mathrm{pt,rex}}\to \InftyCat, \quad ([n],\C)\mapsto \Gap(n,\C).
\]

\sssec{} \label{ngaprelative} 
We now introduce a relative version of the gap construction. Suppose that $p\colon X \rightarrow S$ is a coCartesian fibration such that:
\begin{itemize}
	\item for every $s\in S$, the fiber $X_s$ is a pointed $\infty$-category with finite colimits;
	\item for every $\alpha\colon t\to s$ in $S$, the functor $\alpha_*\colon X_t\to X_s$ preserves finite colimits.
\end{itemize} 
A morphism $f$ in $X$ is called \emph{vertical} if $p(f)$ is an equivalence. 
As a special case, we call a morphism in $\Ar_n$ \emph{vertical} if it is of the form $(i,j) \rightarrow (i',j)$.

\begin{defn} \label{def:ngaprelative} The $\infty$-category of \emph{relative $n$-gapped objects}, $\Gap_S(n,X)$, is the full subcategory of $\Fun(\Ar_n, X)$ spanned by functors $F$ such that:
\begin{enumerate}
\item $F$ takes vertical morphisms to vertical morphisms;
\item $F(i,i)$ is a zero object in its fiber for all $i$;
\item For $i \leq j \leq k$, the square
\begin{equation*}
\begin{tikzcd}
F(i,j) \ar{r} \ar{d} & F(i, k) \ar{d}\\
F(j, j) \ar{r} & F(j, k)
\end{tikzcd}
\end{equation*}
is a relative pushout square.
\end{enumerate}
\end{defn}

Given conditions (1) and (2), condition (3) states that the square
\[
\begin{tikzcd}
\alpha_*F(i,j) \ar{r} \ar{d} & F(i, k) \ar{d}\\
0 \ar{r} & F(j, k)
\end{tikzcd}
\]
is a pushout square in $X_{(p\circ F)(j,k)}$, where $\alpha=(p\circ F)((j,j)\to (j,k))$.

Just as in the usual setting, for any monotone map $q\colon [n] \rightarrow [m]$, the functor $q^*\colon \Fun(\Ar_m,X)\to \Fun(\Ar_n,X)$ preserves the subcategories of relative gapped objects.
Moreover, for any strict morphism of coCartesian fibrations 
\begin{equation*}
\begin{tikzcd}
 X \ar{rr}{f} \ar{dr} &[-15pt]  &[-15pt] X' \ar{dl} \\
  & S &
\end{tikzcd}
\end{equation*}
 that preserves finite colimits fiberwise, the induced functor $f_*\colon \Fun(\Ar_n,X)\to \Fun(\Ar_n,X')$ also preserves relative gapped objects.
We therefore obtain a functor
$$\Gap_S\colon\Delta^{\op} \times \coCart_S^{\mathrm{pt,rex}} \rightarrow \InftyCat,\quad ([n], X\to S) \mapsto \Gap_S(n, X),$$
where $\coCart_S^{\mathrm{pt,rex}}$ corresponds to $\Fun(S,\InftyCat^{\mathrm{pt,rex}})$ under the Grothendieck construction.

\sssec{} 
For the coCartesian fibration $\C \rightarrow *$, we recover the usual notion of an $n$-gapped object in $\C$ (\sssecref{sssec:ngapusual}). 
For the coCartesian fibration $\id_S\colon S \rightarrow S$, we have the following identification. 
We denote by $\iota_n$ the functor
\[
\iota_n\colon \Delta^{n}\hookrightarrow\Ar_n,\quad i\mapsto (0,i).
\]

\begin{lem} \label{lem:Gap-id}
	Let $S$ be an arbitrary $\infty$-category. 
	Then the composition
	 $$\Gap_S(n, S) \hookrightarrow \Fun(\Ar_n,S) \xrightarrow{\iota_n^*} \Fun(\Delta^n, S)$$
	 is an equivalence of $\infty$-categories.
\end{lem}

\begin{proof} 
	Since the fibers of $\id_S$ are contractible, $\Gap_S(n,S)$ is the full subcategory of $\Fun(\Ar_n,S)$ spanned by the functors sending vertical morphisms to equivalences. It is clear that the functor $\iota_n^*$ has a fully faithful left adjoint identifying $\Fun(\Delta^n,S)$ with this subcategory.
\end{proof}

\sssec{} Given a coCartesian fibration $p\colon X \rightarrow S$ in $\coCart_S^{\mathrm{pt,rex}}$, we can view $p$ as a strict morphism of coCartesian fibrations over $S$:
\begin{equation*}
\begin{tikzcd}
 X \ar{rr}{p} \ar[swap]{dr}{p} &[-15pt]  &[-15pt] S \ar{dl}{\id_S} \\
  & S\rlap. &
\end{tikzcd}
\end{equation*}
By Lemma~\ref{lem:Gap-id}, we obtain a functor 
 $$p_*\colon\Gap_S(n, X) \rightarrow \Gap_S(n, S) \stackrel{\iota_n^*}\simeq \Fun(\Delta^n, S).$$

\begin{prop} \label{prop:gapcocart}
	Suppose that $p\colon X \rightarrow S$ is a coCartesian fibration in $\coCart_S^{\mathrm{pt,rex}}$.
	
	\noindent{\em(i)}
	For every $n$, the functor $$p_*\colon \Gap_S(n, X) \rightarrow \Fun(\Delta^n, S)$$ is a coCartesian fibration.
	
	\noindent{\em(ii)}
	For every monotone map $\phi\colon [m]\to [n]$, the functor $$\phi^*\colon \Gap_S(n, X)\to\Gap_S(m,X)$$ preserves coCartesian edges.
\end{prop}

\begin{proof} 
	Since the map $p\colon X \rightarrow S$ is a coCartesian fibration, the functor $p_*\colon\Fun(\Ar_n, X) \rightarrow \Fun(\Ar_n, S)$ is a coCartesian fibration
	and the functor $\phi^*\colon \Fun(\Ar_n,X)\to\Fun(\Ar_m,X)$ preserves coCartesian edges \cite[Proposition 3.1.2.1]{HTT}. 
	Since $\Gap_S(n, X)$ and $\Gap_S(n, S)$ are full subcategories of $\Fun(\Ar_n, X)$ and $\Fun(\Ar_n, S)$ respectively, 
	it suffices to show that if $f\colon F \rightarrow G$ is a $p_*$-coCartesian arrow in $\Fun(\Ar_n, X)$ with $F$ being $n$-gapped and $p_*G$ being $n$-gapped, then $G$ is also $n$-gapped.

Firstly, since $p_*G$ is $n$-gapped, it takes vertical morphisms to equivalences, hence $G$ takes vertical morphisms to vertical morphisms. Secondly, the assumption that $f$ is $p_*$-coCartesian means that $G(i,j) \simeq \alpha_{ij*}F(i,j)$ where $\alpha_{ij}=p(f(i,j))$. Since $\alpha_{ij*}$ preserves finite colimits and $F$ is $n$-gapped, we immediately deduce that $G$ satisfies conditions (2) and (3) of Definition~\ref{def:ngaprelative}.
\end{proof}

\sssec{} Let $\coCart\subset \Fun(\Delta^1,\InftyCat)$ be the subcategory whose objects are the coCartesian fibrations $p\colon X\to S$ and whose morphisms are the squares
\[
\begin{tikzcd}
	X \ar{r}{f} \ar[swap]{d}{p} & X' \ar{d}{p'} \\ 
	S \ar{r} & S'
\end{tikzcd}
\]
such that $f$ sends $p$-coCartesian edges to $p'$-coCartesian edges. 

For example, by Proposition~\ref{prop:gapcocart}, if $p\colon X\to S$ is a coCartesian fibration in $\coCart_S^{\mathrm{pt,rex}}$, we have a functor
\[
\Delta^\op \to \coCart,\quad [n]\mapsto (p_*\colon \Gap_S(n,X) \to \Fun(\Delta^n,S)).
\]

\begin{lem}\label{lem:unstraightening}
	Let $C$ be an $\infty$-category and $F\colon C\to\Fun(\Delta^1,\InftyCat)$ a functor.
	Then $F$ lands in $\coCart$ if and only if the induced functor
	\[
	\int_C d_1\circ F \to \int_C d_0\circ F
	\]
	is a coCartesian fibration.
\end{lem}

\begin{proof}
	This follows from \cite[Lemma 1.4.14]{LurieGoodwillie}.
\end{proof}

\sssec{} 
We now specialize to the coCartesian fibration $p\colon\QCoh \rightarrow \Sch^{\op}$, whose fiber over a scheme $S$ is the stable $\infty$-category $\QCoh(S)$. 
By Proposition~\ref{prop:gapcocart}, we have a morphism of simplicial $\infty$-categories
 \begin{equation}\label{eqn:simp-cocart}
 p_*\colon \Gap_{\Sch^\op}(\bullet, \QCoh) \rightarrow \Gap_{\Sch^\op}(\bullet, \Sch^\op) \stackrel{\iota_\bullet^*}\simeq \Fun(\Delta^{\bullet}, \Sch^{\op}),
 \end{equation}
 which is a coCartesian fibration in each degree.

\sssec{}\label{sssec:Lbullet}
Recall that the functoriality of the cotangent complex is encoded by the triangle
\begin{equation*}
  \begin{tikzcd}
     & \QCoh \ar{d}{p} \\
    \Fun(\Delta^1,\Sch)^\op \ar{r}{s} \ar{ur}{\sL} &  \Sch^\op,
  \end{tikzcd}
\end{equation*}
where $s$ is the source functor. We define a simplicial functor
\[
\sL_\bullet\colon \Fun(\Delta^\bullet,\Sch^\op) \to \Fun(\Ar_\bullet,\QCoh)
\]
by the diagram
\[
\begin{tikzcd}
	 &[+1.5em] & \Fun(\Ar_n,\QCoh) \ar{d}{p} \\
	\Fun(\Delta^n,\Sch^\op) \ar{r}{\Fun(\Delta^1,-)} \ar[bend left=10]{urr}{\sL_n} \ar[swap,bend right=10]{drr}{\id} & \Fun(\Ar_n,\Fun(\Delta^1,\Sch)^\op) \ar{ur}{\sL} \ar{r}{s} & \Fun(\Ar_n,\Sch^\op) \ar{d}{\iota_n^*} \\
	& & \Fun(\Delta^n,\Sch^\op),
\end{tikzcd}
\]
where we identify $\Fun(\Delta^1,\Sch^\op)$ with $\Fun(\Delta^1,\Sch)^\op$ using the canonical isomorphism $\Delta^1\simeq(\Delta^1)^\op$. 
Here is the functor $\sL_2$:
\[
\begin{tikzcd}[column sep={3em,between origins},row sep={3em,between origins}]
X \ar{r}{f^\op} & Y \ar{r}{g^\op} & Z
\end{tikzcd}
\quad\mapsto\quad
\begin{tikzcd}[column sep={3em,between origins},row sep={3em,between origins}]
\id_X \ar{r} & f \ar{r} \ar{d}  & fg \ar{d} \\
 & \id_Y \ar{r} & g \ar{d} \\
 & & \id_Z
\end{tikzcd}
\quad\mapsto\quad
\begin{tikzcd}[column sep={3em,between origins},row sep={3em,between origins}]
0 \ar{r} & \sL_f \ar{r} \ar{d}  & \sL_{fg} \ar{d} \\
 & 0 \ar{r} & \sL_g \ar{d} \\
 & & 0\rlap.
\end{tikzcd}
\]
Since $g^*\sL_f \to \sL_{fg} \to \sL_g$ is a cofiber sequence in $\QCoh(Z)$, the last diagram is a relative $2$-gapped object. 
Since a functor $F\colon \Ar_n\to \QCoh$ is a relative $n$-gapped object if and only if $\phi^*(F)$ is a relative $2$-gapped object for every $\phi\colon [2]\to [n]$, 
this implies that $\sL_n$ lands in the full subcategory $\Gap_{\Sch^\op}(n,\QCoh)\subset \Fun(\Ar_n,\QCoh)$. In other words, $\sL_\bullet$ is a section of the simplicial coCartesian fibration~\eqref{eqn:simp-cocart}.

Recall that $\sL\colon \Fun(\Delta^1,\Sch)^\op\to\QCoh$ takes Tor-independent Cartesian squares to coCartesian edges.
Since the cotangent complex of an lci morphism is perfect and since flat morphisms are Tor-independent from any other morphism, we have
\begin{equation}\label{eqn:L-factorization}
\begin{tikzcd}
	\Fun(\Delta^\bullet,\Sch^\op) \ar{r}{\sL_\bullet} & \Gap_{\Sch^\op}(\bullet,\QCoh) \\
	\Phi_{\bullet^\op}(\Sch,\mathrm{syn})^\op \ar[hookrightarrow]{u} \ar[dashed]{r}{\sL_\bullet} & \Gap_{\Sch^\op}(\bullet,\Perf)^\cocart. \ar[hookrightarrow]{u}
\end{tikzcd}
\end{equation}
Here, $\Phi_{\bullet^\op}$ denotes the composition of $\Phi_\bullet$ with the canonical involution of $\Delta^\op$. This is needed because
\[
\Phi_{\bullet}(\Sch,\mathrm{syn})^\op \subset \Fun(\Delta^\bullet,\Sch)^\op\simeq \Fun((\Delta^\bullet)^\op,\Sch^\op).
\]

\sssec{} $K$-theory is taken in the sense of \cite{blumberg2013universal}. In particular, it is a functor $$K\colon \InftyCat^\mathrm{st} \rightarrow \Spc,$$
where $\InftyCat^\mathrm{st}$ is the $\infty$-category of stable $\infty$-categories, with a natural transformation
$$\eta\colon (-)^{\simeq} \rightarrow K.$$
The defining property of $K$-theory is \emph{additivity}: 
 \begin{enumerate}
 	\item $K(0)\simeq *$;
 	\item given a split exact sequence of stable $\infty$-categories
	 $$A \rightarrow B \rightarrow C,$$ 
	 the induced map $K(B) \to K(A) \times K(C)$ is an equivalence;
	 \item the structure of $\Einfty$-space on $K(-)$ induced by (1) and (2) is grouplike.
 \end{enumerate}
 If $X$ is a scheme, we set
  \[
  K(X) = K(\Perf(X)).
  \]

 \sssec{}
For a functor $\sigma\colon \Delta^n\to \Sch^\op$, denote by $\Gap_\sigma(\Perf)$ the fiber of
 \[
 p_*\colon \Gap_{\Sch^\op}(n,\Perf)\to \Fun(\Delta^n,\Sch^\op)
 \]
 over $\sigma$. Then $\Gap_\sigma(\Perf)$ is a stable $\infty$-category, and there is a split exact sequence of stable $\infty$-categories
 \[
 \Gap_{\sigma\circ\delta^n}(\Perf) \to \Gap_\sigma(\Perf) \to \Perf(\sigma(n))
 \]
 where $\delta^n\colon \Delta^{n-1}\hookrightarrow\Delta^n$ is the initial segment and the second map is evaluation at $(n-1,n)\in\Ar_n$.
 By additivity and induction, we obtain a canonical equivalence
 \begin{equation}\label{eqn:KGap}
 K(\Gap_\sigma(\Perf)) \simeq K(\sigma(1)) \times K(\sigma(2)) \times \cdots \times K(\sigma(n))
 \end{equation}
 induced by evaluation at $(i-1,i)\in\Ar_n$, for $1\leq i\leq n$.
 
 \sssec{}\label{sssec:intK}
 Applying Lemma~\ref{lem:unstraightening} to the simplicial coCartesian fibration
 \[
 \Gap_{\Sch^\op}(\bullet,\Perf)\to \Fun(\Delta^\bullet,\Sch^\op),
 \]
 we obtain a coCartesian fibration
 \begin{equation}\label{eqn:intGap}
 \int_{\Delta^\op}  \Gap_{\Sch^\op}(\bullet,\Perf)\to  \int_{\Delta^\op} \Fun(\Delta^\bullet,\Sch^\op),
 \end{equation}
 which classifies a functor to $\InftyCat^\mathrm{st}$.
 Whiskering this functor to the natural transformation $\eta\colon (-)^\simeq \to K$ and unstraightening, we obtain a strict morphism of coCartesian fibrations in $\infty$-groupoids
 \begin{equation}\label{eqn:GaptoK}
 \begin{tikzcd}
 	 \int\Gap_{\Sch^\op}(\bullet,\Perf)^\cocart \ar{rr}{\eta} \ar{dr} &[-4em] &[-4em]  \int K\Gap_{\Sch^\op}(\bullet,\Perf) \ar{dl} \\
	 & \int\Fun(\Delta^\bullet,\Sch^\op)\rlap. & 
 \end{tikzcd}
 \end{equation}
Its fiber over $\sigma\colon \Delta^n\to\Sch^\op$ is the map
\[
\eta\colon \Gap_\sigma(\Perf)^\simeq \to K(\Gap_\sigma(\Perf)).
\]

\sssec{}\label{sssec:vop}
Given a coCartesian fibration $p\colon X \rightarrow S$ classified by $F\colon S\to\InftyCat$, its \emph{vertical opposite} $p^\vop\colon X^\vop\to S$ is the coCartesian fibration classified by $F(-)^\op$. 
The coCartesian fibration in spaces $\Tw(\C)^\op\to \C^\op\times \C$ (see \cite[\sectsign 5.2.1]{HA}) is natural in $\C\in\InftyCat$, and whiskering it by $F$ yields a functor
\[
\Fun(S,\InftyCat) \to \Fun(S,\Fun(\Delta^1,\InftyCat)).
\]
Applying the Grothendieck construction and Lemma~\ref{lem:unstraightening}, we obtain a functor
\[
\coCart_S \to \coCart, \quad (X\to S) \mapsto (\Tw_S(X)^\vop\to X^\vop\times_SX).
\]
The \emph{fiberwise mapping space functor} 
$$\Maps\colon X^\vop \times_S X \rightarrow \Spc$$
is the functor classifying the coCartesian fibration $\Tw_S(X)^\vop\to X^\vop\times_SX$. By construction, for every $s\in S$, the composite
\[
X_s^\op \times X_s \to X^\vop \times_S X \xrightarrow{\Maps} \Spc
\]
is the usual mapping space functor. 

Moreover, if
\begin{equation*}
\begin{tikzcd}
 X \ar{rr}{f} \ar[swap]{dr} &[-15pt]  &[-15pt] Y \ar{dl} \\
  & S &
\end{tikzcd}
\end{equation*}
is a strict morphism of coCartesian fibrations over $S$,
we have by functoriality a commutative square
\[
\begin{tikzcd}
	\Tw_S(X)^\vop \ar{r} \ar{d} &[+10pt] \Tw_S(Y)^\vop \ar{d} \\
	X^\vop\times_SX \ar{r}{f^\vop\times f} & Y^\vop \times_SY,
\end{tikzcd}
\]
whence by straightening an induced natural transformation
\begin{equation*}
\begin{tikzcd}
 X^\vop\times_SX \ar{r}{f^\vop\times f} \ar[swap,bend right=10]{dr}{\Maps} &[+10pt] Y^\vop\times_SY \ar{d}{\Maps} \\
 \ar[ur,phantom,near end, "\Rightarrow" right] & \Spc.
\end{tikzcd}
\end{equation*}

 If $x\colon S\to X^\vop$ is a section of $p^\vop$, we denote by $\Maps(x,-)\colon X\to\Spc$ the composite
 \[
 X \simeq S\times_SX \xrightarrow{x\times\id} X^\vop\times_SX \xrightarrow{\Maps} \Spc.
 \]
 Given $f\colon X\to Y$ as above, we obtain a natural transformation
 \begin{equation}\label{eqn:Maps-functoriality}
 \Maps(x,-) \to \Maps(f^\vop\circ x,f(-))\colon X\to \Spc.
 \end{equation}
 
 \sssec{}\label{sssec:relative-maps}
Applying the construction of~\sssecref{sssec:vop} to the coCartesian fibration
 \[
 \int K\Gap_{\Sch^\op}(\bullet,\Perf) \to \int\Fun(\Delta^\bullet,\Sch^\op)
 \] 
from~\eqref{eqn:GaptoK}, we obtain a functor
\begin{equation}\label{eqn:Maps(0,-)}
\Maps(0,-)\colon  \int K\Gap_{\Sch^\op}(\bullet,\Perf) \to \Spc.
\end{equation}
Here, $0$ refers to the section induced by the zero section of the stable coCartesian fibration~\eqref{eqn:intGap}.
Note that the restriction to the fiber over $\sigma\colon \Delta^n\to\Sch^\op$ is the functor
\[
\Maps(0,-)\colon K(\Gap_\sigma(\Perf)) \to \Spc.
\]

\sssec{}\label{sssec:definition-of-fr}
We can now define the labeling functor $\fr$ on the pair $(\Sch,\syn)$. It is the composition
\begin{multline}\label{eqn:fr}
\fr\colon \int\Phi_{\bullet}(\Sch,\syn)^\op\simeq \int\Phi_{\bullet^\op}(\Sch,\syn)^\op \\
\xrightarrow{\int\sL_\bullet} \int\Gap_{\Sch^\op}(\bullet,\Perf)^\cocart
 \xrightarrow{\eta} \int K\Gap_{\Sch^\op}(\bullet,\Perf)\xrightarrow{\Maps(0,-)} \Spc.
\end{multline}
The first equivalence is the base change of the canonical involution of $\Delta^\op$, the second map is~\eqref{eqn:L-factorization}, the third is~\eqref{eqn:GaptoK}, and the last is~\eqref{eqn:Maps(0,-)}.

Explicitly, if $\sigma\colon\Delta^n \to \Sch^\op$ is a sequence of syntomic morphisms, then
\[
\fr(\sigma) = \Maps_{K(\Gap_\sigma(\Perf))}(0, \sL_n(\sigma)).
\]
 
\begin{prop}\label{prop:fr-is-Segal}
	The functor $\fr$ defined in~\eqref{eqn:fr} is a labeling functor on $(\Sch,\syn)$.
\end{prop}

\begin{proof}
	We must check that $\fr$ is a Segal presheaf on $\Phi_{\bullet^\op}(\Sch,\syn)$. 
	Let $\sigma\colon\Delta^n \to \Sch^\op$ be a sequence of syntomic morphisms
	\[
	X_0 \xleftarrow{f_1} X_{1} \xleftarrow{f_2} \cdots \xleftarrow{f_{n}} X_n.
	\]
	By~\eqref{eqn:KGap}, we have
	\begin{align*}
	\fr(\sigma) &= \Maps_{K(\Gap_\sigma(\Perf))}(0, \sL_n(\sigma)) \\
	&\simeq 
	\Maps_{K(X_1)}(0,\sL_n(\sigma)(0,1)) \times \cdots\times\Maps_{K(X_n)}(0,\sL_n(\sigma)(n-1,n))\\
	&\simeq
	\Maps_{K(X_1)}(0,\sL_{f_1}) \times\cdots\times\Maps_{K(X_{n})}(0,\sL_{f_{n}}),
	\end{align*}
	where the middle equivalence is induced by evaluation at $(i-1,i)\in\Ar_n$ for $1\leq i\leq n$.
	We similarly have equivalences
	\[
	\fr(\rho_i^*(\sigma)) \simeq \Maps_{K(X_i)}(0,\sL_{f_i})
	\]
	induced by evaluation at $(0,1)\in\Ar_1$.
	Putting these facts together, we deduce that the Segal map
	\[
	\fr(\sigma) \to \fr(\rho_1^*(\sigma)) \times\cdots\times \fr(\rho_n^*(\sigma))
	\]
	is an equivalence, as desired.
\end{proof}

\begin{rem}\label{rem:lci-flat}
	The same construction produces a labeling functor on $(\Sch,\mathrm{lci},\mathrm{flat})$. Indeed, in the square~\eqref{eqn:L-factorization}, we can replace the lower left corner by $\Phi_{\bullet^\op}(\Sch,\mathrm{lci},\mathrm{flat})^\op$.
\end{rem}

\sssec{}\label{sssec:framedcorr}
Here is our main definition.

\begin{defn} \label{defn:main-def}
	Let $S$ be a scheme.	The \emph{$\infty$-category of framed correspondences} over $S$, $\Span^{\fr}(\Sch_S)$, is the complete Segal space $\Span_{\bullet}^{\fr}(\Sch_S, \fsyn)^+$. 
	If $C\subset\Sch_S$ is a full subcategory, we denote by $\Span^\fr(C)$ the corresponding full subcategory of $\Span^\fr(\Sch_S)$.
\end{defn}

By \sssecref{sssec:complete-segal-space}, the objects of $\Span^\fr(\Sch_S)$ are $S$-schemes, and we have an equivalence
\[
\Maps_{\Span^\fr(\Sch_S)}(X,Y) \simeq \Corr^\fr_S(X,Y)
\]
natural in $(X,Y)\in \Sch_S^\op\times\Sch_S$.

\begin{rem}
By Remark~\ref{rmk:not-complete}, the Segal space $\Span^{\fr}_{\bullet}(\Sch_S,\fsyn)$ is not complete: any nontrivial loop in $K(X)$ defines an automorphism of $X$ that is not visible in $\Span^{\fr}_{0}(\Sch_S,\fsyn)\simeq \Sch_S^\simeq$.
\end{rem}

\sssec{} \label{sssec:etale}
 Consider the $\infty$-category $\Span^{\fet}(\Sch_S)$ of finite étale correspondences, given by the complete Segal space $\Span_{\bullet}(\Sch_S, \fet)$ where $\fet$ denotes the collection of finite étale morphisms. We want to construct a functor of $\infty$-categories:
$$\Span^{\fet}(\Sch_S) \to \Span^{\fr}(\Sch_S).$$

To do this, we consider a $K$-theory-free variant of~\eqref{eqn:fr}:
\begin{multline*}
\fr'\colon \int\Phi_{\bullet}(\Sch,\syn)^\op\simeq \int\Phi_{\bullet^\op}(\Sch,\syn)^\op \\
\xrightarrow{\int\sL_\bullet} \int\Gap_{\Sch^\op}(\bullet,\Perf)^\cocart
 \xrightarrow{\Maps(0,-)} \Spc.
\end{multline*}
If $\sigma\colon\Delta^n \to \Sch^\op$ is a sequence of syntomic morphisms, then
\[
\fr'(\sigma) = \Maps_{\Gap_\sigma(\Perf)^\simeq}(0, \sL_n(\sigma)).
\]
By the functoriality of fiberwise mapping spaces described in~\eqref{eqn:Maps-functoriality}, $\eta$ induces a natural transformation $\fr'\to\fr$.

Since $0$ is an initial object of $\Gap_\sigma(\Perf)$, $\fr'(\sigma)$ is contractible if $\sL_n(\sigma)\simeq 0$ and empty otherwise. But the cotangent complex of a syntomic morphism is zero if and only if that morphism is étale, so
\[
\fr'(\sigma)\simeq\begin{cases}
* & \text{if every edge of $\sigma$ is étale,}\\
\emptyset & \text{otherwise.}
\end{cases}
\]
In particular, $\fr'$ is a labeling functor on $(\Sch,\syn)$ and its restriction to $(\Sch,\et)$ is the contractible labeling functor. 
We therefore have morphisms in $\Lab\Trip$
\[
(\Sch_S,\et;*)\simeq (\Sch_S,\et;\fr') \to (\Sch_S,\syn;\fr') \to (\Sch_S,\syn;\fr).
\]
By \sssecref{sssec:Span-functoriality}, we obtain in particular a morphism of Segal spaces
\[
\Span_\bullet(\Sch_S,\fet) \to \Span_\bullet^\fr(\Sch_S,\fsyn).
\]
Applying Segal completion, we obtain the desired functor $\Span^\fet(\Sch_S)\to\Span^\fr(\Sch_S)$.

\sssec{}
Let $\Span^{\fsyn}(\Sch_S)$ be the $\infty$-category of finite syntomic correspondences, given by the complete Segal space $\Span_\bullet(\Sch_S,\fsyn)$. Then the morphism in $\Lab\Trip$
\[
(\Sch_S,\fsyn;\fr) \to (\Sch_S,\fsyn;*)
\]
induces by \sssecref{sssec:Span-functoriality} a morphism of Segal spaces
\[
\Span^\fr_\bullet(\Sch_S,\fsyn) \to \Span_\bullet(\Sch_S,\fsyn).
\]
Applying Segal completion, we obtain a functor
\[
\Span^\fr(\Sch_S) \to \Span^\fsyn(\Sch_S).
\]


\ssec{The symmetric monoidal structure}
\label{ssec:symmon}

In this section, we construct a symmetric monoidal structure on the $\infty$-category $\Span^{\fr}(\Sch_S)$ such that the tensor product of two $S$-schemes $X$ and $Y$ is given by the Cartesian product $X\times_S Y$.
We recall the definition of a commutative monoid:

\begin{defn}
	Let $\C$ be an $\infty$-category with finite products. A \emph{commutative monoid} or \emph{$\Einfty$-object} in $\C$ is a functor $X\colon \Fin_* \to \C$ such that for every $n \geq 0$ the Segal maps
$\{1, \ldots, n\}_+ \to \{i\}_+$ induce an equivalence
$$X(\{1, \ldots, n\}_+) \simeq \prod_{i=1}^n X(\{i\}_+) .$$
A \emph{symmetric monoidal $\infty$-category} is a commutative monoid in the $\infty$-category $\InftyCat$.
\end{defn}

\sssec{}\label{sssec:Span-monoidal}
We first discuss symmetric monoidal structures in the general setting of \ssecref{ssec:labeled-corr}. 
Note that the $\infty$-category $\Trip$ has limits that are created by the forgetful functor $\Trip\to\InftyCat$. A commutative monoid in $\Trip$ will be called a \emph{symmetric monoidal triple}. Concretely, a symmetric monoidal triple is a triple $(\C,M,N)$ together with a symmetric monoidal structure on $\C$ such that:
\begin{itemize}
	\item $M$ and $N$ are stable under tensor products;
	\item the tensor product $\otimes\colon \C\times\C\to\C$ preserves pullbacks of morphisms in $M\times M$ along morphisms in $N\times N$.
\end{itemize}
The second condition is automatic if the symmetric monoidal structure on $\C$ is Cartesian (and it rarely holds otherwise).

The functor $\Span_\bullet\colon \Trip \to \Fun(\Delta^\op,\Spc)$ preserves limits since it is corepresentable by $(\Sigma_\bullet,\Sigma_\bullet^L,\Sigma_\bullet^R)$, hence it preserves commutative monoids.
Thus, if $(\C,M,N)$ is a symmetric monoidal triple, then the complete Segal space $\Span_\bullet(\C,M,N)$ is a symmetric monoidal $\infty$-category.

\sssec{} We now introduce the notion of a symmetric monoidal structure on a labeling functor $F$, which will allow us to construct a symmetric monoidal structure on the $\infty$-category $\Span_\bullet^F(\C,M,N)^+$.

\begin{lem}\label{lem:monoidal-labeling}
	Let $A$ be an $\infty$-category and $f\colon A\to\Trip$ a functor. Then there is an equivalence of $\infty$-categories between:
	\begin{itemize}
		\item Lifts of $f$ to $\Lab\Trip$;
		\item Functors $$F\colon \int_{\Delta^\op\times A} \Phi_\bullet^\op\circ f \to \Spc $$ such that, for every $a\in A$, the restriction of $F$ to $\int_{\Delta^\op}\Phi_\bullet(f(a))^\op$ is a Segal presheaf. 
	\end{itemize}
\end{lem}

\begin{proof}
	Recall that $\Lab\Trip\to\Trip$ is the Cartesian fibration classified by
	\[
	\Trip^\op\to\InftyCat,\quad (\C,M,N)\mapsto \Seg\Pre(\Phi_\bullet(\C,M,N)).
	\]
	Let $E\to\Trip$ be the Cartesian fibration classified by
	\[
	\Trip^\op\to\InftyCat,\quad (\C,M,N)\mapsto \textstyle\Fun(\int\Phi_\bullet(\C,M,N)^\op,\Spc),
	\]
	so that $\Lab\Trip$ is a full subcategory of $E$.
	The $\infty$-category of lifts of $f$ to $E$ is equivalent to the $\infty$-category of lifts of
	$$ A^\op \to \InftyCat,\quad a\mapsto \textstyle\Fun(\int\Phi_\bullet(f(a))^\op,\Spc),$$
	to the universal Cartesian fibration.
	Since $\Fun(-,\Spc)$ transforms left-lax colimits into right-lax limits (see the proof of \cite[Lemma 16.16]{norms}), we have an equivalence of $\infty$-categories
	\[
	\Fun_\Trip(A,E) \simeq \textstyle\Fun(\int_{\Delta^\op\times A} \Phi_\bullet^\op\circ f,\Spc).
	\]
	It remains to observe that a lift of $f$ to $E$ lands in $\Lab\Trip$ if and only if the corresponding functor $\int_{\Delta^\op\times A} \Phi_\bullet^\op\circ f\to \Spc$ satisfies the given condition.
\end{proof}

In the definition of Segal presheaf (Definition~\ref{def:segal-presheaf}), we can replace $\Delta^\op$ by any category with ``Segal maps'', e.g., $\Fin_\pt$.

\begin{defn}
	Let $(\C,M,N)$ be a symmetric monoidal triple. A \emph{symmetric monoidal labeling functor} on $(\C,M,N)$ is a functor
	\[
	F\colon \int_{([n],I_+)\in \Delta^\op\times\Fin_\pt} \Phi_n(\C^I,M^I,N^I)^\op \to \Spc
	\]
	such that:
	\begin{enumerate}
		\item for every $I_+\in\Fin_\pt$, the restriction of $F$ to $\int_{[n]\in\Delta^\op} \Phi_n(\C^I,M^I,N^I)^\op$ is a Segal presheaf;
		\item the restriction of $F$ to $\int_{I_+\in\Fin_\pt}\Phi_1(\C^I,M^I,N^I)^\op$ is a Segal presheaf.
	\end{enumerate}
\end{defn}

\sssec{}
Let $(\C,M,N)$ be a symmetric monoidal triple.
By Lemma~\ref{lem:monoidal-labeling}, a symmetric monoidal labeling functor $F$ on $(\C,M,N)$ is the same thing as a lift of 
\[
\Fin_\pt\to\Trip, \quad I_+\mapsto (\C^I,M^I,N^I),
\]
 to $\Lab\Trip$ satisfying a certain Segal condition. Composing with the functor $\Lab\Trip\to\Fun(\Delta^\op,\Spc)$ constructed in \sssecref{sssec:Span-functoriality}, we obtain in particular a $\Fin_\pt$-object
\begin{equation}\label{eqn:SpanF-monoidal}
\Fin_\pt\to \Fun(\Delta^\op,\Spc),\quad I_+\mapsto \Span_\bullet^{F_I}(\C^I,M^I,N^I),
\end{equation}
where $F_I$ is the restriction of $F$ to $\int_{\Delta^\op} \Phi_\bullet(\C^I,M^I,N^I)^\op$.

\begin{prop}\label{prop:SpanF-monoidal}
	Let $(\C,M,N)$ be a symmetric monoidal triple and let $F$ be a symmetric monoidal labeling functor on $(\C,M,N)$. Then the functor~\eqref{eqn:SpanF-monoidal} is a commutative monoid in Segal spaces.
\end{prop}

\begin{proof}
	Each $\Span_\bullet^{F_I}(\C^I,M^I,N^I)$ is a Segal space by Theorem~\ref{thm:segal}. 
	We have a commutative square
	\[
	\begin{tikzcd}
	\Span^{F_I}_n(\C^I,M^I,N^I) \ar{r} \ar{d} & \prod_{i\in I}\Span^{F_{\{i\}}}_n(\C,M,N)  \ar{d} \\
	\Span_n(\C^I,M^I,N^I) \ar{r} & \prod_{i\in I}\Span_n(\C,M,N),
	\end{tikzcd}
	\]
	and we must show that the top arrow is an equivalence. Since the bottom one is, it suffices to check that this square is Cartesian. Let $\sigma=(\sigma_i)_{i\in I}\in \Span_n(\C^I,M^I,N^I)$. The vertical fiber of this square over $\sigma$ is the map
	\begin{equation}\label{eqn:monoidal-Segal-map}
	\lim (F_I)_\sigma \to \prod_{i\in I} \lim (F_{\{i\}})_{\sigma_i},
	\end{equation}
	where
	\[
	(F_I)_\sigma\to (F_{\{i\}})_{\sigma_i}\colon \int_{\Delta^\op}\Phi_\bullet(\Sigma_n,\Sigma_n^L,\Sigma_n^R)^\op\to\Spc.
	\]
	The assumption that $F\colon \int_{I_+\in\Fin_\pt}\Phi_1(\C^I,M^I,N^I)^\op\to\Spc$ is a Segal presheaf means that for every collection $(\tau_i\colon \Delta^1\to\C)_{i\in I}$ of morphisms in $M$, the map
	\[
	F_I((\tau_i)_{i\in I}) \to \prod_{i\in I} F_{\{i\}}(\tau_i)
	\]
	is an equivalence. Since $F_I$ and $F_{\{i\}}$ are Segal presheaves, this holds more generally for any collection $(\tau_i\colon \Delta^m\to\C)_{i\in I}$ of sequences of morphisms in $M$.
	 It follows that 
	$(F_I)_\sigma \simeq \prod_{i\in I} (F_{\{i\}})_{\sigma_i}$.
	 Hence, the map~\eqref{eqn:monoidal-Segal-map} is an equivalence, as desired.
\end{proof}

Recall from \cite[Section 14]{rezk-css} that the Segal completion $X^+$ of a Segal space $X$ has the explicit description
\[
X^+\simeq \colim_{n\in\Delta^\op} X^{E(n)},
\]
for some $E(\bullet)\colon \Delta\to\Fun(\Delta^\op,\Spc)$. It follows that Segal completion preserves finite products. Under the assumptions of Proposition~\ref{prop:SpanF-monoidal}, we deduce that the functor
\[
\Fin_\pt \to \Fun(\Delta^\op,\Spc),\quad I_+\mapsto \Span_\bullet^{F_I}(\C^I,M^I,N^I)^+,
\]
is a commutative monoid in complete Segal spaces, i.e., it is a symmetric monoidal structure on the $\infty$-category $\Span_\bullet^{F_*}(\C,M,N)^+$.

\sssec{}\label{sssec:Corrfr-monoidal}
For a scheme $S$, let $p_S\colon \QCoh_S\to \Sch_S^\op$ be the pullback of the coCartesian fibration $p\colon \QCoh\to \Sch^\op$.
The $\infty$-category $\QCoh_S$ has finite colimits that are preserved by $p_S$: if $(X_i,\sF_i)$ is a finite diagram in $\QCoh_S$, then
\[
\colim_{i} (X_i,\sF_i) = (\lim_i X_i, \colim_i\pi_i^*\sF_i),
\]
where $\pi_i\colon \lim_i X_i\to X_i$ is the canonical projection. 

Moreover, the cotangent complex functor $\sL\colon \Fun(\Delta^1,\Sch_S)^\op\to \QCoh_S$ preserves finite coproducts when restricted to the full subcategory $\Fun^\mathrm{flat}(\Delta^1,\Sch_S)^\op$ spanned by the flat morphisms. Indeed, $\sL_{\id_S}=(S,0)$ is the initial object of $\QCoh_S$ and the canonical map
\[
\pi_1^*\sL_f \oplus \pi_2^*\sL_g \to \sL_{f\times_S g}
\]
is an equivalence when $f$ and $g$ are flat. It follows that the diagram
\begin{equation*}
  \begin{tikzcd}
     & \QCoh_S \ar{d}{p_S} \\
    \Fun^\mathrm{flat}(\Delta^1,\Sch_S)^\op \ar{r}{s} \ar{ur}{\sL} &  \Sch_S^\op
  \end{tikzcd}
\end{equation*}
can be promoted to a diagram of symmetric monoidal functors for the coCartesian symmetric monoidal structures. 

Note that a functor $\Ar_n\to \QCoh_S^I$ is a relative $n$-gapped object over $(\Sch_S^\op)^I$ if and only if it is so componentwise.
Consequently, we obtain a simplicial symmetric monoidal functor
$$I_+\mapsto p_{S*}^I\colon \Gap_{(\Sch_S^\op)^I}(\bullet, \QCoh_S^I) \to \Fun(\Delta^\bullet,(\Sch_S^\op)^I).$$ 
Since syntomic morphisms in $\Sch_S$ are preserved by finite products \cite[Tag 01UI]{stacks} and perfect complexes are preserved by finite sums, we obtain as in \sssecref{sssec:Lbullet} a simplicial symmetric monoidal functor
\[
I_+\mapsto \sL_\bullet^I\colon  \Phi_{\bullet^\op}(\Sch_S^I,\syn^I)^\op \to \Gap_{(\Sch_S^\op)^I}(\bullet,\Perf_S^I)^\cocart,
\]
which is a section of the previous one.

\sssec{}
We can now repeat the constructions of \sssecref{sssec:intK}, \sssecref{sssec:relative-maps}, and \sssecref{sssec:definition-of-fr} with the simplicial coCartesian fibration 
 \[
 p_*\colon \Gap_{\Sch^\op}(\bullet,\Perf)\to \Fun(\Delta^\bullet,\Sch^\op)
 \]
 replaced by the simplicial coCartesian fibration
\[
\int_{I_+\in \Fin_*}p_{S*}^I\colon \int_{I_+\in\Fin_*} \Gap_{(\Sch_S^\op)^I}(\bullet, \Perf_S^I) \to \int_{I_+\in\Fin_*} \Fun(\Delta^\bullet,(\Sch_S^\op)^I)
\]
(see Lemma~\ref{lem:unstraightening}), whose fiber over $(\sigma_i\colon \Delta^n\to\Sch_S^\op)_{i\in I}$ is the stable $\infty$-category
\[
\prod_{i\in I}\Gap_{\sigma_i}(\Perf).
\]
We arrive at the functor
\begin{multline}\label{eqn:fr-monoidal}
\fr_S^\otimes\colon \int\Phi_{\bullet}(\Sch_S^?,\syn^?)^\op\simeq \int\Phi_{\bullet^\op}(\Sch_S^?,\syn^?)^\op \\
\xrightarrow{\int\sL_\bullet^?} \int\Gap_{(\Sch_S^\op)^?}(\bullet,\Perf_S^?)^\cocart
 \xrightarrow{\eta} \int K\Gap_{(\Sch_S^\op)^?}(\bullet,\Perf_S^?)\xrightarrow{\Maps(0,-)} \Spc,
\end{multline}
where the Grothendieck constructions are taken over $(\bullet,?_+)\in \Delta^\op\times\Fin_*$. Explicitly, for $(\sigma_i)_{i\in I}$ a finite collection of sequences of syntomic morphisms $\sigma_i\colon \Delta^n\to \Sch_S^\op$, we have
\[
\fr_S^\otimes((\sigma_i)_{i\in I}) \simeq \Maps_{K(\prod_{i\in I}\Gap_{\sigma_i}(\Perf))}(0,(\sL_n(\sigma_i))_{i\in I}).
\]

\begin{prop}\label{prop:fr-monoidal}
	Let $S$ be a scheme. The functor $\fr_S^\otimes$ defined in~\eqref{eqn:fr-monoidal} is a symmetric monoidal labeling functor on $(\Sch_S,\syn)$.
\end{prop}

\begin{proof}
	This follows immediately from Proposition~\ref{prop:fr-is-Segal} and the fact that $K\colon \InftyCat^\mathrm{st}\to \Spc$ preserves finite products.
\end{proof}

By Proposition~\ref{prop:SpanF-monoidal}, we therefore obtain a symmetric monoidal structure on the $\infty$-category $\Span^\fr(\Sch_S)$. By construction, given two spans
\begin{equation*}
\begin{tikzcd}
 & Z \ar[swap]{dl}{f} \ar{dr}{g} & \\
 X  & & Y
\end{tikzcd}
\qquad
\begin{tikzcd}
 & Z' \ar[swap]{dl}{f'} \ar{dr}{g'} & \\
 X'  & & Y'
\end{tikzcd}
\end{equation*}
with trivializations $\tau\in\Maps_{K(Z)}(0,\sL_f)$ and $\tau'\in\Maps_{K(Z')}(0,\sL_{f'})$, their tensor product is the span
\[
\begin{tikzcd}
 & Z\times_SZ' \ar[swap]{dl}{f\times_Sf'} \ar{dr}{g\times_Sg'} & \\
 X\times_SX'  & & Y\times_SY',
\end{tikzcd}
\]
with trivialization of $\sL_{f\times_Sf'}\simeq \pi_1^*\sL_{f}\oplus\pi_2^*\sL_{f'}$ given by the image of $(\tau,\tau')$ by the functor
\[
K(Z)\times K(Z') \to K(Z\times_SZ'), \quad (x,y)\mapsto \pi_1^*(x)+\pi_2^*(y).
\]

\begin{rem}
	The same construction produces a symmetric monoidal labeling functor on the symmetric monoidal triple $(\mathrm{Fl}_S,\mathrm{lci},\mathrm{flat})$, where $\mathrm{Fl}_S$ is the category of flat $S$-schemes (see Remark \ref{rem:lci-flat}).
\end{rem}

\sssec{}\label{sssec:Corrfr-monoidal-functors}
The variant $\fr'$ from \sssecref{sssec:etale} can similarly be promoted to a symmetric monoidal labeling functor $\fr_S^{\prime\otimes}$ on $(\Sch_S,\syn)$, together with a natural transformation $\fr_S^{\prime\otimes}\to \fr_S^{\otimes}$. 
The morphisms in $\Lab\Trip$
\[
(\Sch_S,\et;*)\simeq (\Sch_S,\et;\fr') \to (\Sch_S,\syn;\fr') \to (\Sch_S,\syn;\fr)
\]
are thus promoted to morphisms of $\Fin_\pt$-objects. As a result, we obtain a symmetric monoidal functor
\[
\Span^\fet(\Sch_S) \to \Span^\fr(\Sch_S),
\]
where the symmetric monoidal structure on $\Span^\fet(\Sch_S)$ is that described in \sssecref{sssec:Span-monoidal}.

The natural transformation $\fr_S^\otimes\to *$ also yields a symmetric monoidal functor
\[
\Span^\fr(\Sch_S)\to \Span^\fsyn(\Sch_S).
\]


\section{Applications} \label{sect:app}

In this section, we give some immediate applications of the recognition principle. In \ssecref{ssec:rep-results}, we prove that the presheaf of normally framed correspondences on a smooth $S$-scheme $X$ is representable by an ind-smooth scheme. As a result, we obtain a model for the zeroth space of the motivic sphere spectrum in terms of the group completion of a smooth Hilbert scheme of ``framed points'' in $\A^\infty$. In \ssecref{ssec:bar}, we discuss models for $\T$- and $\G$-delooping in $\H^{\fr}(k)$. In direct analogy with topology, we call these the $\T$- and \emph{$\G$-bar constructions}; the formulas are also given by $\Delta^{\op}$-colimits.

Lastly, in \ssecref{ssec:fr-to-cyc} and \ssecref{ssec:burn-to-cyc}, we define two functors --- one into and the other out of the $\infty$-category $\Span^{\fr}(\Sm_S)$. The first is a functor $$\cyc\colon \Span^{\fr}(\Sm_S) \rightarrow \Span^{\mathrm{cyc}}(\Sm_S),$$ where $S$ is Noetherian and the target is Voevodsky's category of finite correspondences. The second is a functor $$\fix\colon \Span(\Fin_G) \rightarrow \Span^{\fr}(\Sm_S),$$ where $S$ is connected with absolute Galois group $G$ and the domain is Barwick's effective Burnside $2$-category of finite $G$-sets.
If $S$ is the spectrum of a perfect field, these functors give rise to adjunctions
\[\cyc^*: \SH^\eff(S) \rightleftarrows \DM^\eff(S):\cyc_*
\]
and 
\[
\fix^*: \Spt_{G} \rightleftarrows \SH^\eff(S): \fix_*,
\]
which recover the well-known adjunctions.


\ssec{Representability of the motivic sphere spectrum} \label{ssec:rep-results}

\sssec{} 
Let $k$ be a perfect field and $X$ a smooth $k$-scheme. By Corollary~\ref{cor:main}, we have
\[
\Omega_\T^\infty\Sigma_\T^\infty X_+ \simeq L_\zar\Lhtp \h^\fr(X)^\gp \simeq L_\zar(\Lhtp \h^{\nfr}(X))^\gp\simeq  L_\zar(\Lhtp\colim_{n} \h^{\nfr,n}(X))^\gp.
\]
In this section, we will show that if $X$ is quasi-projective, then the presheaves $\h^{\nfr,n}(X)$ are representable by explicit smooth $k$-schemes. For example, $\h^{\nfr,n}(\Spec k)$ is a torsor under a smooth affine group scheme over the Hilbert scheme of zero-dimensional local complete intersections in $\A^n_k$.

\sssec{}
Let $S$ be a scheme and let $X\in\Sch_S$. We denote by $\Hilb^\fin(X/S)\colon\Sch_S^\op\to\Set$ the Hilbert functor of points of $X$ and by $\Hilb^\flci(X/S)\subset\Hilb^\fin(X/S)$ the subfunctor of local complete intersections. By definition, $\Hilb^\fin(X/S)(Y)$ (resp.\ $\Hilb^\flci(X/S)(Y)$) is the set of closed subschemes of $X\times_SY$ that are finite locally free (resp.\ finite syntomic) over $Y$. 
The degree of a finite locally free morphism induces coproduct decompositions
\[
\Hilb^\fin(X/S)\simeq\bigcoprod_{d\geq 0}\Hilb^\fin_d(X/S) \quad\text{and}\quad \Hilb^\flci(X/S)\simeq \bigcoprod_{d\geq 0}\Hilb^\flci_d(X/S)
\]
in the category of product-preserving presheaves on $\Sch_S$.
We recall some well-known representability properties of these Hilbert functors:

\begin{lem}\label{lem:Hilb-properties}
	Let $S$ be a scheme, let $X\in\Sch_S$, and let $d\geq 0$.
	
	\noindent{\em(i)}
	If $X\to S$ is separated, then $\Hilb^\fin(X/S)$ is representable by a separated algebraic space over $S$, which is locally of finite presentation if $X\to S$ is. It is a scheme if every finite set of points of every fiber of $X\to S$ is contained in an affine open subset of $X$, for example if $X\to S$ is locally quasi-projective.
	
	\noindent{\em(ii)}
	If $X\to S$ is finitely presented and locally (resp.\ strongly) quasi-projective, then $\Hilb^\fin_d(X/S)$ is finitely presented and locally (resp.\ strongly) quasi-projective over $S$.
	
	\noindent{\em(iii)}
	$\Hilb^\flci(X/S)\subset \Hilb^\fin(X/S)$ is an open subfunctor.
	
	\noindent{\em(iv)}
	If $X\to S$ is smooth, then $\Hilb^\flci(X/S)$ is formally smooth over $S$.
\end{lem}

\begin{proof}
	The first part of assertion (i) follows from \cite[Theorem 4.1]{RydhHilb} and \cite[Theorem 1.1]{OlssonStarr}, and the second part is \cite[Tag 0B9A]{stacks}.
	Assertion (ii) is \cite[Corollaries 2.7 and 2.8]{AltmanKleiman}.
	Assertion (iii) follows at once from \cite[Corollary 19.3.8]{EGA4-4}. 
	Let us prove (iv).
	Let $V\subset V'$ be a first-order thickening of affine schemes and let $V'\to S$ be a morphism. We must show that every closed subscheme $Z\subset V\times_SX$ that is finite syntomic over $V$ can be lifted a closed subscheme $Z'\subset V'\times_SX$ that is finite syntomic over $V'$. Since $X\to S$ is smooth, the immersion $i\colon Z\hookrightarrow V\times_SX$ is lci \cite[Tag 069M]{stacks}, and in particular the conormal sheaf $\sN_{Z/V\times_SX}$ is finite locally free.
Since $Z$ is affine, the canonical sequence of conormal sheaves
\[
0 \to i^*(\sN_{V\times_SX/V'\times_SX}) \to \sN_{Z/V'\times_SX} \to \sN_{Z/V\times_SX}\to 0
\]
is split exact. Choosing a splitting, we can identify $\sN_{Z/V\times_SX}$ with a subsheaf of $\sN_{Z/V'\times_SX}$.
Let $\tilde Z$ be the universal first-order infinitesimal neighborhood of $Z$ in $V'\times_SX$, so that $\sN_{Z/V'\times_SX}$ can be identified with the ideal sheaf of $Z$ in $\tilde Z$, and let $Z'\subset \tilde Z$ be the closed subscheme cut out by the subsheaf $\sN_{Z/V\times_SX}$. By construction, $\sN_{Z/Z'}$ is the pullback of $\sN_{V\times_SX/V'\times_SX}$ to $Z$, hence also the pullback of $\sN_{V/V'}$ to $Z$. By \cite[Tag 06BH]{stacks}, we deduce that $Z'$ is flat over $V'$ and that $Z=Z'\times_{V'}V$. By \cite[Tag 06AG(23)]{stacks}, we further deduce that $Z'$ is finite locally free over $V'$. Thus, $Z'$ defines a $V'$-point of $\Hilb^\fin(X/S)$ lifting $Z$. Finally, since $V\to V'$ is surjective, we immediately conclude by (iii) that $Z'\to V'$ is lci.
\end{proof}

\sssec{}
\label{sssec:representability}
Let $\Hilb^\flci(\A^n)$ be the Hilbert scheme of finite local complete intersections in $\A^n$ (over $\Z$). By Lemma~\ref{lem:Hilb-properties}, $\Hilb^\flci(\A^n)$ is a countable sum of smooth quasi-projective $\Z$-schemes. Let
	\[
\begin{tikzcd}
   \sZ \ar[swap]{dr}{p} \ar[hookrightarrow]{r}{i}
     & \A^n\times \Hilb^\flci(\A^n) \ar{d} \\
     & \Hilb^\flci(\A^n)
\end{tikzcd}
\]
be the universal finite local complete intersection in $\A^n$, and let $\Isom(\sN_i,\sO_\sZ^n)\to \sZ$ be the $\GL_n$-torsor classifying trivializations of the conormal sheaf $\sN_i$.

If $S$ is a scheme and $Y\in\Sch_S$, then $\h^{\nfr,n}_S(Y)\colon \Sch_S^\op\to\Set$ is the Weil restriction of $\Isom(\sN_i,\sO_\sZ^n)\times Y$ along $p\times\id_S\colon \sZ\times S\to \Hilb^\flci(\A^n)\times S$. Indeed, an $S$-morphism
\[
X\to (p\times\id_S)_*(\Isom(\sN_i,\sO_\sZ^n)\times Y)
\]
consists of a morphism $X\to \Hilb^\flci(\A^n)$ together with a $(\sZ\times S)$-morphism
\[
\sZ\times_{\Hilb^\flci(\A^n)} X  \to \Isom(\sN_i,\sO_\sZ^n)\times Y.
\]
By definition of the Hilbert scheme, this is precisely a closed subscheme $Z\subset\A^n_X$ that is finite syntomic over $X$, together with a trivialization of its conormal sheaf and an arbitrary $S$-morphism $Z\to Y$.

\begin{thm}
	Let $S$ be a scheme, $Y$ an $S$-scheme, and $n\geq 0$.
	
	\noindent{\em(i)}
	The functor $\h^{\nfr,n}_S(Y)\colon \Sch_S^\op\to\Set$ is representable by an algebraic space over $S$. It is a scheme if every finite set of points in every fiber of $Y\to S$ is contained in an affine open subset of $Y$, for example if $Y$ is locally quasi-projective over $S$.
	
	\noindent{\em(ii)}
	If $S$ is quasi-compact and $Y$ is quasi-projective over $S$, then $\h^{\nfr,n}_S(Y)$ is a countable sum of quasi-projective $S$-schemes.
	
	\noindent{\em(iii)}
	If $Y$ is smooth over $S$, so is $\h^{\nfr,n}_S(Y)$.
	
	\noindent{\em(iv)}
	For every morphism of schemes $S'\to S$, $\h^{\nfr,n}_{S'}(Y\times_SS') \simeq \h^{\nfr,n}_S(Y)\times_SS'$.
\end{thm}

\begin{proof}
 Since $p\colon\sZ\to\Hilb^\flci(\A^n)$ is finite, flat, and finitely presented, it follows from \cite[Theorem 3.7(b)]{RydhHilb} that $\h^{\nfr,n}_S(Y)$ is representable by an algebraic space over $S$.
The second part of (i) follows from \cite[\sectsign7.6, Theorem 4]{NeronModels}.
If $Y$ is quasi-projective over $S$, then $\h^{\nfr,n}_S(Y)$ is quasi-projective over $\Hilb^\flci(\A^n)\times S$ by \cite[Lemma 5.6]{norms}, hence it is a countable sum of quasi-projective $S$-schemes if $S$ is quasi-compact, by \cite[Proposition 5.3.4(ii)]{EGA2}.
	If $Y$ is smooth, then $\h^{\nfr,n}_S(Y)$ is smooth over $\Hilb^\flci(\A^n)\times S$ by \cite[\sectsign7.6, Proposition 5(h)]{NeronModels}, whence over $S$.
	Assertion (iv) is immediate from the definition of $\h^{\nfr,n}$.
\end{proof}

\sssec{} We discuss in more details the case of the sphere spectrum.

\begin{defn}\label{defn:Hilbfr}
	Let $X$ be a smooth $S$-scheme and $Z\in \Hilb^\fin(X/S)(T)$. A \emph{framing} of $Z$ is an isomorphism $\phi\colon \sN_{i}\simeq i^*(\Omega_{X\times_ST/T})$, where $i\colon Z\hook X\times_ST$. We define the \emph{Hilbert scheme of framed points} of $X/S$ to be the functor
	\[
	\Hilb^\fr(X/S)\colon \Sch_S^\op \to \Set
	\]
	sending $T$ to the set of pairs $(Z,\phi)$ with $Z\in \Hilb^\fin(X/S)(T)$ and $\phi$ a framing of $Z$.
\end{defn}

The degree induces a coproduct decomposition
\[
\Hilb^\fr(X/S) \simeq \bigcoprod_{d\geq 0} \Hilb^\fr_d(X/S)
\]
in the category of product-preserving presheaves on $\Sch_S$.

Observe that if $Z\in \Hilb^\fin(X/S)(T)$ admits a framing, then by Nakayama's lemma it is locally cut out by a minimal number of equations, so it is a local complete intersection over $T$. We therefore have a forgetful map
\[
\Hilb^\fr(X/S) \to \Hilb^{\flci}(X/S).
\]
As in \sssecref{sssec:representability}, $\Hilb^\fr(X/S)$ is then the Weil restriction of the torsor of framings of $\sZ\subset X\times_S\Hilb^\flci(X/S)$, hence it is a smooth $S$-scheme if $X$ is quasi-projective over $S$, by Lemma~\ref{lem:Hilb-properties}(ii,iv). If moreover $X$ has relative dimension $n$ over $S$, then $\Hilb^\fr_d(X/S)$ is smooth of relative dimension $dn(n+1)$ over $S$.

By definition, we have $\h^{\nfr,n}_S(S)\simeq \Hilb^\fr(\A^n_S/S)$. Define the ind-$S$-scheme $\Hilb^\fr(\A^\infty_S/S)$ by
\[
\Hilb^\fr(\A^\infty_S/S) = \colim_n\Hilb^\fr(\A^n_S/S),
\]
so that $\h^\nfr_S(S) \simeq \Hilb^\fr(\A^\infty_S/S)$. Proposition~\ref{prop:nfr-E_infinity} states that $\Lhtp \Hilb^\fr(\A^\infty_S/S)$ has a canonical $\Einfty$-structure, which can be described informally as forming disjoint unions of finite subschemes of $\A^\infty$. The second equivalence of Corollary~\ref{cor:main} becomes:

\begin{thm}\label{thm:sphere-representability}
	Suppose that $S$ is pro-smooth over a field. Then there is a canonical equivalence of grouplike $\Einfty$-spaces
	\[
	\Omega^\infty_\T \mathbf 1_S \simeq L_\zar(\Lhtp \Hilb^\fr(\A^\infty_S/S))^\gp.
	\]
\end{thm}

Theorem~\ref{thm:sphere-representability} can be regarded as an algebro-geometric analog of the description of the topological sphere spectrum in terms of framed $0$-dimensional manifolds and cobordisms.


\ssec{Motivic bar constructions} \label{ssec:bar}

In this section, we introduce various models for $\T$-deloopings and $\G$-deloopings in motivic homotopy theory, in other words, a ``motivic version'' of the bar construction. 

\sssec{}

Since the functor $\gamma^*\colon \H(S)_* \rightarrow \H^{\fr}(S)$ preserves colimits, we can write $\G^{\fr} \in \H^{\fr}(S)$ as the geometric realization of the simplicial diagram
\begin{equation} \label{eqn:bar-x}
\begin{tikzcd}[sep=small]
\gamma^*(\G_{m+}) 
& \gamma^*(\G_{m++})\arrow[l, shift left]
\arrow[l, shift right]
&  \gamma^*(\G_{m+++}) \arrow[l]
\arrow[l, shift left=
2
]
\arrow[l, shift right=
2
]
& \cdots
\arrow[l, shift left]
\arrow[l, shift right]
\arrow[l, shift left=3]
\arrow[l, shift right=3].
\end{tikzcd}
\end{equation}
Similarly, $\T^{\fr}$ is the geometric realization of the simplicial diagram
\begin{equation} \label{eqn:bar}
\begin{tikzcd}[sep=small]
\gamma^*\A^1_+ 
& \gamma^*((\A^1 \coprod \G_{m})_+) \arrow[l, shift left]
\arrow[l, shift right]
& \gamma^*((\A^1 \coprod \G_m \coprod \G_{m})_+) \arrow[l]
\arrow[l, shift left=
2
]
\arrow[l, shift right=
2
]
& \cdots
\arrow[l, shift left]
\arrow[l, shift right]
\arrow[l, shift left=3]
\arrow[l, shift right=3]
\end{tikzcd}
\end{equation}
Since colimits commutes with the monoidal structure in $\H^{\fr}(S)$, $\G^{\fr}\otimes \sF$ can be computed by tensoring the diagram~\eqref{eqn:bar-x} with $\sF$ levelwise and taking the colimit, and similarly for $\T^{\fr}\otimes\sF$.

\sssec{} If $\sF \in \H^{\fr}(S)$ we define the \emph{$\G$-bar construction} to be the diagram
\[\mathrm{Bar}_\G(\sF)\colon \Delta^{\op} \rightarrow \H(S)_*\] 
given by
\begin{equation*} \label{eqn:bar-x1}
\begin{tikzcd}[sep=small]
\gamma_*(\gamma^*(\G_{m+}) \otimes \sF)) 
& \gamma_*(\gamma^*(\G_{m++}) \otimes \sF)) \arrow[l, shift left]
\arrow[l, shift right]
&  \gamma_*(\gamma^*(\G_{m+++}) \otimes \sF)) \arrow[l]
\arrow[l, shift left=
2
]
\arrow[l, shift right=
2
]
& \cdots
\arrow[l, shift left]
\arrow[l, shift right]
\arrow[l, shift left=3]
\arrow[l, shift right=3].
\end{tikzcd}
\end{equation*}
Similarly, we define the \emph{$\T$-bar construction} to be the diagram 
\[\mathrm{Bar}_\T(\sF)\colon \Delta^{\op} \rightarrow \H(S)_*\] 
given by
\begin{equation*} \label{eqn:bar-x2}
\begin{tikzcd}[sep=small]
\gamma_*(\gamma^*(\A^1_+) \otimes \sF) 
& \gamma_*(\gamma^*((\A^1 \coprod \G_{m})_+) \otimes \sF) \arrow[l, shift left]
\arrow[l, shift right]
&  \gamma_*(\gamma^*((\A^1 \coprod \G_m \coprod \G_{m})_+) \otimes \sF) \arrow[l]
\arrow[l, shift left=
2
]
\arrow[l, shift right=
2
]
& \cdots
\arrow[l, shift left]
\arrow[l, shift right]
\arrow[l, shift left=3]
\arrow[l, shift right=3].
\end{tikzcd}
\end{equation*}
Since $\gamma_*$ preserves sifted colimits by Proposition~\ref{prop:gamma-basic}, the colimit of $\mathrm{Bar}_{\G}(\sF)$ (resp.\ $\mathrm{Bar}_{\T}(\sF)$) computes $\gamma_*(\G^{\fr} \otimes \sF)$ (resp.\ $\gamma_*(\T^{\fr} \otimes \sF)$).

\begin{rem}
	There is similarly a $\P$-bar construction $\mathrm{Bar}_\P(\sF)$, where $\P=(\P^1,\infty)$, such that $\lvert \mathrm{Bar}_\P(\sF)\rvert \simeq \lvert \mathrm{Bar}_\T(\sF)\rvert$.
\end{rem}

\sssec{}\label{sect:inv-to-loops} Although we can define the bar construction over an arbitrary base, we have to specialize to the hypotheses of Theorem~\ref{thm:main} in order to interpret its colimit as a delooping. We begin with the following consequence of Theorem~\ref{thm:main}.

\begin{prop} \label{prop:full-faith} Let $k$ be a perfect field. Then, for any $\sF \in \H^{\fr}(k)^{\gp}$, the unit maps 
\begin{equation*}
\sF \rightarrow \Omega_{\G}(\G^{\fr} \otimes \sF)
\end{equation*}
and
\begin{equation*}
\sF \rightarrow \Omega_{\T}(\T^{\fr} \otimes \sF)
\end{equation*}
are equivalences.
\end{prop}

\begin{proof} 
	This follows immediately from Theorem~\ref{thm:main}(i).
\end{proof}

By Proposition~\ref{prop:full-faith}, if $\sF \in \H^{\fr}(k)^{\gp}$, we obtain equivalences
\begin{align*}
\gamma_*\sF &  \simeq \gamma_*\Omega_{\G}(\G^{\fr} \otimes \sF) \\
& \simeq  \Omega_{\G}\gamma_*(\G^{\fr} \otimes \sF) \\
& \simeq  \Omega_{\G}\lvert\mathrm{Bar}_{\G}(\sF)\rvert.
\end{align*}
In a similar manner,
\[\gamma_*\sF \simeq \Omega_{\T}\lvert \mathrm{Bar}_\T(\sF)\rvert .\]


\ssec{From framed correspondences to Voevodsky's finite correspondences} \label{ssec:fr-to-cyc}

\nc{\sft}{\mrm{ft}}
\nc{\flf}{\mrm{flf}}

Let $S$ be a Noetherian scheme and let $\Sch_S^\sft\subset \Sch_S$ be the full subcategory of \emph{separated} $S$-schemes of finite type.
We denote by $\Span^\cyc(\Sch_S^\sft)$ Voevodsky's category of finite correspondences over $S$, as defined in \cite[\sectsign 9]{CD}. 
The goal of this subsection is to construct a symmetric monoidal functor
\[
\cyc\colon \Span^\fr(\Sch_S^\sft) \to \Span^\cyc(\Sch_S^\sft).
\]

\sssec{}
Recall that the objects of $\Span^\cyc(\Sch_S^\sft)$ are the separated $S$-schemes of finite type and that the set of morphisms from $X$ to $Y$ is the abelian group
\[
c_0(X\times_SY/X,\Z)
\]
of cycles on $X\times_SY$ that are finite and universally integral over $X$.

\sssec{}
Let $\Span^\flf(\Sch_S^\sft)=\Span(\Sch_S^\sft,\flf)$ be the $2$-category of finite locally free correspondences between $S$-schemes of finite type.
Suppose that
\begin{equation*}
  \begin{tikzcd}
    & Z \ar[swap]{ld}{f}\ar{rd}{g}
  \\
  X & & Y
  \end{tikzcd}
\end{equation*}
is a span over $S$ with $f$ finite flat, and let $(f,g)\colon Z\to X\times_SY$.
Since every irreducible component of $Z$ is flat and equidimensional over $X$, the fundamental cycle $[Z]$ on $Z$ is a Hilbert cycle over $X$ in the sense of \cite[Definition 8.1.10]{CD}. In particular, by \cite[8.1.35(P3)]{CD}, it is universally integral over $X$. It follows from \cite[Corollary 8.2.10]{CD} that the cycle $(f,g)_*[Z]$ on $X\times_SY$ is universally integral over $X$. Since it is also finite over $X$, we have
 \[
 (f,g)_*[Z]\in c_0(X\times_SY/X,\Z).
 \]
It is obvious that the cycle $(f,g)_*[Z]$ depends only on the isomorphism class of $Z$ over $X$ and $Y$. Hence, we obtain a map
\[
\cyc_{X,Y}\colon \Corr^\flf_S(X,Y) \to c_0(X\times_SY/X,\Z),\quad (X\xleftarrow fZ\xrightarrow gY) \mapsto (f,g)_*[Z].
\]

\begin{lem} \label{lem:cyc is functorial}
Let $S$ be a Noetherian scheme.

\noindent{\em(i)}
Let $X,Y,Z\in\Sch_S^\sft$ and let $X\xleftarrow f T\xrightarrow g Y$ and $Y\xleftarrow h W\xrightarrow k Z$ be finite flat correspondences. Then
\[
\cyc_{Y,Z}(W)\circ \cyc_{X,Y}(T) = \cyc_{X,Z}(T\times_{Y}W)
\]
in $c_0(X\times_SZ/X,\Z)$.

\noindent{\em(ii)}
Let $X,Y,X',Y'\in\Sch_S^\sft$ and let $X\xleftarrow f Z\xrightarrow g Y$ and $X'\xleftarrow{f'} Z'\xrightarrow{g'} Y'$ be finite flat correspondences. Then
\[
\cyc_{X,Y}(Z)\times_S\cyc_{X',Y'}(Z')=\cyc_{X\times_SX',Y\times_SY'}(Z\times_SZ')
\]
in $c_0(X\times_SX'\times_SY\times_SY'/X\times_SX',\Z)$.
\end{lem}

\begin{proof}
	(i) By definition of the composition of finite correspondences \cite[Definition 9.1.5]{CD}, we have
	\[
	\cyc_{Y,Z}(W)\circ \cyc_{X,Y}(T) = \pr_{XZ*}((f,g)_*[T] \otimes_Y (h,k)_*[W]),
	\]
	where $\pr_{XZ}\colon X\times_SY\times_SZ\to X\times_SZ$. 
	We have
	\begin{align*}
		(f,g)_*[T] \otimes_Y (h,k)_*[W] &= (\id_X\times_S(h,k))_*((f,g)_*[T]\otimes_Y[W]) \\
		&=(\id_X\times_S(h,k))_*((f,g)\times_Y\id_W)_*([T]\otimes_Y[W]) \\
		&=((f,g)\times_Y(h,k))_*([T]\otimes_Y[W]) \\
		&= ((f,g)\times_Y(h,k))_*[T\times_YW].
	\end{align*}
	The first and second equalities use \cite[Corollary 8.2.10]{CD},
	the third equality is the functoriality of the pushforward of cycles, and the last equality is \cite[8.1.35(P3)]{CD}, using that $[W]$ is a Hilbert cycle over $Y$. Hence,
	\[
	\cyc_{Y,Z}(W)\circ \cyc_{X,Y}(T) = \pr_{XZ*}((f,g)\times_Y(h,k))_*[T\times_YW] = \cyc_{X,Z}(T\times_YW),
	\]
	as desired.
	
	(ii) Let $\alpha=\cyc_{X,Y}(Z)$ and $\alpha'=\cyc_{X',Y'}(Z')$, and let us denote the product of $S$-schemes by concatenation. By definition \cite[Definition 9.2.1]{CD}, we have
	\[
	\alpha\times_S\alpha' = (\alpha\otimes_X[XX']) \otimes_{XX'} ([XX']\otimes_{X'}\alpha').
	\]
	By \cite[Corollary 8.2.10]{CD}, we have
	\begin{gather*}
	\alpha\otimes_X[XX'] = ((f,g)\times_X\id_{XX'})_*[ZX'], \\
	[XX']\otimes_{X'}\alpha' = (\id_{XX'}\times_{X'}(f',g'))_*[XZ'].
	\end{gather*}
	Using the same steps as in the proof of (i), we then have
	\begin{align*}
		\alpha\times_S\alpha' &=((f,g)\times_X\id_{XX'})_*[ZX'] \otimes_{XX'} (\id_{XX'}\times_{X'}(f',g'))_*[XZ'] \\
		&= (\id_{XX'Y}\times_{X'}(f',g'))_* (((f,g)\times_X\id_{XX'})_*[ZX'] \otimes_{XX'} [XZ']) \\
		&=(\id_{XX'Y}\times_{X'}(f',g'))_* ((f,g)\times_S\id_{Z'})_* ([ZX'] \otimes_{XX'} [XZ']) \\
		&=(ff',gg')_* ([ZX'] \otimes_{XX'} [XZ']) \\
		&= (ff',gg')_* [ZX'\times_{XX'} XZ'] = (ff',gg')_* [ZZ'],
	\end{align*}
	as desired.
\end{proof}
 
 \sssec{}\label{sssec:cyc-flf}
 It follows easily from Lemma~\ref{lem:cyc is functorial} that the maps $\cyc_{X,Y}$ assemble into a symmetric monoidal functor
 \[
\cyc\colon \Span^\flf(\Sch_S^\sft)\to \Span^\cyc(\Sch_S^\sft),
 \]
 such that the composite
 \[
 \Sch_{S+}^\sft \simeq \Span^\clopen(\Sch_S^\sft) \hook \Span^\flf(\Sch_S^\sft)\xrightarrow{\cyc} \Span^\cyc(\Sch_S^\sft)
 \]
 is the usual graph functor $\Gamma$ \cite[(9.1.8.1)]{CD}.
 
 \sssec{} \label{sec:frames-to-dm}
 We define the symmetric monoidal functor
\[
\cyc\colon \Span^\fr(\Sch_S^\sft) \to \Span^\cyc(\Sch_S^\sft)
\]
as the composition of the forgetful functor $\Span^\fr(\Sch_S^\sft)\to\Span^\flf(\Sch_S^\sft)$ from \sssecref{sssec:Corrfr-monoidal-functors} and the functor $\cyc$ from \sssecref{sssec:cyc-flf}.
Restricting $\cyc$ to smooth $S$-schemes and applying the usual constructions, we obtain an adjunction
\[
\cyc^*: \SH^\fr(S) \rightleftarrows \DM(S) : \cyc_*,
\]
where $\cyc^*$ is symmetric monoidal.
 
\begin{prop}\label{prop:commute}
	Let $S$ be a Noetherian scheme. Then the following triangle of symmetric monoidal functors commutes:
	\[
   \begin{tikzcd}
     \SH(S) \ar[swap]{d}{\gamma^*} \ar{rd}{\Gamma^*}
       & 
     \\
     \SH^\fr(S) \ar[swap]{r}{\cyc^*}
       & \DM(S).
   \end{tikzcd}
	\]
\end{prop}

\begin{proof}
	This follows from the commutativity of the triangle
	\[
	\begin{tikzpicture}[commutative diagrams/every diagram]
	\matrix[matrix of math nodes, name=m, commutative diagrams/every cell] {
	\Sm_{S+}^\sft & \\ \Span^\fr(\Sm_S^\sft) & \Span^\cyc(\Sm_S^\sft). \\};
	  \path[commutative diagrams/.cd, every arrow, every label]
	    (m-1-1) edge node[swap]{$\gamma$} (m-2-1)
	            edge node{$\Gamma$} (m-2-2)
	    (m-2-1) edge node[swap]{$\cyc$} (m-2-2);
 		\node [anchor=east,overlay,inner sep=0, outer sep=0] at (m-2-1 -| 0.5\textwidth,0) {\qedhere};
	\end{tikzpicture}
	\]
\end{proof}

When $k$ is a perfect field, it follows from Theorem~\ref{thm:main}(ii) that the symmetric monoidal functor
\[
\SH^\eff(k) \to \DM^\eff(k)
\]
can be identified with the left Kan extension of $\cyc\colon \Span^\fr(\Sm_k^\sft) \to \Span^\cyc(\Sm_k^\sft)$ to the $\infty$-categories of $\A^1$-invariant Nisnevich-local presheaves of spectra.

\subsection{From the Burnside category to framed correspondences} \label{ssec:burn-to-cyc} 
Let $S$ be a connected scheme with a separably closed geometric point $x\colon \Spec k\to S$, and let $G=\hat\pi_1^\et(S,x)$ be the profinite étale fundamental group of $S$ at $x$. Let $\Fin_G$ be the category of finite $G$-sets and $\Span(\Fin_G)$ Barwick's effective Burnside $2$-category of finite $G$-sets \cite[Section 3]{BarwickMackey}.  Here, we construct a symmetric monoidal functor 
\[
\fix\colon \Span(\Fin_G) \rightarrow \Span^{\fr}(\Sm_S),
\]
which sends a finite quotient of $G$ to the corresponding finite étale cover of $S$.

\sssec{} 
By Grothendieck's Galois theory of schemes, there is a canonical equivalence of categories
\[
\Fin_G \simeq \FEt_S,
\]
where $\FEt_S\subset \Sch_S$ is the full subcategory of finite étale $S$-schemes (see for example \cite[Tag 0BND]{stacks}).

We then define $\fix\colon \Span(\Fin_G) \rightarrow \Span^{\fr}(\Sm_S)$ as the composite
\[
\Span(\Fin_G) \simeq \Span(\FEt_S)\subset \Span^\fet(\Sm_S) \to \Span^\fr(\Sm_S),
\]
where the last functor is the symmetric monoidal functor defined in \sssecref{sssec:Corrfr-monoidal-functors}.

\sssec{}\label{sssec:SptG}
The symmetric monoidal $\infty$-category of $G$-spectra, for $G$ a profinite group, is defined in \cite[Example B]{BarwickMackey} as
\[
\Spt_G= \Pre_\Sigma(\Span(\Fin_G),\Spt).
\]
We will also need the alternative definition from \cite[\sectsign 9.2]{norms}, which characterizes $\Spt_G$ by the following universal property: if $\C$ is a symmetric monoidal $\infty$-category with sifted colimits such that $\otimes\colon\C\times\C\to\C$ preserves sifted colimits, then the functor $\Fin_{G+}\to \Spt_G$ induces a fully faithful functor
\[
\Fun^{\otimes,\mathrm{sift}}(\Spt_G,\C) \hook \Fun^\otimes(\Fin_{G+},\C),
\]
whose essential image consists of the functors whose extension to $\Pre_\Sigma(\Fin_G)_*$ inverts all $G$-representation spheres.

\sssec{} The functor $\fix$ induces a symmetric monoidal functor
\[
\fix^*\colon \Pre_\Sigma(\Span(\Fin_G)) \to \H^\fr(S)
\]
by left Kan extension.
Since $\fix^*$ sends $S^1$ to $S^1$, we have an induced symmetric monoidal functor
\[
\fix^*\colon \Spt_G \to \SH^{S^1,\fr}(S)
\]
by the symmetric monoidal universal property of stabilization \cite[Remark 2.25]{Robalo}. Finally, one may compose with $\Sigma^\infty_{\G,\fr}$ to obtain a symmetric monoidal functor
\[
\fix^*\colon \Spt_G \to \SH^{\fr}(S).
\]
Note that each version of $\fix^*$ preserves colimits and hence admits a right adjoint $\fix_*$.

\sssec{}
Recall from \cite[\sectsign 10]{norms} that we have a symmetric monoidal functor
\[
c\colon \Spt_G\to \SH(S),
\]
which is the unique symmetric monoidal colimit-preserving extension of
\[
\Fin_{G+} \simeq \FEt_{S+} \to \H(S)_* \xrightarrow{\Sigma^\infty_\T} \SH(S).
\]

\begin{prop}
	Let $S$ be a connected scheme with a separably closed geometric point $x$ and let $G=\hat\pi_1^\et(S,x)$. Then the following triangle of symmetric monoidal functors commutes:
	\[
	\begin{tikzcd}
		\Spt_G \ar{r}{c} \ar[swap]{dr}{\fix^*} & \SH(S) \ar{d}{\gamma^*} \\
		& \SH^\fr(S).
	\end{tikzcd}
	\]
\end{prop}

\begin{proof}
	By construction, both composites agree as symmetric monoidal functors when restricted to $\Fin_{G+}$, hence they agree by the universal property of $\Fin_{G+}\to\Spt_G$ recalled in~\sssecref{sssec:SptG}.
\end{proof}

When $k$ is a perfect field with separable closure $k^\mathrm{sep}$ and $G=\mathrm{Gal}(k^\mathrm{sep}/k)$, it follows from Theorem~\ref{thm:main}(ii) that the symmetric monoidal functor
\[
c\colon \Spt_G \to \SH^\eff(k)
\]
can be identified with the left Kan extension of $\fix\colon \Span(\Fin_G)\to\Span^\fr(\Sm_k)$.

\begin{rem} 
	The above discussion applies also to disconnected schemes $S$ if we replace $G=\hat\pi_1^\et(S,x)$ with the profinite \'{e}tale fundamental groupoid $\widehat\Pi_1^{\et}(S)$, see \cite[\sectsign 9–10]{norms}. 
\end{rem}


\appendix
\section{Voevodsky's Lemma}
\label{app:voevodsky-lemma}


\ssec{Voevodsky's Lemma}
In this section we compute the Nisnevich sheaf associated with the quotient of a scheme by an open subscheme. 
This computation first appeared in \cite[\sectsign 3]{garkusha2014framed}, where it is attributed to Voevodsky.
The result explains the definition of equationally framed correspondences (Definition~\ref{defn:Voevodsky framed correspondence}).

\sssec{}
Let $X$ be a scheme and $Z\subset X$ a closed subset. 
An \emph{étale neighborhood of $Z$ in $X$} is an algebraic space $U$ with an étale morphism $U \rightarrow X$ such that the projection $U \times_X Z \rightarrow Z$ is an isomorphism, for some (hence for any) closed subscheme structure on $Z$. 
We denote by $X^h_Z$ the pro-algebraic space of étale neighborhoods of $Z$ in $X$. 

\begin{lem}\label{lem:henselization}
	Let $X$ be a qcqs scheme and $Z\subset X$ a closed subset.
	
	\noindent{\em(i)}
	$X^h_Z$ is a limit of qcqs algebraic spaces.
	
	\noindent{\em(ii)}
	If $X$ is affine, then $X^h_Z$ is a limit of affine schemes.
	
	\noindent{\em(iii)}
	If $(X,Z)$ is a henselian pair, then $X^h_Z=X$.
	
	\noindent{\em(iv)}
	If $X$ is geometrically unibranch and has finitely many irreducible components, then, for every separated $X$-scheme $T$, every $X$-morphism $X^h_Z\to T$ factors uniquely through the pro-scheme of separated quasi-compact étale neighborhoods of $Z$.
\end{lem}

\begin{proof}
	(i) is a consequence of the fact that the category of étale algebraic spaces over $X$ is compactly generated, with compact objects the qcqs étale algebraic spaces.
	(ii) follows from (iii) by considering the usual henselization of $X$ along $Z$.
	(iii) follows from \cite[Proposition A.7]{Rydh:2010}.
	 Let us prove (iv).
	Let $V$ be a qcqs étale neighborhood of $Z$, and let $R\rightrightarrows U$ be a presentation of $V$ by finitely presented étale $X$-schemes.
	We will show that there exists a factorization $R\hook \bar R\hook U\times_XU$ where $\bar R$ is a closed étale equivalence relation on $U$ and $R\hook\bar R$ is a dense open immersion. This result and (i) complete the proof as every $V\to T$ factors uniquely through $U/\bar R$, which is a separated quasi-compact étale neighborhood of $Z$ (hence a scheme).
	By nil-invariance of the étale site \cite[Proposition A.4]{Rydh:2010}, we can assume that $X$ is reduced. Then every quasi-compact étale $X$-scheme has integral local rings and finitely many irreducible components. Let $\bar R$ be the scheme-theoretic closure of $R$ in $U\times_XU$. 
	This has the desired properties by \cite[Propositions A.5.5 and B.2.1]{ferrand-kahn}.
\end{proof}

\sssec{}
\label{sssec:Q}
Let $X$ be an $S$-scheme and let $U\subset X$ be an open subscheme.
Let $Q(X,U)$ be the pointed presheaf on $S$-schemes whose sections over $Y$ are pairs $(Z,\varphi)$ where $Z\subset Y$ is a closed subset and $\varphi\colon Y^h_Z \to X$ is a morphism of pro-algebraic spaces over $S$ such that $\varphi^{-1}(X-U)=Z$. The base point is the unique pair with $Z=\initial$. There is a factorization
\[
X/U\to Q(X,U) \to L_\nis(X/U)
\]
given as follows. The first map sends $f\colon Y\to X$ to 
\[(f^{-1}(X-U), Y^h_{f^{-1}(X-U)} \to Y \to X).\]
The second map sends $(Z\subset Y,\varphi\colon Y^h_Z\to X)$ to
\[Y\simeq L_\nis(Y^h_Z\coprod_{Y^h_Z-Z}(Y-Z))\xrightarrow{\varphi} L_\nis(X/U),\]
using the fact that the inclusion of schemes into algebraic spaces is cocontinuous for the Nisnevich topology. By inspection, this construction gives a commutative diagram
\begin{equation*}\label{eqn:Q(X,U)}
	\begin{tikzcd}
		X/U \ar{r} \ar{d} & Q(X,U) \ar{d} \ar{dl} \\
		L_\nis(X/U) \ar{r} & L_\nis Q(X,U)
	\end{tikzcd}
\end{equation*}
of pointed presheaves on $S$-schemes.

\begin{prop}\label{prop:VoevodskyLemma}
	Let $X$ be an $S$-scheme and let $U\subset X$ be an open subscheme.
	Then the maps $Q(X,U)\to L_\nis(X/U)\to L_\et(X/U)$ are isomorphisms.
\end{prop}

\begin{proof}
		By the above commutative diagram, it suffices to show that $Q(X,U)$ is an étale sheaf.
		It is obvious that $Q(X,U)$ transforms coproducts into products. 
		Let $V\to Y$ be a surjective étale morphism. We must show that
		\[
		Q(X,U)(Y) \to Q(X,U)(V) \rightrightarrows Q(X,U)(V\times_YV)
		\]
		is an equalizer diagram. Since closed subsets satisfy fpqc descent, it suffices to show that, given $Z\subset Y$ with preimage $W\subset V$, the diagram
		\[
		(V\times_YV)^h_{W\times_ZW} \rightrightarrows V^h_W \to Y^h_Z,
		\]
		is a coequalizer of pro-algebraic spaces over $S$.
		Note that $V^h_W \to Y^h_Z$ is an epimorphism since it is a limit of surjective étale morphisms.
		Suppose given a commutative diagram
		\[
		(V\times_YV)^h_{W\times_ZW} \rightrightarrows V^h_W \to C
		\]
		where $C$ is an algebraic space over $S$, represented by a commutative diagram $A \rightrightarrows B \to C$ where $A$ and $B$ are étale algebraic spaces over $Y$. Recall that the category of étale algebraic spaces over $Y$ is cocomplete (being equivalent to the étale topos of $Y$) and that the inclusion of étale algebraic spaces over $Y$ into algebraic spaces over $S$ preserves colimits (since it has a right adjoint).
		The coequalizer of $A\rightrightarrows B$ is then an étale neighborhood of $Z$ in $Y$, which shows that the morphism $V^h_W\to C$ factors through $Y^h_Z$.
\end{proof}

\begin{cor}\label{cor:VoevodskyGeneral}
	Let $S$ be a scheme and let $\nis\leq \tau\leq\et$.
	Let $X$ and $Y$ be $S$-schemes, $U\subset X$ an open subscheme, and $W\subset Y$ a sieve. 
	Then
	\[
	\Maps_{\Pre(\Sch_{S})_\pt}(Y/W,L_\tau(X/U))
	\]
	is the pointed set of pairs $(Z,\varphi)$, where $Z\subset Y$ is a closed subset such that $W\subset Y-Z$ and $\varphi\colon Y^h_Z\to X$ is a morphism of pro-algebraic spaces over $S$ such that $\varphi^{-1}(X-U)=Z$.
\end{cor}

\begin{proof}
	By Proposition~\ref{prop:VoevodskyLemma}, $L_\tau(X/U)\simeq Q(X,U)$, giving us the desired description.
\end{proof}

\begin{rem}\label{rem:stupid}
If $Z'\subset Y$ is a closed subscheme such that $Z'\subset Y-Z$, then $Y^h_Z\simeq (Y-Z')^h_Z$.
\end{rem}

\begin{cor}[Voevodsky]\label{cor:VoevodskyLemma}
	Let $S$ be a scheme, let $X,Y\in\Sch_{S}$, and let $\nis\leq \tau\leq\et$.
	For every $n \geq 0$, there is a natural bijection of pointed sets 
  \begin{equation*}
    \Corr^{\efr,n}_S(X,Y) \simeq \Maps_{\Pre(\Sch_{S})_\pt}(X_+ \wedge (\P^1)^{\wedge n}, L_\tau(Y_+ \wedge \T^{\wedge n})),
  \end{equation*}
  where $\P^1$ is pointed at $\infty$ and $\T=\A^1/\A^1-0$.
\end{cor}

\begin{proof}
Given Remark~\ref{rem:stupid}, this is a special case of Corollary~\ref{cor:VoevodskyGeneral}.
\end{proof}

\ssec{Closed gluing}
In the category of schemes, every cospan $A\hookleftarrow C\hookrightarrow B$ of closed immersions admits a pushout $A\coprod_CB$, which can be computed in the category of ringed topological spaces \cite[Théorème 7.1]{ferrand}. In particular, the canonical maps $A\to A\coprod_CB$ and $B\to A\coprod_CB$ are closed immersions.

\begin{defn}\label{defn:closedgluing}
Let $\sF$ be a presheaf on $\Sch_{S}$. We say that $\sF$ \emph{satisfies closed gluing} if it sends every pushout square of closed immersions to a pullback square and if $\sF(\emptyset)\simeq *$.
\end{defn}

\begin{rem}
	Closed gluing is weaker than closed descent and is not a topological condition. A presheaf on $\Sch_{S}$ is a sheaf for the topology of finite closed covers if and only if it satisfies closed gluing and is nil-invariant.
	For every $X\in\Sch_{S}$, the presheaf $Y\mapsto \Maps_S(Y,X)$ satisfies closed gluing, but it is nil-invariant if and only if $X\to S$ is formally étale.
\end{rem}

\begin{lem}\label{lem:lifting-neighborhoods}
	Let $X$ be a scheme, $i\colon X_0\hook X$ a closed immersion, $Z\subset X$ a closed subset, and $Z_0=i^{-1}(Z)$. Then the pullback functor
	\[
	i^*\colon\{\text{étale neighborhoods of $Z$ in $X$}\} \to \{\text{étale neighborhoods of $Z_0$ in $X_0$}\}
	\]
	has a fully faithful right adjoint.
\end{lem}

\begin{proof}
	Identifying étale algebraic spaces with sheaves of sets on the small étale site, the right adjoint is given by the pushforward functor $i_*$.
\end{proof}

\begin{prop}\label{prop:closedgluing}
	Let $X$ be an $S$-scheme and $U\subset X$ an open subscheme. Then $Q(X,U)$ satisfies closed gluing.
\end{prop}

\begin{proof}
	It is obvious that $Q(X,U)(\emptyset)\simeq *$.
	Let $Y_{01}\hook Y_0$ and $Y_{01}\hook Y_1$ be closed immersions with pushout $Y=Y_0\coprod_{Y_{01}} Y_1$.
	Note that the presheaf sending $Y$ to its set of closed subsets satisfies closed gluing.
	Let $Z\subset Y$ be a closed subset and set $Z_*=Y_*\cap Z\subset Y_*$ for $*\in\{0,1,01\}$. It remains to show that
	\[
   \begin{tikzcd}
    (Y_{01})^h_{Z_{01}} \ar{r}\ar{d}
       & (Y_1)^h_{Z_1} \ar{d}
     \\
     (Y_0)^h_{Z_0} \ar{r}
       & Y^h_Z
   \end{tikzcd}
	\]
	is a pushout square of pro-algebraic spaces. This follows immediately from Lemma~\ref{lem:lifting-neighborhoods}, since left adjoint functors are coinitial and pushouts of pro-objects can be computed levelwise.
\end{proof}

\sssec{}
Let $K$ be a finite nonsingular simplicial set. Then $K$ admits a geometric realization $\lvert K\rvert_S$ in the category of $S$-schemes with respect to the cosimplicial object $\A^\bullet_S$, which can be built inductively using pushouts of closed immersions. For example, $\lvert\partial\Delta^n\rvert_S=\partial\A^n_S$. 

Our interest in closed gluing stems from the following tautological lemma:

\begin{lem}\label{lem:sset-as-schemes}
  Let $\sF$ be a presheaf on $\Sch_S$ satisfying closed gluing and let $K$ be a finite nonsingular simplicial set. Then the canonical map of presheaves
  \[
  \sF(-\times_S\lvert K\rvert_S) \to \Maps(K,\sF(-\times\bA^\bullet))
  \]
  is an equivalence, where the right-hand side is a mapping space in $\Fun(\Delta^\op,\Spc)$.
\end{lem}

\section{Finite fields}
\label{app:finite-fields}


In this appendix, we develop some techniques that can be used to extend results about infinite perfect fields to finite fields.
Similar results were obtained independently by A.\ Druzhinin and J.\ Kylling \cite[Section 3]{DruzhininKylling}.

\ssec{Framed correspondences induced by monic polynomials}

\sssec{}
Let $S$ be a scheme and $p\in \sO(S)[x]$ a monic polynomial cutting out a closed subscheme $T\subset \A^1_S$. Then the map $f\colon T\to S$ is finite syntomic and the $\sO_T$-module $(p)/(p^2)$ is free of rank $1$ generated by $p$. We denote by $\phi(p)$ the tangentially framed correspondence
\begin{equation*}
  \begin{tikzcd}
     & T \ar[swap]{ld}{f,\alpha}\ar{rd}{f} & \\
    S &   & S
  \end{tikzcd}
\end{equation*}
over $S$, where $\alpha$ is the trivialization of $\sL_{f}$ in $K(T)$ induced by $p$. Note that $\phi(p)\simeq\id_S$ if $p$ has degree $1$.

\sssec{}
Recall that there is a canonical action of $\Omega K(S)$ on $S$ in the $\infty$-category $\Span^\fr(\Sch_S)$. If $a$ is a unit in $\sO(S)$, let us write $\langle a\rangle$ for the image of $a$ by the canonical map $\sO(S)^\times\to\Omega K(S)$. 
For $n\geq 0$, we define $n_\epsilon$ to be the formal alternating sum
\[
n_\epsilon = 1 + \langle -1\rangle + 1 + \dotsb
\]
with $n$ terms.

\sssec{}
We shall say that the a full subcategory $C\subset \Sch_S$ is \emph{admissible} if it contains $S$, if it is closed under finite coproducts, and if $X\in C$ implies $\A^1\times X\in C$.

\begin{prop}\label{prop:monic-polynomial}
	Let $S$ be a scheme, $C\subset\Sch_S$ an admissible subcategory, and $p\in \sO(S)[x]$ a monic polynomial of degree $d$.
	Then, for every $\sF\in \Pre_{\Sigma,\A^1}(\Span^\fr(C))$, the map $\phi(p)^*\colon \sF(S)\to \sF(S)$ is multiplication by $d_\epsilon$.
\end{prop}

\begin{proof}
	Write $p=x^d+q$. If $V\subset \A^1\times\A^1_S$ is the vanishing locus of $x^d+tq$, then $V$ is finite syntomic over $\A^1\times S$. Let $\beta$ be the trivialization of $\sL_{V/\A^1\times S}$ induced by $x^d+tq$. Then $V$ and $\beta$ define a framed correspondence $\psi\colon \A^1\times S\to S$ such that $\psi\circ i_1\simeq\phi(p)$ and $\psi\circ i_0\simeq\phi(x^d)$. Since $\sF$ is $\A^1$-invariant, we have $\phi(p)^*\simeq \phi(x^d)^*\colon\sF(S)\to \sF(S)$. Thus, we may assume that $p=x^d$, and we prove the result by induction on $d$. If $d=0$, the result is trivial. By considering the vanishing locus of $(x-t)x^{d-1}$ in $\A^1\times \A^1_S$, we deduce as before that $\phi(x^d)^*\simeq \psi^* + \chi^*\colon\sF(S)\to \sF(S)$, where $\psi$ and $\chi$ are the framed correspondences
	\[
    \begin{tikzcd}
       & S[x]/(x-1) \ar{ld} \ar{rd} & \\
      S &   & S
    \end{tikzcd}
	\qquad
    \begin{tikzcd}
       & S[x]/(x^{d-1}) \ar{ld} \ar{rd} & \\
      S &   & S
    \end{tikzcd}
	\]
	with trivializations of the cotangent complex induced by the generator $(x-1)x^{d-1}$ in the respective conormal sheaves $(x-1)/(x-1)^2$ and $(x^{d-1})/(x^{d-1})^2$. In the first conormal sheaf, this generator is equal to $x-1$, so $\psi=\phi(x-1)\simeq\id_S$. In the second conormal sheaf, this generator differs from $x^{d-1}$ by the unit $x-1$ in $\sO(S)[x]/(x^{d-1})$. Observe that the units $x-1$ and $-1$ are $\A^1$-homotopic via $xt-1\in \sO(S)[t,x]/(x^{d-1})^\times$. The framed correspondence
	\[
   \begin{tikzcd}
      & \A^1\times S[x]/(x^{d-1}) \ar{ld} \ar{rd} & \\
     \A^1\times S &   & S
   \end{tikzcd}
	\]
	with cotangent complex trivialized by $(xt-1)x^{d-1}$ is thus an $\A^1$-homotopy from $\phi(x^{d-1}) \circ \langle -1\rangle$ to $\chi$.
	Hence, we have $\phi(x^d)^*\simeq \id + \langle -1\rangle \phi(x^{d-1})^{ *}\colon\sF(S)\to \sF(S)$, and we conclude using the induction hypothesis and the recursive formula $d_\epsilon = 1 + \langle -1\rangle (d-1)_\epsilon$.
\end{proof}

\ssec{Conservativity results}

\begin{prop}\label{prop:coprime-injective}
	Let $S$ be a scheme, $C\subset\Sch_S$ an admissible subcategory, and $p_1,p_2\in \sO(S)[x]$ monic polynomials of coprime degrees whose vanishing loci $T_1$ and $T_2$ belong to $C$. Then, for every $\A^1$-invariant additive presheaf of abelian groups on $\Span^\fr(C)$, the pullback map $\sF(S)\to \sF(T_1)\times \sF(T_2)$ is injective.
\end{prop}

\begin{proof}
	Let $d_i$ be the degree of $p_i$. By Proposition~\ref{prop:monic-polynomial}, $\phi(p_i)^*\colon \sF(S)\to \sF(S)$ is multiplication by $(d_i)_\epsilon$. 
	Let $x$ be in the kernel of $\sF(S)\to \sF(T_1)\times \sF(T_2)$. Since the map $\phi(p_i)^*$ factors through the pullback $\sF(S)\to \sF(T_i)$, we have $(d_i)_\epsilon x=0$. By Bézout's identity, we can find integers $u_1,u_2\in\Z$ such that $u_1d_1+u_2d_2=1$. Then $u_1(d_1)_\epsilon + u_2 (d_2)_\epsilon = 1 + n(1-\langle -1\rangle)$ for some $n\in\Z$, hence $x=n(\langle -1\rangle-1)x$. Multiplying through by $1+\langle -1\rangle$, we find $(1+\langle -1\rangle)x=0$, that is, $\langle -1\rangle x = -x$. If $d_i$ is odd (and at least one of $d_1$ and $d_2$ is), this implies $(d_i)_\epsilon x = x$, so $x=0$.
\end{proof}

\begin{rem}\label{rem:finite-efr}
	If $S$ is affine and $T_1$ and $T_2$ are unramified over $S$, the conclusion of Proposition~\ref{prop:coprime-injective} also holds for $\sF$ an $\A^1$-invariant stable radditive presheaf of abelian groups on $\Span^\efr_*(C)$ (see \sssecref{sssec:hitr}). Indeed, we can replace $C$ by the smallest admissible subcategory of $\Sch_S$ containing $T_1$ and $T_2$, and in this case it follows from Corollaries~\ref{cor:vfr-vs-nfr} and \ref{cor:nfr-vs-dfr} that such presheaves factor uniquely through $\Span^\fr(C)$.
\end{rem}

\begin{prop}\label{prop:coprime-conservative}
	Let $S$ be a scheme and $p_1,p_2\in \sO(S)[x]$ monic polynomials of coprime degrees whose vanishing loci $T_1$ and $T_2$ are smooth over $S$. Then the following pullback functors are conservative and detect homotopic maps:
	\begin{enumerate}
		\item $\Pre_{\Sigma,\A^1}(\Span^\fr(\Sm_S))^\gp \to \Pre_{\Sigma,\A^1}(\Span^\fr(\Sm_{T_1}))^\gp \times \Pre_{\Sigma,\A^1}(\Span^\fr(\Sm_{T_2}))^\gp$,
		\item $\H^\fr(S)^\gp \to \H^\fr(T_1)^\gp \times \H^\fr(T_2)^\gp$,
		\item $\SH^{S^1,\fr}(S) \to \SH^{S^1,\fr}(T_1)\times \SH^{S^1,\fr}(T_2)$,
		\item $\SH^\fr(S) \to \SH^\fr(T_1)\times \SH^\fr(T_2)$,
		\item $\SH(S)\to \SH(T_1)\times \SH(T_2)$.
	\end{enumerate}
\end{prop}

\begin{proof}
	Since $\A^1$-invariance is preserved by delooping, the $\infty$-category $\Pre_{\Sigma,\A^1}(\Span^\fr(\Sm_S))^\gp$ is prestable.
	In particular, a morphism is an equivalence if and only if the identity map of its cofiber is nullhomotopic.
	To prove (1), it therefore suffices to show that that the given functor detects homotopic maps, and this follows from Proposition~\ref{prop:coprime-injective} applied to the presheaves $\pi_0\Hom(\sF,\sF')$ for $\sF,\sF'\in\Pre_{\Sigma,\A^1}(\Span^\fr(\Sm_S))^\gp$.
	Since $T_1$ and $T_2$ are smooth over $S$, the pullback functor in (1) preserves Nisnevich-local objects, hence (2) follows directly from (1).
	The proofs of (3) and (4) are identical to that of (1), as is the proof of (5) once we know that for $E,E'\in\SH(S)$, the presheaf $\pi_0\Omega^\infty_\T\Hom(E,E')$ on $\Sm_S$ extends to $\Span^\fr(\Sm_S)$. This is indeed the case by the theory of fundamental classes developed in \cite{DJKFundamental}.
\end{proof}

\begin{rem}
	If $S$ is the spectrum of a field, one can also prove Proposition~\ref{prop:coprime-conservative}(5) using Morel's identification of the endomorphisms of the motivic sphere spectrum with the Grothendieck–Witt group and \cite[Proposition 5.2]{HoyoisGLV}.
\end{rem}

\begin{cor}\label{cor:infinite-conservative}
	Let $k$ be a field, $p\neq q$ prime numbers, and $k_p$ (resp.\ $k_q$) a separable algebraic $p$-extension (resp.\ $q$-extension) of $k$. Then the following pullback functors are conservative:
	\begin{enumerate}
		\item $\Pre_{\Sigma,\A^1}(\Span^\fr(\Sm_k))^\gp \to \Pre_{\Sigma,\A^1}(\Span^\fr(\Sm_{k_p}))^\gp\times\Pre_{\Sigma,\A^1}(\Span^\fr(\Sm_{k_q}))^\gp$,
		\item $\H^\fr(k)^\gp \to \H^\fr(k_p)^\gp \times \H^\fr(k_q)^\gp$,
		\item $\SH^{S^1,\fr}(k) \to \SH^{S^1,\fr}(k_p)\times \SH^{S^1,\fr}(k_q)$,
		\item $\SH^\fr(k) \to \SH^\fr(k_p)\times \SH^\fr(k_q)$,
		\item $\SH(k)\to \SH(k_p)\times \SH(k_q)$.
	\end{enumerate}
\end{cor}

\begin{proof}
	Since $k_p$ and $k_q$ are filtered colimits of smooth $k$-algebras, the pullback functor in (1) preserves Nisnevich-local objects, so that (2) follows from (1).
	We give the proof of (5) in such a way that the proofs of (1), (3), and (4) are obtained by change of notation.
	Since $\SH(k)$ is prestable, it suffices to show that if $E\in \SH(k)$ becomes zero in $\SH(k_p)\times \SH(k_q)$, then $E\simeq 0$. Since moreover $\SH(k)$ is compactly generated, it suffices to show that if $C\in\SH(k)$ is compact and $n\geq 0$, any map $c\colon C\to \Sigma^n E$ is nullhomotopic. Write $k_p=\colim_\alpha k_\alpha$ where each $k_\alpha$ is a finite $p$-extension of $k$. Then, by compactness of $C$,
	\[
	\Maps_{\SH(k_p)}(C_{k_p},\Sigma^nE_{k_p}) \simeq\colim_\alpha \Maps_{\SH(k_\alpha)}(C_{k_\alpha}, \Sigma^n E_{k_\alpha}),
	\]
	so $c$ becomes nullhomotopic over some finite $p$-extension of $k$, and similarly for $q$. We then conclude using Proposition~\ref{prop:coprime-conservative}.
\end{proof}

\begin{rem}[Bachmann]
	If $k$ is a field, the pullback functors $\SH^{S^1}(k)\to \SH^{S^1}(k(t))$ and hence $\SH(k)\to \SH(k(t))$ are conservative. Indeed, for every $E\in\SH^{S^1}(k)$ and $n\in \Z$, the sheaf $\pi_{n}^\nis(E)$ is $\A^1$-invariant and unramified by \cite[Theorem 6.1 and Corollary 6.9]{Morel}, and so for any smooth $k$-scheme $X$, we have $\pi_{n}^\nis(E)(X)\simeq \pi_{n}^\nis(E)(X\times\A^1)\subset \pi_{n}^\nis(E)(X_{k(t)})$.
\end{rem}


\let\mathbb=\mathbf

{\small
\newcommand{\etalchar}[1]{$^{#1}$}
\providecommand{\bysame}{\leavevmode\hbox to3em{\hrulefill}\thinspace}

}

\parskip 0pt

\end{document}